\documentclass[nonblindrev]{informs3_noinforms2c}
\RequirePackage{tgtermes}
\RequirePackage{newtxtext}
\RequirePackage{newtxmath}
\RequirePackage{bm}
\RequirePackage{endnotes}

\OneAndAHalfSpacedXII % Current default line spacing
\usepackage{algpseudocode}
\usepackage{tikz}
\usepackage{tabularray}
\usepackage{amsfonts,amscd,epsfig}
\usepackage{mathtools}
\usepackage[colorlinks = true,
            linkcolor = blue,
            urlcolor  = blue,
            citecolor = blue,
            anchorcolor = blue]{hyperref}
\usepackage{graphicx}
\usepackage{subfig}
\usepackage{natbib}
\usepackage{caption}
\usepackage[capposition=top]{floatrow}
\usepackage{tabularx}

\usepackage[ruled]{algorithm2e}

\usepackage{bm}
\usepackage{diagbox}
\usepackage{booktabs}
\usepackage{multirow}
\usepackage{ltxtable}
\usepackage{booktabs,caption}
\usepackage{flafter}
\usepackage{float}
\usepackage{bbm}
\usepackage[flushleft]{threeparttable}
\usepackage{url}
\usepackage{accents}
\usepackage{lscape}
\usepackage{changes}
\usepackage{makecell}

\makeatletter
\g@addto@macro{\UrlBreaks}{\UrlOrds}
\makeatother
\usepackage{natbib}
\DeclareMathOperator*{\Prob}{\mathbb{P}}

\newcommand{\indep}{\perp\!\!\!\!\perp} 
\newcommand{\R}{\mathbb{R}}
\newcommand{\E}{\mathbb{E}}

\newcommand{\Gcal}{\mathcal{G}}

\newcommand{\Mcal}{\mathcal{M}}

\newcommand{\Xcal}{\mathcal{X}}

\newcommand{\Var}{\mathrm{Var}}

\newcommand{\gap}{\vspace{0.1in}}

\newcommand{\onebld}{{\bf 1}}
\newcommand{\wt}{\widetilde}
\newcommand{\wh}{\widehat}

\newcommand{\ma}[1]

\allowdisplaybreaks

\bibpunct[, ]{(}{)}{,}{a}{}{,}%
\def\bibfont{\small}%
\EquationsNumberedThrough    % Default: (1), (2), ...
\TheoremsNumberedThrough     % Preferred (Theorem 1, Lemma 1, Theorem 2)
\ECRepeatTheorems  %  

\begin{document}
%%%%%%%%%%%%%%%%

%\setcounter{page}{1}
%\DOI{10.1287/xxxx.0000.0000}%

\RUNAUTHOR{}

\RUNTITLE{Weight-Clipping-Based Policy Learning}
\TITLE{Offline Policy Learning with Weight Clipping and Heaviside Composite Optimization}
\ARTICLEAUTHORS{
       \AUTHOR{Jingren Liu}
		\AFF{Institute of Operations Research and Analytics, National University of Singapore, Singapore 117602, \EMAIL{jingren.liu@u.nus.edu}}
        \AUTHOR{Hanzhang Qin}
        \AFF{
        Department of Industrial Systems Engineering and Management, National University of Singapore, Singapore 117576,
        \EMAIL{hzqin@nus.edu.sg}}
        \AUTHOR{Junyi Liu}
		\AFF{
       Department of Industrial Engineering, Tsinghua University, Beijing, China 100084,
        \EMAIL{junyiliu@tsinghua.edu.cn}}
         \AUTHOR{Mabel C. Chou}
		\AFF{
        Department of Analytics and Operations, NUS Business School, Singapore 119245,
        \EMAIL{mabelchou@nus.edu.sg}}
        \AUTHOR{Jong-Shi Pang}
		\AFF{
        Department of Industrial and Systems Engineering, University of Southern California, 90089
        Los Angeles, USA,
        \EMAIL{jongship@usc.edu}}
        % \today
	} % end of the block
    
\ABSTRACT{Offline policy learning aims to use historical data to learn an optimal personalized decision rule. In the standard estimate-then-optimize framework, reweighting-based methods (e.g., inverse propensity weighting or doubly robust estimators) are widely used to produce unbiased estimates of policy values. However, when the propensity scores of some treatments are small, these reweighting-based methods suffer from high variance in policy value estimation, which may mislead the downstream policy optimization and yield a learned policy with inferior value. In this paper, we systematically develop an offline policy learning algorithm based on a weight-clipping estimator that truncates small propensity scores via a clipping threshold chosen to minimize the mean squared error (MSE) in policy value estimation. Focusing on linear policies, we address the bilevel and discontinuous objective induced by weight-clipping-based policy optimization by reformulating the problem as a Heaviside composite optimization problem, which provides a rigorous computational framework. The reformulated policy optimization problem is then solved efficiently using the progressive integer programming method, making practical policy learning tractable. We establish an upper bound for the suboptimality of the proposed algorithm, which reveals how the reduction in MSE of policy value estimation, enabled by our proposed weight-clipping estimator, leads to improved policy learning performance.
}%


\maketitle
\newpage
\section{Introduction}
In a wide range of domains, the abundance of individual-level data creates opportunities to learn personalized decision-making rules. For example, in precision medicine, an optimal therapy should be tailored to a patient’s clinical profile \citep{bertsimas2017personalized,cui2017tree}. In online marketplaces, personalized pricing or assortment can be designed based on customers’ demographic characteristics  \citep{chen2022statistical, elmachtoub2023balanced,xie2025personalized,tang2025offline}. In public policy, governments can determine the allocation of public services or program access to individuals based on their socioeconomic characteristics \citep{kitagawa2018should,mbakop2021model}. Motivated by such applications, researchers in causal inference, operations research, and statistics have developed a framework known as \textit{offline policy learning}, which uses only pre-existing historical data to learn a mapping from individual characteristics to treatment assignments with the goal of maximizing an outcome of interest \citep[e.g.,][]{athey2021policy,dudik2011doubly, fang2025treatment, jin2025policy, kallus2018balanced, kallus2021minimax, kitagawa2018should, qi2019estimation, qi2023robustness, ye2025deep, zhan2024policy, zhang2025personalized, zhou2023offline}.

Typical offline policy learning algorithms consist of two steps: (i) policy evaluation and (ii) policy optimization. The first step aims to construct accurate value estimators for all admissible policies, and the second step finds a policy with the highest estimated value over a given policy class. Therefore, the performance of the downstream policy optimization critically depends on the accuracy of policy evaluation. A key challenge in policy evaluation is that an individual’s reward is observed only under the treatment actually received, but counterfactual rewards under alternative treatments are unobserved. Inverse propensity weighting (IPW) \citep{horvitz1952generalization} and doubly robust (DR) \citep{robins1994estimation} estimators\endnote{The DR estimator is also referred to as the Augmented IPW (AIPW) estimator, reflecting its construction as an IPW estimator augmented with a reward regression term.} are two classic reweighting-based estimators used to produce unbiased policy value estimates. Due to their unbiasedness, a number of policy learning algorithms have been developed based on these estimators (\citealp[e.g.,][]{athey2021policy,dudik2011doubly,fang2025treatment,kitagawa2018should,zhou2023offline}). However, as recognized in prior studies \citep{kallus2018balanced,kang2007demystifying,zhan2021off,zhan2024policy}, IPW and DR suffer from high variance, particularly when evaluating a target policy that is not well explored in the data, i.e., when the probability of observing the treatments assigned by a target policy is small in the data set. This lack of exploration can happen when some treatments are costly or restricted by ethical considerations; see, e.g., Table 9 of \cite{wei2023transfer}, which presents an electronic medical record for vasopressor therapy where treated patients are far fewer than untreated patients.

\begin{figure}[t]
    \centering
    \caption{
        A Toy Example 
    }
    \includegraphics[width=0.8\linewidth]{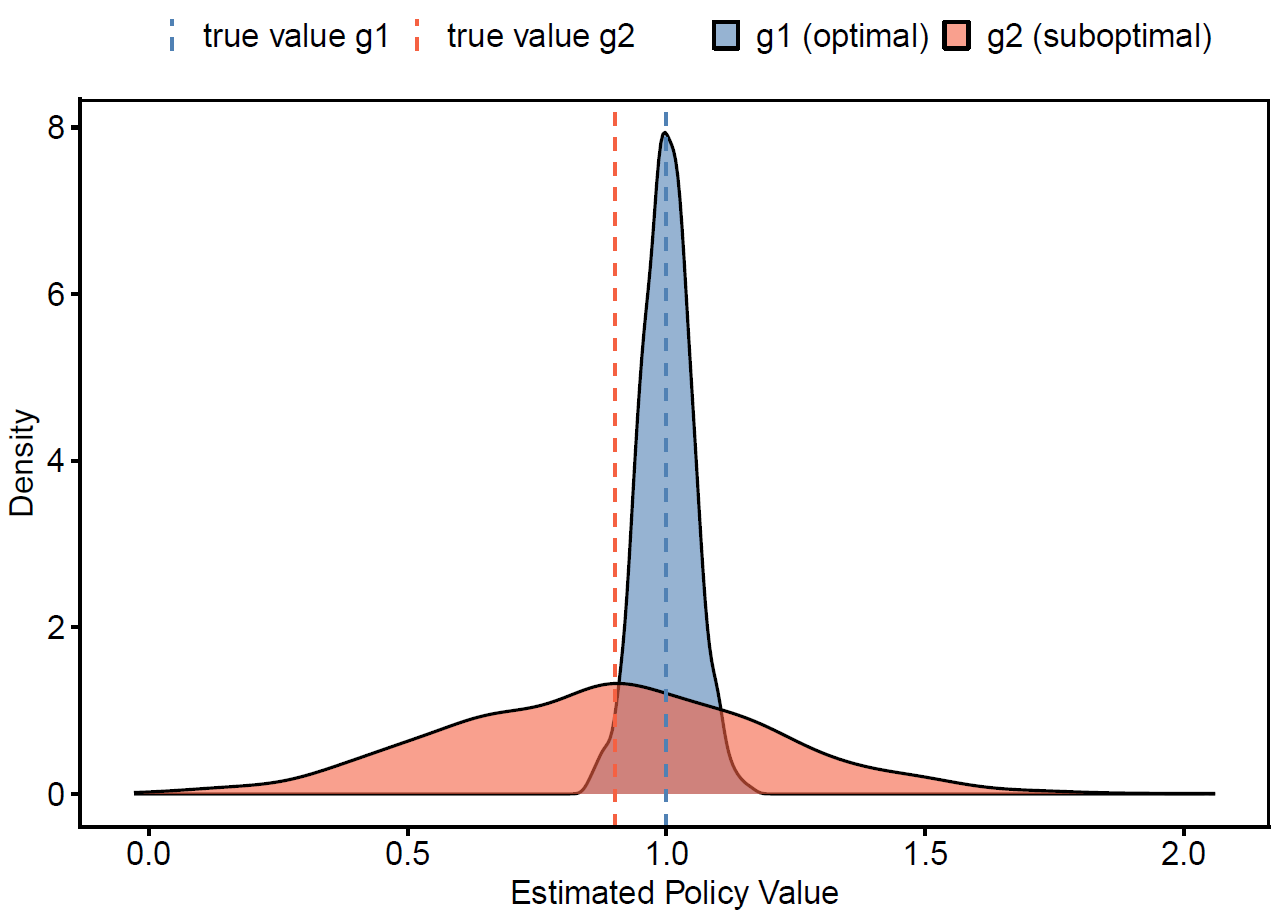}
    
    \label{fig:toy_example}
   \vspace{-8mm}
   \floatfoot{\textit{Note.}  The shaded region depicts the distribution of the estimated policy values under two estimators for 1{,}000 sampled individuals. In this instance, the learning algorithm selects the suboptimal policy $g_2$ with probability 36.5\%.}    
\end{figure}

From the perspective of a policy learner, the consequence of high variance in policy evaluation is that the learning algorithm may be misled into selecting a policy that appears promising but is in fact suboptimal. For illustration, consider a toy example in Figure \ref{fig:toy_example}. Suppose there are two policies, $g_1$ and $g_2$, where $g_1$ is optimal and $g_2$ is suboptimal. Because $g_1$ is well explored while $g_2$ is under-explored, the value estimator for $g_2$ has larger variance, and the algorithm may select $g_2$ due to the seemingly larger estimate, resulting in unsatisfactory realized rewards.

Recognizing that high variance in policy evaluation is driven by small propensity scores, a line of research has adopted the idea of \textit{weight clipping}, which truncates small propensity scores in the reweighting-based estimators to reduce variance (at the expense of a small bias) in policy evaluation \citep{bottou2013counterfactual,ma2022minimax,swaminathan2015batch, su2019cab,su2020doubly, wang2017optimal}. However, existing weight-clipping estimators typically select the clipping threshold via data-dependent tuning, rather than providing a closed-form expression. This makes it difficult to embed weight clipping into downstream policy optimization. Specifically, weight clipping induces a bilevel structure: we seek a policy that maximizes the clipped value estimator, while the (policy-dependent) clipping threshold is chosen to minimize estimation error. This entanglement can be resolved if the optimal clipping threshold admits an explicit expression as a function of a given policy. A second difficulty in applying weight-clipping estimators to policy optimization arises on the computational side: because both the clipping rule and the treatment rule are discrete, the resulting objective is discontinuous, raising concerns about the tractability of computing high-quality solutions. Together, these issues hinder the direct use of weight-clipping estimators for efficient policy learning. In this paper, we systematically address these challenges in weight-clipping-based policy learning over a class of linear policies. Below, we summarize our main results and contributions.

\subsection{Main Results and Contributions}
Our proposed policy learning framework is built upon a newly developed weight-clipping estimator that enables an optimal balance of the bias-variance trade-off in policy evaluation. Meanwhile, we resolve the optimization and computational difficulties and provide statistical guarantees. Our contributions are fivefold:

% \vspace{0.1in}

(a) \textbf{A novel estimator with optimized weight clipping.}  To mitigate the high variance in policy evaluation --- which may mislead the learning algorithm into selecting an optimistically estimated policy  --- we propose a novel weight-clipping estimator that truncates small propensity scores in the DR estimator. Unlike existing weight-clipping estimators, which find the clipping via data-dependent tuning, our estimator determines this threshold as the \textit{closed-form} minimizer of the mean squared error (MSE) in policy evaluation, thereby rigorously balancing the bias-variance trade-off. It is important to note that the closed-form expression for the optimal clipping threshold is derived for any admissible policy, which is crucial for solving the bilevel problem in policy optimization, as will be discussed shortly.
We refer to the proposed estimator as the Optimized Clipped Doubly Robust (OCDR) estimator, which is defined in Section \ref{sec:ocdr_definition}.

% \vspace{0.1in}

(b) \textbf{Optimization.} Although the proposed OCDR estimator is appealing for policy evaluation, the downstream policy optimization problem is bilevel and discontinuous, making it difficult to compute the optimal policy within a linear policy class, in which the policies are not enumerable (as discussed in Section \ref{sec:cdr_based_policy_optimization}). To address these challenges, we first use the closed-form expression of the optimal clipping threshold to collapse the bilevel optimization problem into a single-level problem in Section \ref{sec:tau_optimization}. We then show that this problem can be reformulated as a Heaviside composite optimization problem (HSCOP) \citep{cui2023minimization, han2024analysis} in Section \ref{sec:HSOCP_formulation}, which provides a rigorous mathematical formulation for learning multi-treatment rules and efficient computation.

% \vspace{0.1in}

(c) \textbf{Computation. }Instead of solving the reformulated policy optimization problem via a direct but potentially impractical mixed-integer program (MIP) formulation, we adapt the progressive integer programming (PIP) algorithm \citep{fang2025treatment} to our weight-clipping-based policy optimization setting, which involves a more complicated reformulation (as formulated in Section \ref{sec:pip_formulation}). The PIP procedure iteratively solves a sequence of MIP subproblems of 
reduced sizes, improving solution quality over iterations. This enables us to efficiently compute a desirable policy at termination of the procedure.

% \vspace{0.1in}

(d) \textbf{Statistical properties of the proposed policy learning algorithm.} By resolving the optimization and computational challenges, we complete the development of a practically tractable policy learning algorithm built upon the OCDR estimator, as presented in Section \ref{sec:policy_learning_algo}. In Section \ref{sec:stats_performance}, we analyze the statistical properties by establishing an upper bound on the suboptimality of the learned policy, which reveals how the generalization error in policy learning benefits from reducing the (worst-case) MSE in policy evaluation, as achieved by the proposed OCDR estimator.

% \vspace{0.1in}
(e) \textbf{Numerical and practical performance.} In Section~\ref{sec:numerical}, we conduct a series of experiments to demonstrate the computational efficiency and statistical superiority of the proposed framework. We further apply our weight-clipping-based policy learning framework to a real-world problem of learning targeting policies for insurance products, showing that it achieves higher expected profit than baseline approaches.

\subsection{Related Work}
\textit{Weight clipping for offline policy learning.} To address the high variance issue in the reweighting-based approaches, weight clipping is designed to exclude small propensities (\citealp[e.g.,][]{bottou2013counterfactual,ma2022minimax, swaminathan2015batch, su2019cab,su2020doubly, wang2017optimal}). We identify two research gaps in the existing literature. First, existing weight-clipping approaches typically select the clipping threshold via data-dependent tuning, rather than providing a closed-form expression for the optimal threshold. Without such a closed-form characterization, applying weight-clipping estimators to policy optimization is difficult, since the optimization objective depends on an additional threshold that must be tuned alongside the policy. Second, there is very limited analysis of policy learning performance based on weight-clipping estimators. A more relevant exception is \cite{su2019cab}, which evaluates learning performance through numerical experiments without providing theoretical guarantees.  We fill these gaps by (i) proposing a new weight-clipping estimator adapted from the SWITCH estimator \citep{wang2017optimal}, with the key improvement that the clipping threshold (which minimizes the MSE in policy evaluation) admits a closed-form expression as a function of any policy within an arbitrary policy class, and (ii) providing the first theoretical analysis for weight-clipping-based policy learning.

\textit{Techniques for policy optimization.} The exact techniques for policy optimization depend on the choice of policy class and the estimator. For example, \cite{zhou2023offline} propose a tree-search algorithm to maximize the DR estimator under a $k$-depth tree policy class. \cite{kallus2018balanced} develops a bilevel policy optimization formulation that maximizes a reweighting-based estimator, where the weights minimize the policy estimation error under smoothness assumptions on the reward function. Such smoothness assumptions are not required in our setting. The work most closely related to ours is \cite{fang2025treatment}, which maximizes the reward estimated by the IPW estimator over a linear policy class. In particular, the authors show that such a policy optimization problem can be reformulated as a Heaviside composite optimization problem (HSCOP) and propose a progressive integer programming (PIP) algorithm to solve it efficiently. Inspired by \cite{fang2025treatment}, we reformulate the more complex weight-clipping optimization problem as an HSCOP. Unlike IPW-based formulations, we incorporate the weight-clipping rule and handle the undetermined signs of the constant coefficients associated with the Heaviside functions via an upper-semicontinuous approximation. We also tailor the PIP algorithm to our formulation to compute the optimal policy efficiently.

\section{Preliminaries on Offline Policy Learning}\label{sec:preliminaries}
\subsection{Offline Data Description}
We consider offline policy learning from an observational data set consisting of tuples of random variables, denoted by $\cal{D}= \{X^s, D_s, Y_s\}_{s=1}^N$, where $N$ is the sample size of the data set. Specifically, $X^s \overset{\text{i.i.d.}}{\sim} \Prob_X$ are random covariate vectors that lie in the covariate space $\cal{X} \in \R^{p}$. Given a realized value $\xi$ of $X^s$, a treatment $D_s$ is sampled from the set $[J] \triangleq \{1, \cdots, J\}$ according to a behavior policy, which is characterized by the propensity score $e(\xi,j) \triangleq \Prob(D_s = j \mid X^s = \xi)$. In our setting, $e(\bullet, \bullet)$ is unchanged for all samples and assumed to be known. Associated with the chosen treatment $D_s$, an outcome (reward) $Y_s = \mu(X^s, D_s) + \epsilon_s$ is observed, where $\mu(\bullet, \bullet)$ is an unknown function and $\{\epsilon_s\}_{s=1}^N$ are i.i.d. mean-zero random variables. Under the potential outcomes framework \citep{neyman1923application,rubin1974estimating}, we let $\{Y_s(j)\}_{j \in [J]}$ denote the potential outcomes that we would have observed by taking treatment $j$. We note that $Y_s = Y_s(D_s)$, meaning that only the potential outcome corresponding to the actually assigned treatment is observed. This feature creates the key difficulty in answering counterfactual questions about the potential outcome under a hypothetical treatment assignment. 

We assume that the observed data satisfy the following standard assumptions in causal inference.
\begin{assumption}\label{asp:dgp}
The data-generating process satisfies the following properties:
\begin{itemize}
    \item[(A)] Unconfoundedness: for every treatment $j \in [J]$, the outcome $Y_s(j)$ is independent of the treatment $D_s$ given the covariate $X^s$, i.e., $Y_s(j) \indep D_s\mid X^s$.
    \item[(B)] Overlap: there exists some $\eta > 0$, such that for every treatment $j \in [J]$, the propensity score $e(\xi,j) \geq \eta$, $\forall \xi \in \Xcal$. 
    \item[(C)] Bounded reward: $Y_s \in [0,M]$ (for $M \geq 1$) almost surely.
\end{itemize}
\end{assumption}
Assumption \ref{asp:dgp}(A) means that the treatment taken is independent of all the reward outcomes conditioned on the covariate vector. Assumption \ref{asp:dgp}(B) requires that each treatment has a minimum probability of being selected for any covariate value. It is worth noting that the level of overlap $\eta$ may be small (often referred to as ``weak overlap"), which is a primary source of large inverse propensity scores and high variance in policy value estimation. Assumption \ref{asp:dgp}(C) is imposed to simplify analysis, and it aligns with many practical applications. 

\subsection{Policy Learning and Performance Metric}\label{sec:Policy Learning and Performance Metric}
Within a (prespecified) policy class $\Gcal$, a policy $g: \Xcal \to [J]$ is a mapping that maps an individual's covariates $\xi \in \Xcal$ to a treatment $j \in [J]$. We define the policy value of a given policy $g$ as the expected reward received when treatments are assigned according to $g$: $\Psi(g) \triangleq \E[Y(g(X))] = \E[\mu(X, g(X))]$,
where the expectation is taken with respect to $\Prob_X$. The optimal policy $g^*$ is a policy that maximizes the policy value, i.e., $g^* = \argmax_{g \in \Gcal}\Psi(g)$.

For computational and practical tractability, we restrict attention to the case where $\mathcal{G}$ is a linear policy class (LPC), denoted by $\mathcal{G}_{\mathrm{LPC}}$. Specifically, 
\begin{equation*}
    \mathcal{G}_{\mathrm{LPC}} \triangleq \left \{g: g(\xi) = \argmax_{j \in [J]}~(\xi^\top \beta^j + b_j), \xi \in \R^{p}, \beta^j \in \R^{p}, b_j \in \R \right\},
\end{equation*}
where $b_j$ is a known scalar considered as the ``base score" of treatment $j$ and each $\beta^j \in \R^p$ is a treatment-dependent parameter that needs to be learned. In $\mathcal{G}_{\mathrm{LPC}}$, the coefficient parameters $\beta^j$ measure the importance of a given covariate for a specific treatment. The benefit of using $\mathcal{G}_{\mathrm{LPC}}$ is that it provides an interpretable decision rule, which is implemented in many applications. For example, in public policy, job training programs are often allocated using simple, interpretable linear policies based on individuals' characteristics such as education and prior earnings \citep{kitagawa2018should}. In precision medicine, linear policies help medical practitioners understand the extent to which each patient characteristic influences the effect of treatment \citep{chen2024robust}. In personalized pricing, \cite{shen2025proxy} consider the application of linear policies that depend on product brands.

Using the offline data, our objective is to learn a policy $\wh{g} \in \Gcal_{\mathrm{LPC}}$ that minimizes the suboptimality gap $\cal{L}(\wh{g}, \Gcal_{\mathrm{LPC}}) \triangleq \Psi(g^*) - \Psi(\wh{g})$, which is the difference between the reward values generated by policies $\wh{g}$ and $g^*$ from $\Gcal_{\rm LPC}$.

\subsection{Policy Evaluation Using Reweighting-Based Methods}
Consider a widely used estimate-then-optimize framework for policy learning. The policy evaluation step aims to construct a reliable estimator $\wh{\Psi}(g)$ of the policy value for each policy $g \in \mathcal{G}_\mathrm{LPC}$, and the policy optimization step finds a policy $\wh{g} \in \cal{G}_{\mathrm{LPC}}$ by solving the maximization problem $\displaystyle{\max_{g \in \mathcal{G}_\mathrm{LPC}}} \wh{\Psi}(g)$. Apparently, good learning performance relies on accurate evaluation of all policies in $\mathcal{G}_\mathrm{LPC}$, including policies not actually implemented (i.e., the counterfactual policies). Reweighting-based methods are widely used for this purpose, among which the DR estimator \citep{dudik2011doubly} is one of the most well-known representatives in policy evaluation.\endnote{The DR estimator is defined as $ \widehat{\Psi}^{\text{DR}}(g) = \frac{1}{N} \sum_{s=1}^N \widehat{\Gamma}^{\mathrm{DR}}_s (g)$, where the objects being averaged are doubly robust scores, such that $\widehat{\Gamma}^{\mathrm{DR}}_s (g) = \widehat{\mu}(X^s,g(X^s)) + \frac{\onebld\{g(X^s) = D_s\}}{e(X^s,g(X^s))} \cdot (Y_s - \widehat{\mu}(X^s, g(X^s)))$,
where $\widehat{\mu}(\bullet, \bullet): \Xcal \times [J] \to [0,M]$ is an estimator of the expected reward $\mu(\xi, j)$ and $\onebld(\bullet)$ is the indicator function.} However, when the overlap is weak in practice (i.e., $\eta$ in Assumption \ref{asp:dgp}(B) is small), the variance of the DR estimator may explode, leading to high MSE. Such a behavior is formally characterized in Proposition \ref{prop: Dominant Term of Variance}.\endnote{Proposition \ref{prop: Dominant Term of Variance} naturally applies to the IPW estimator \citep{horvitz1952generalization}.}

\begin{proposition}\label{prop: Dominant Term of Variance}
    Denote $\widehat{\Gamma}^{\mathrm{DR}}_s(g)$ by the doubly robust score in the DR estimator. Under Assumption \ref{asp:dgp}, suppose $L' \leq \Var(Y_s \mid X^s)  \leq U'$ for all $s \in [N]$ and $U'\geq L' \geq 0$.
    For any policy $g \in \Gcal_\mathrm{LPC}$, there exist positive constants $U \geq L \geq 0$, such that $\frac{L}{e(X^s,g(X^s))} \leq \Var(\widehat{\Gamma}^{\mathrm{DR}}_s(g) \mid X^s) \leq \frac{U}{e(X^s,g(X^s))}, ~\forall s \in [N].$
\end{proposition}

Proposition \ref{prop: Dominant Term of Variance} shows that the variance of the DR estimator scales proportionally with the inverse propensity score (IPS) $e(X^s,g(X^s))^{-1}$. Therefore, if a treatment $j$ prescribed by a policy $g \in \Gcal_\mathrm{LPC}$ is rarely observed in the offline data, the IPS becomes extremely large, inflating the variance of the DR estimator. 

\section{Weight-Clipping-Based Policy Learning: Formulation and Challenges}\label{sec:formulation}
In light of Proposition \ref{prop: Dominant Term of Variance}, the inferior performance of the DR estimator in policy evaluation can be attributed to the high variance caused by the large IPSs. Therefore, a natural remedy for variance reduction is to clip the large IPSs; see the references in the previous sections. However, to the best of our knowledge, existing work has not systematically developed an offline policy learning framework based on weight-clipping estimators, particularly one that addresses the associated optimization challenges and establishes theoretical properties. Motivated by this gap, we propose a framework for weight-clipping-based policy learning, called Optimized Clipped Doubly Robust Learning (OCDRL), for efficient learning under the weak overlap condition. The OCDRL framework consists of two key components: (i) policy evaluation using an Optimized Clipped Doubly Robust (OCDR) estimator, whose clipping threshold is optimally chosen to minimize the MSE in policy evaluation; and (ii) policy optimization formulated as a Heaviside composite optimization problem, together with an efficient computational method for solving it. We discuss each component in detail below.

\subsection{(Optimized) Clipped Doubly Robust Estimator}\label{sec:ocdr_definition}
We first introduce the clipped doubly robust (CDR) estimator, which is a direct adaptation of the SWITCH estimator \citep{wang2017optimal} without optimizing the clipping threshold. Analyzing the estimation error of this baseline CDR estimator then motivates the construction of an optimized clipping threshold.

\begin{definition}[Clipped Doubly Robust Estimator]\label{def:cdr}
Given a clipping threshold $\tau \in [0, 1/\eta]$, the CDR estimator is defined as
\begin{align*}
\widehat{\Psi}^{\mathrm{CDR}}(g) = \frac{1}{N} \sum_{s=1}^N  \wh{\Gamma}_s^{\mathrm{CDR}}(g,\tau),
\end{align*}
where the objects being averaged are clipped doubly robust scores, such that
    \begin{align}\label{eq:CDR score}
    \wh{\Gamma}_s^{\mathrm{CDR}}(g,\tau) &= \widehat{\mu}(X^s,g(X^s)) + \frac{\onebld\{g(X^s) = D_s\}}{e(X^s,g(X^s))} \cdot \onebld\left\{\frac{1}{e(X^s,g(X^s))} \leq \tau \right\} \cdot (Y_s - \widehat{\mu}(X^s, g(X^s))),
    \end{align}
    where $\widehat{\mu}(\bullet, \bullet) \, : \,
    \Xcal \, \times \, [J] \, 
    \to \, [0,M]$ is an estimator for the expected reward $\mu(\xi, j)$ and $\onebld(\bullet)$ is the indicator function.
\end{definition}

From Definition \ref{def:cdr}, we see that the CDR estimator interpolates between the direct estimate of the reward model $\widehat{\mu}(\bullet, \bullet)$ (a.k.a. the Direct Method (DM)) and the DR estimator.
The design principle of the CDR estimator is straightforward: given a policy $g$, when the IPS $e(X^s,g(X^s))^{-1}$ is small, we use the DR estimator to take advantage of its unbiasedness. However, when the IPS $e(X^s,g(X^s))^{-1}$ is large, we switch to DM to reduce variance with only a small bias cost. It is worth highlighting that, here, ``small" and ``large" are determined by the clipping threshold $\tau$, which needs to be determined carefully. Conceptually, a good threshold parameter $\tau$ should minimize the MSE in policy evaluation, which is defined as $\E\big[\big(\widehat{\Psi}(g) - \Psi(g)\big)^2\big]$, with the expectation being taken with respect to the randomness of the offline data.

\begin{proposition}\label{prop:MSE_CDR}
For a given $\tau \in [0, 1/\eta]$ and a fixed policy $g \in \Gcal_{\rm LPC}$, the MSE of the CDR estimator in policy evaluation is upper bounded by
\begin{align}
   \label{MSE upper bound}M^2 \left(\underbrace{\E\left[\onebld\left\{\frac{1}{e(X,g(X))}>\tau\right\}\right]^2}_{(\mathrm{i})} + \underbrace{\frac{2}{N} \E\left[\frac{\onebld\{g(X) = D\}}{e(X,g(X))^2} \cdot \onebld\left\{\frac{1}{e(X,g(X))} \leq \tau \right\}\right] + \frac{2}{N}}_{(\mathrm{ii})}\right).
\end{align}
\end{proposition}

As shown in Proposition \ref{prop:MSE_CDR}, the MSE upper bound consists of two dominant terms: term (i) is caused by the bias in estimation, and term (ii) is due to variance in estimation. In addition, 
Proposition~\ref{prop:MSE_CDR} also reveals how the parameter $\tau$ balances the bias-variance trade-off. As 
$\tau \, \downarrow \, 0$, 
term (ii) in the MSE upper 
bound (\ref{MSE upper bound}) vanishes, leaving 
only the bias dominating the MSE. In contrast, as 
$\tau \, \uparrow \, \infty$,
term (i) vanishes, and the MSE becomes dominated by variance. Therefore, it is natural to consider selecting an appropriate $\tau$ that strikes the balance of bias and variance. 

In practice, since the expected values cannot be computed exactly, we minimize the sample-average approximation of the MSE upper bound. Formally, given a policy $g$, we solve the following optimization problem:
\begin{align}\label{program: optimal tau}
    \min_{\tau \in [0, 1/\eta]} ~\left(\frac{1}{N} \sum_{s=1}^N \onebld\left\{\frac{1}{e(X^s,g(X^s))}>\tau\right\}\right)^2 + \frac{2}{N^2} \sum_{s = 1}^N \frac{\onebld\{g(X^s) = D_s\}}{e(X^s,g(X^s))^2} \cdot \onebld\left\{\frac{1}{e(X^s,g(X^s))} \leq \tau \right\}.
\end{align}

We denote by $\wh{\tau}(g)$ the minimizer of 
problem~\eqref{program: optimal tau}, 
which will be shown to be unique under a mild assumption on the
propensity scores. We are now ready to formally define our proposed Optimized Clipped Doubly Robust (OCDR) estimator.
\begin{definition}[Optimized Clipped Doubly Robust Estimator]\label{def:OCDR}
For any given policy $g$, let $\wh{\tau}(g)$ be the minimizer of problem~\eqref{program: optimal tau}, the OCDR estimator is defined as
   \begin{align*}
    \widehat{\Psi}^{\mathrm{OCDR}}(g) = \frac{1}{N} \sum_{s=1}^N  \wh{\Gamma}_s^{\mathrm{OCDR}}(g,\wh{\tau}(g)),
\end{align*}
where the objects being averaged are the optimized clipped doubly robust scores, such that
\begin{equation}
\label{eq:OCDR_scores}
        \wh{\Gamma}_s^{\mathrm{OCDR}}(g,\wh{\tau}(g)) = \widehat{\mu}(X^s,g(X^s))  + \displaystyle{
        \frac{\onebld\{g(X^s) = D_s\}}{e(X^s,g(X^s))} 
        } \times \onebld\left\{\frac{1}{e(X^s,g(X^s))} \leq \wh{\tau}(g) \right\} \times (Y_s - \widehat{\mu}(X^s, g(X^s))),
 \end{equation}
    where $\widehat{\mu}(\bullet, \bullet) \, : \,
    \Xcal \, \times \, [J] \, 
    \to \, [0,M]$ is the estimator for expected reward $\mu(\xi, j)$
    and $\onebld(\bullet)$ is the indicator function.
\end{definition}

We derive an explicit and computable form of the MSE minimizer $\wh{\tau}(g)$ (in Section \ref{sec:tau_optimization}), which makes the dependence of the clipping threshold on the policy $g$ transparent. This closed-form characterization of the optimal clipping threshold is a distinguishing feature that highlights the departure of our approach to existing weight-clipping estimators. We discuss the motivation for this design choice shortly, when outlining the associated challenges in policy optimization.

\subsection{Challenges in OCDR-Based Policy Optimization}\label{sec:cdr_based_policy_optimization}
Given that the OCDR estimator is promising in reducing the variance in policy evaluation, it is tempting to find the optimal policy by maximizing the value 
determined by the OCDR 
estimator. Unfortunately, the discontinuity of the objective function, together with the need to optimize the threshold parameter 
$\wh{\tau}(g)$, makes the optimization problem challenging. 

In policy optimization, focusing on the linear policy class $\mathcal{G}_{\mathrm{LPC}}$, we aim to find a solution to an empirical reward maximization problem, defined as
\begin{align}\label{program:policy optimization}
    \max_{g \in \mathcal{G}_{\mathrm{LPC}}} \widehat{\Psi}^{\mathrm{OCDR}}(g, \wh{\tau}(g)) = \max_{g \in \mathcal{G}_{\mathrm{LPC}}} \frac{1}{N} \sum_{s=1}^N \Gamma_s^{\mathrm{OCDR}}(g,\wh{\tau}(g)).
\end{align}
The first practical challenge we recognize is that problem \eqref{program:policy optimization} is a bilevel optimization problem. In the upper-level problem, we aim to find 
$\wh{g} \in \displaystyle{
\argmax_{g \in \mathcal{G}_{\mathrm{LPC}}}
} \, \wh{\Psi}^{\mathrm{OCDR}}(g, \wh{\tau}(g))$,
provided that $\wh{\tau}(g)$, the solution to the MSE minimization problem \eqref{program: optimal tau}, can be evaluated for any given policy 
$g \in \mathcal{G}_{\mathrm{LPC}}$. However, since the policies in $\mathcal{G}_{\mathrm{LPC}}$ are not enumerable, simple brute-force methods (e.g., grid search) fail to solve such a bilevel optimization problem. We address this issue by deriving an explicit
expression for the optimal solution $\wh{\tau}(g)$,
and then substituting this solution back into the upper-level objective 
\eqref{program:policy optimization}; 
the result is a transformation of the bilevel optimization problem into a 
single-level formulation. 

The second challenge we need to tackle is the discontinuous objective in policy optimization, incurred by the discrete nature of the OCDR score defined in Equation \eqref{eq:OCDR_scores} as well as the treatment rule in policy class $\mathcal{G}_{\mathrm{LPC}}$. Inspired by \cite{fang2025treatment}, we address this challenge by reformulating the policy optimization problem as a \textit{Heaviside composite optimization problem} (HSCOP), whose definition (for our purpose) is given below.
\begin{definition}[Heaviside Composite Optimization Problem]\label{def:HSCOP}
\begin{align}\label{General HSCOP}
    \max_{x \in \cal{P}} f_1(x) + f_2(x), \quad \text{where} f_2(x) \triangleq \sum_{i=1}^N \psi_i \onebld_{[0,\infty)}(\phi_i(x)),
\end{align}
where $\cal{P}$ is a given polyhedral subset of $\R^d$, $f_1: \R^d \to \R$ is a Bouligand differentiable function, the $\psi_i$'s are nonnegative constants, each $\phi_i : \R^d \to \R$ 
is a (multivariate) piecewise affine function, and
$\onebld_{[ \, 0,\infty )}(t) \triangleq \left\{
\begin{array}{ll}
1 & \mbox{if $t \geq 0$} \\
0 & \mbox{if $t < 0$}
\end{array} \right\}$ is the (binary-valued)
Heaviside function
of the closed interval $[ \, 0,\infty )$.
\end{definition}

From a computational perspective, an attractive feature of the HSCOP is that it can be solved via a sequence of restricted mixed-integer program (MIP) formulations \citep{fang2025treatment}. This significantly reduces the number of integer variables required,  compared with a large-scale full MIP formulation. Specifically, given a feasible $\Bar{x}$ to problem \eqref{General HSCOP}, define three complementary index sets:
$\cal{I}^{<}_{\delta_2} \triangleq \{i \mid \phi_i(\Bar{x}) < -\delta_2\},
\cal{I}^{\mathrm{in}}_{\delta_{1,2}} \triangleq \{i \mid -\delta_2 \leq \phi_i(\Bar{x}) \leq \delta_1\},
\cal{I}^{>}_{\delta_{1}} \triangleq \{i \mid \phi_i(\Bar{x}) > \delta_1\},$ where $\delta_1, \delta_2 \geq 0$. An important observation is that if $\phi_i(\bar{x})$ is nonzero, then for all $x$ sufficiently close to $\bar{x}$, the sign of $\phi_i(x)$ must coincide with the sign of $\phi_i(\bar{x})$ by continuity. Conversely, if $\phi_i(\bar{x}) = 0$, the sign of $\phi_i(x)$ in a neighborhood of $\bar{x}$ cannot be determined. Therefore, when introducing an integer variable $z_i \in \{0,1\}$ to represent the Heaviside function $\onebld_{[0,\infty)}(\phi_i(x))$, it suffices to fix $z_i = 0$ for indices $i$ in $\cal{I}^{<}_{\delta_2}$ and $z_i = 1$ for indices $i$ in $\cal{I}^{>}_{\delta_{1}}$, while retaining $z_i$ as a binary decision variable only for indices $i \in \cal{I}^{\mathrm{in}}_{\delta_{1,2}}$. This observation motivates the algorithmic development of the progressive integer programming (PIP) method \citep{fang2025treatment}, which iteratively solves a sequence of MIPs with adaptive control of the parameter pair $\{\delta_1, \delta_2\}$ until the objective values do not improve. With guarantees of local optimality, the PIP method makes large-scale policy optimization computationally tractable and efficient.

\section{Mathematical Reformulations for Policy Optimization}\label{sec:math formulation}
In this section, we describe how to reformulate 
the policy optimization 
problem \eqref{program:policy optimization} for practical computation. The first step is the
derivation of the optimal clipping threshold $\wh{\tau}(g)$
for any $g \in \cal{G}_{\mathrm{LPC}}$, which transforms the bilevel weight-clipping-based policy optimization problem into a single-level one. It is worth noting that, while we focus on the linear policy class $\mathcal{G}_{\mathrm{LPC}}$ in this paper, the explicit form of $\wh{\tau}(g)$ derived applies to arbitrary policy classes.

\subsection{Explicit Solution to the MSE Minimization Problem (\ref{program: optimal tau})}\label{sec:tau_optimization}
Our strategy for characterizing the explicit solution $\wh{\tau}(g)$ proceeds in three steps. First, we reformulate problem \eqref{program: optimal tau} by constructing a new objective whose summands are indexed by the sorted IPS values. This reformulation enables a precise characterization of the optimality conditions and a computable expression for $\wh{\tau}(g)$, which are presented in the second step. Finally, we use a pruning strategy to reduce the computational burden.

\textbf{Step 1: Reformulate problem (\ref{program: optimal tau}). }
For notational simplicity, define the IPS as 
$C_s \triangleq [ e(X^s,g(X^s)) ]^{-1}$ and $\varphi_s(g) \triangleq C_s \cdot \onebld\{g(X^s) = D_s\}$. We sort the samples according to the IPS values $\{C_s\}_{s=1}^N$ in ascending order and add $C_{(0)}=0$, such that $0 = C_{(0)} \leq C_{(1)} \leq  C_{(2)} \leq  \cdots \leq  C_{(N)}$. Define $\varphi_{(s)}(g) \triangleq  (C_{(s)})^2 \cdot \onebld \{g(X^{(s)}) = D_{(s)}\}$,  corresponding to the sample with the $s$-th smallest IPS $C_{(s)}$. Given a policy $g$, the objective of problem \eqref{program: optimal tau} can be reformulated as follows,
\[
\begin{array}{ll}
\theta(\tau;g) & = \displaystyle{\frac{1}{N^2}} \left( \sum_{s=1}^N \onebld\{\tau < C_{(s)}\} + 2 \sum_{s<t}\onebld \{\tau < C_{(s)} \} \right) + \frac{2}{N^2} \sum_{s=1}^N \varphi_{(s)}(g)\left(1-  \onebld \{ \tau < C_{(s)}\} \right)\\[0.25in]
& = \displaystyle{\frac{1}{N^2}} \Bigg( \underbrace{\sum_{s=1}^N  \Big( 2 (N-s) +1 - 2 \varphi_{(s)}(g) \Big) \onebld \{\tau < C_{(s)} \}}_{ \mbox{denoted by $\wt{\theta}(\tau;g)$}  } \Bigg) + \frac{2}{N^2} \sum_{s=1}^N \varphi_{(s)}(g).
\end{array}
\]
Here, we have reformulated 
problem \eqref{program: optimal tau} with a new objective function $\wt{\theta}(\tau;g)$. Observe that the optimal solution to problem \eqref{program: optimal tau} is also 
a minimizer of $\wt{\theta}(\tau;g)$ for a given
policy $g$.

\textbf{Step 2: Characterize the optimality conditions.}
Let $\Phi_m(g) \triangleq \displaystyle{
\sum_{s=m}^N
} \, (2(N-s)+1 - 2 \varphi_{(s)}(g))$ and $\bar{\Phi}_{m,\ell}(g) \triangleq \Phi_m(g) -\Phi_{\ell}(g)$ %= \sum_{s=m}^{\ell-1} (2(N-s)+1 - 2 \varphi_{(s)}(g))$, 
for  $m, \ell =1,\ldots, N$.  It can be shown that
\[
\wt{\theta}(\tau;g) = \begin{cases} \Phi_m(g) & \mbox{ if } C_{(m-1)} \leq \tau < C_{(m)}, \quad \mbox{ for } m=1,2, \ldots, N, \\
\Phi_{N+1}(g) \equiv 0 & \mbox{ if } \tau \geq C_{(N)},\end{cases}
\]
which implies that the optimal value of 
$\wt{\theta}(\bullet;g)$ can be obtained as one
of the sorted IPSs $\{C_{(s)}\}_{s=0}^N$. The following mild assumption ensures the uniqueness of the minimizer of $\wt{\theta}(\bullet;g)$, which simplifies the expression for $\wh{\tau}(g)$.\endnote{In practice, $C_s$ may take integer values (e.g., under randomized assignment with fixed treatment probabilities). Nevertheless, by applying an arbitrarily tiny perturbation to the propensity scores, we can ensure that Assumption~\ref{asp:non-integer} holds, which in turn guarantees that the minimizer of $\wh{\theta}(\bullet;g)$ is unique.} We impose this assumption for the rest of the paper. The explicit expression of $\wh{\tau}(g)$ is given in Proposition \ref{prop:explict_expression_hat_tau}.

\begin{assumption}\label{asp:non-integer}
    The summation of any subset of 
    $\{2(C_s)^2\}_{s=1}^N$ is not an integer.
\end{assumption}

\begin{proposition}
\label{prop:explict_expression_hat_tau}
    Let $\wh{\tau}(g)$ be the minimizer of problem \eqref{program: optimal tau}.  
    Under Assumption \ref{asp:non-integer}, we have
\begin{equation} 
\label{eq:optimal_tau_simplified}
\begin{array}{ll}
\wh{\tau}(g)   
% & =\displaystyle{\sum_{m=0}^N C_{(m)}} \cdot \displaystyle{\onebld_{[0,\infty)}} \left(\min_{  t \in[N+1], t \neq m+1} \bar \Phi_{t,m+1}(g) \right) \\
=  \displaystyle{
\sum_{m=0}^N 
} \, C_{(m)} \cdot \displaystyle{\onebld_{(0,\infty)}} \left(\min_{ t \in[N+1], t \neq m+1}  \Phi_{t}(g) - \Phi_{m+1}(g) \right).
\end{array}
\end{equation}
\end{proposition}

Proposition \ref{prop:explict_expression_hat_tau} identifies a computable expression for $\wh{\tau}(g)$. However, computing $\Phi_m(g)$ for all $m = 1, \ldots, N + 1$ is computationally intensive when the sample size $N$ becomes large. This motivates a pruning strategy that reduces the number of values $\Phi_m(g)$ to be examined, which we implement in the final step.

\noindent\textbf{Step 3: Pruning of solution candidates. } The pruning strategy follows from Proposition \ref{prop: pruning}.

\begin{proposition}\label{prop: pruning}
    For the optimization problem with objective function $\wt{\theta}(\tau;g)$, if the minimum of $\Phi_m(g)$ is not achieved at index $m = N+1$, the minimum must be achieved at an index $m$ in the set $\Mcal \triangleq \{m = 1, \cdots, N \mid 2\wt{\varphi}_{(m)} \geq 2(N-m)+1\}$, where $\wt{\varphi}_{(s)}(g) \triangleq (C_{(s)})^2 \geq \varphi_{(s)}(g)$.
\end{proposition}

Proposition \ref{prop: pruning} implies that it suffices to examine the sorted IPS values with indices $m \geq m^*$, where $m^* \triangleq \min \{ m  \in [N] : 2\wt{\varphi}_{(m)} \geq 2(N-m)+1\}$. Under this pruning strategy, the explicit optimal solution $\wh{\tau}(g)$ can be computed efficiently.

\begin{proposition}\label{prop:final_expression_hat_tau} Under the pruning strategy and Assumption \ref{asp:non-integer}, we have
\begin{equation}
\label{eq:optimal_tau_simplified_pruning}
\begin{array}{ll}
\wh{\tau}(g)  %& =  \displaystyle{\sum_{m+1 \in \overline{\Mcal}}}  C_{(m)} \cdot \displaystyle{\onebld_{[0,\infty)}} \left(\min_{t \in \overline{M}, t \neq m+1} \bar \Phi_{t,m+1}(g) \right)\\
& =  \displaystyle{\sum_{m = m^*-1}^N}  C_{(m)} \cdot \displaystyle{\onebld_{(0,\infty)}} \left(\min_{N+1 \geq t \geq m^*, t \neq m+1}  \Phi_{t}(g) - \Phi_{m+1}(g) \right),
\end{array}
\end{equation}
where $m^* \triangleq \min \{ m  \in [N] : 2\wt{\varphi}_{(m)} \geq 2(N-m)+1\}$.
\end{proposition}

By plugging the optimal solution $\wh{\tau}(g)$ given in Equation \eqref{eq:optimal_tau_simplified_pruning} into the policy optimization problem \eqref{program:policy optimization}, we obtain a single-level formulation. It then remains to address the discontinuity in the objective function.

\subsection{HSCOP Formulation of the Policy Optimization Problem (\ref{program:policy optimization})}\label{sec:HSOCP_formulation}
We reformulate the policy optimization problem \eqref{program:policy optimization} as an HSCOP to address the discontinuity. For all $g \in \mathcal{G}_{\mathrm{LPC}}$ and given $\wh{\tau}(g)$ characterized by Equation \eqref{eq:optimal_tau_simplified_pruning}, the policy optimization problem \eqref{program:policy optimization} can be expressed as
\begin{equation}
\label{eq:obj_upper_level}
\begin{array}{ll} 
\displaystyle{ \sum_{s=1}^N} \sum_{j=1}^J \, 
\wh{\mu}(X^s,j) \cdot \onebld \{ g(X^s) = j\}    + \displaystyle{ \sum_{s=1}^{N}}   C_{s} ( Y_{s} -\wh{\mu} (X^{s}, g(X^{s}) )  )    \cdot  \onebld \{g(X^{s}) = D_{s}\} \cdot  \onebld \left\{\wh{\tau}(g)  \geq C_{s} \right\},
\end{array}
\end{equation}
where we invoke the definition of the OCDR estimator in Definition \ref{def:OCDR} and the fact that $\widehat{\mu}(X^s, g(X^s)) = \displaystyle{
\sum_{j=1}^J 
} \, \wh{\mu}(X^s, j) \onebld \{g(X^s) = j\}$. To obtain an HSCOP formulation, we reformulate the indicator functions as Heaviside composite functions, by leveraging the explicit expression of the optimal clipping threshold $\wh{\tau}(g)$ and the treatment rule of $\mathcal{G}_{\mathrm{LPC}}$.

The explicit form of $\wh{\tau}(g)$ given in Proposition \ref{prop:final_expression_hat_tau} implies that $\onebld \left\{\wh{\tau}(g)  \geq C_{(s)} \right\} = 1 $ for $s=1, \ldots, m^*-1$ and for all $g \in \mathcal G_\mathrm{LPC}$. In addition, for $s \geq m^*$, we have
\[
\begin{array}{ll}
  \onebld \left\{\wh{\tau}(g) \geq C_{(s)} \right\} %= \onebld_{[0,\infty)} \left( \displaystyle{\min_{t \in \overline{M}, m^* \leq t \leq s}}  \Phi_{t}(g) - \min_{t \in \overline{M}, t \geq s+1}\Phi_{t}(g) \right)   \\[0.2in] 
  & = \displaystyle{\sum_{N+1 \geq m  \geq s+1}}  \displaystyle{\onebld_{[0,\infty)}} \left(\min_{N+1 \geq t \geq m^*}  \Phi_{t}(g) -\Phi_{m}(g)  \right) \\[0.2in]
  & =  \displaystyle{\onebld_{[0,\infty)}} \left( \displaystyle{\max_{N+1 \geq m  \geq s+1}} 
 \left( \min_{N+1 \geq t \geq m^*}  \Phi_{t}(g) -\Phi_{m}(g) \right) \right)\\[0.2in]
  & =  \displaystyle{\onebld_{[0,\infty)}} \left(  
\min_{s  \geq t \geq m^*}  \Phi_{t}(g) - \displaystyle{\min_{N+1 \geq m  \geq s+1}}\Phi_{m}(g) \right),
 \end{array}
\]
where the second equality holds because $\wh{\tau}(g) \geq C_{(s)}$ if and only if the unique minimum value over 
$\{\Phi_m(g)\}_{m \geq m^*}$ is obtained within the subset $\{\Phi_m(g)\}_{m \geq s+1}$ for 
$s \geq m^*$. 

Recall the treatment rule of policy class $\mathcal{G}_{\mathrm{LPC}}$ specified in Section \ref{sec:Policy Learning and Performance Metric}, which assigns an individual with covariate $\xi \in \Xcal$ to the treatment $j \in [J]$ that maximizes the score function, i.e.,
\begin{equation*}
g(\xi) =j, \quad \text{if and only if $\left\{ \begin{array}{ll} \xi^{\top} \beta^j + b_j > \xi^{\top} \beta^i +b_i, \forall \, 1 \leq i \leq j-1 \\ \xi^{\top} \beta^j + b_j \geq  \xi^{\top} \beta^i +b_i, \forall \,  j+1 \leq i \leq J  \end{array} \right\}$;}
\end{equation*}
and thus
\begin{equation} \label{eq:linear_policy_formula}
\onebld\left\{ g(\xi) = j \right\} \, = \,\onebld_{( \, 0,\infty )}\left(\min_{ 1 \leq i \leq j-1}h_{j,i}(\xi,\boldsymbol{\beta})\right) \cdot  
\onebld_{[ \, 0,\infty )}\left(\min_{ j+1 \leq i \leq J } h_{j,i}(\xi,\boldsymbol{\beta})  \right) 
\end{equation}
with 
$\boldsymbol{\beta} \triangleq \{\beta^{\, j} \}_{j=1}^J$ and $
h_{j,i}(\xi,\boldsymbol{\beta}) \, \triangleq \, 
\xi^{\top} \beta^j + b_j -   ( \, \xi^{\top} \beta^i + b_i \, )$. Thus, the policy optimization problem \eqref{program:policy optimization} is equivalent to the following problem
\begin{equation}
   \label{eq:CDR_policy_optimization}
\max_{\boldsymbol{\beta} \in \mathbb R^p} 
\quad 
\wh{\Psi}_{\rm HSC}(\boldsymbol{\beta}),
\end{equation}
with the objective function $\wh{\Psi}_{\rm HSC}(\boldsymbol{\beta})$ involving Heaviside composite functions, given by
\[
\begin{array}{ll}
& \widehat{\Psi}_{\rm HSC}(\boldsymbol{\beta}) =\displaystyle{ \sum_{s=1}^N} \sum_{j=1}^J \, \wh{\mu}(X^s,j) \cdot \,\onebld_{( \, 0,\infty )}\left(\min_{ 1 \leq i \leq j-1}h_{j,i}(X^s,\boldsymbol{\beta})\right) \cdot  
\onebld_{[ \, 0,\infty )}\left(\min_{j+1 \leq i \leq J } h_{j,i}(X^s,\boldsymbol{\beta})  \right)   \\[0.2in]
     &  + \displaystyle{ \sum_{s=1}^{m^*-1}}   C_{s} ( Y_{s} -\wh{\mu} (X^{s}, g(X^{s}) )  )    \cdot    \,\onebld_{( \, 0,\infty )}\left(\displaystyle{ \min_{ 1 \leq i \leq D_{s}-1}} h_{D_{s},i}(X^{s},\boldsymbol{\beta})\right) \cdot  
\onebld_{[ \, 0,\infty )}\left(\displaystyle{ \min_{ D_{s}+1 \leq i \leq J }} h_{D_{s},i}(X^{s},\boldsymbol{\beta})  \right)   \\[0.2in]
     & + \displaystyle{ \sum_{s=m^*}^{N}}   C_{s} ( Y_{s} -\wh{\mu} (X^{s}, g(X^{s}) )  )    \cdot  \left( \begin{array}{ll} \,\onebld_{( \, 0,\infty )}\left(\displaystyle{ \min_{ 1 \leq i \leq D_{s}-1}} h_{D_{s},i}(X^{s},\boldsymbol{\beta})\right) \cdot  
\onebld_{[ \, 0,\infty )}\left(\displaystyle{ \min_{ D_{s}+1 \leq i \leq J }} h_{D_{s},i}(X^{s},\boldsymbol{\beta})  \right)   \\[0.1in] 
\cdot  \displaystyle{\onebld_{(0,\infty)}} \left(   
\min_{s  \geq t \geq m^*}  \Phi_{t}(g) - \displaystyle{\min_{N+1 \geq m  \geq s+1}}\Phi_{m}(g) \right) \end{array}  \right)  \\[0.2in]
\end{array}
\]
where $\Phi_{N+1}(\boldsymbol{\beta}) =0$, and for $m=1, \ldots, N$, 
\[\Phi_m(\boldsymbol{\beta}) = 2 \displaystyle{\sum_{s=m}^N} \left( (N-s)+\frac{1}{2} - (C_{s})^2 \cdot \onebld_{( \, 0,\infty )}\left(\displaystyle{ \min_{ 1 \leq i \leq D_{s}-1}} h_{D_{s},i}(X^{s},\boldsymbol{\beta})\right) \cdot  
\onebld_{[ \, 0,\infty )}\left(\displaystyle{ \min_{ D_{s}+1 \leq i \leq J }} h_{D_{s},i}(X^{s},\boldsymbol{\beta})  \right)    \right).\]

However, one may observe a mismatch: because the signs of 
$\{Y_{s} - \wh{\mu}(X^{s}, g(X^{s}))\}_{s=1}^N$ are undetermined, $\widehat{\Psi}_{\rm HSC}(\boldsymbol{\beta})$ is not upper semicontinuous, violating the requirement that all constant coefficients are nonnegative in the definition of HSCOP (see Definition \ref{def:HSCOP}). To address this issue, we construct an upper-semicontinuous (u.s.c.) $\varepsilon$-approximation of 
$\wh{\Psi}_{\mathrm{HSC}}(\boldsymbol{\beta})$ by utilizing Heaviside approximation functions  $\onebld_{( \, -\varepsilon,\infty )}(\bullet)$ and $\onebld_{[\,\varepsilon,\infty )}(\bullet)$  for a given small constant
$\varepsilon \geq 0$.  Clearly, 
$\onebld_{( \, -\varepsilon,\infty )}(t) \, \geq \,
\onebld_{( \, -\varepsilon^{\, \prime},\infty )}(t) \, \geq \,
\onebld_{[ \, 0, \infty )}(t),$ $\forall \, \varepsilon \, > \, \varepsilon^{\, \prime} 
\, > \, 0 \ \mbox{ and } \ \forall \, t \, \in \, \mathbb{R}$, and  $\displaystyle{
\lim_{\varepsilon \downarrow 0}
} \, \onebld_{(-\varepsilon,\infty)}(t) 
\, = \, \onebld_{[\,0,\,\infty)}(t)$ for all $t \, \in \, \mathbb{R}$.
Specifically, the right-hand-side of Equation
\eqref{eq:linear_policy_formula}  can be lower approximated by the u.s.c. function 
$\onebld_{[ \, 0,\infty )}(h^{1}_j(\xi, \boldsymbol{\beta};\varepsilon))$, 
and upper approximated by the lower-semicontinuous (l.s.c.) function 
$\onebld_{(\, 0,\infty )} (h^{2}_j(\xi, \boldsymbol{\beta};\varepsilon))$ as follows: 
\[\underbrace{\onebld_{[ \, 0,\infty )} (h^{1}_j(\xi, \boldsymbol{\beta};\varepsilon))}_{\mbox{ u.s.c.}} \leq \onebld_{( \, 0,\infty )}\left(\min_{ 1 \leq i \leq j-1}h_{j,i}(\xi,\boldsymbol{\beta})\right) \times  
\onebld_{[ \, 0,\infty )}\left(\min_{ j+1 \leq i \leq J } h_{j,i}(\xi,\boldsymbol{\beta})  \right) \leq  \underbrace{\onebld_{(\, 0,\infty )} (h^{2}_j(\xi, \boldsymbol{\beta};\varepsilon))}_{\mbox{ l.s.c. }},\] 
with
\[
\begin{array}{ll}
& h^{1}_j(\xi, \boldsymbol{\beta};\varepsilon)  \triangleq \min \left\{  \displaystyle{\min_{ 1 \leq i \leq j-1}h_{j,i}(\xi,\boldsymbol{\beta}) - \varepsilon, \min_{ j+1 \leq i \leq J } h_{j,i}(\xi,\boldsymbol{\beta})} \right\} ,  \\[0.15in]
& h^{2}_j(\xi, \boldsymbol{\beta};\varepsilon)  \triangleq \min \left\{  \displaystyle{\min_{ 1 \leq i \leq j-1}h_{j,i}(\xi,\boldsymbol{\beta}) , \min_{ j+1 \leq i \leq J } h_{j,i}(\xi,\boldsymbol{\beta})+ \varepsilon} \right\}. 
\end{array}
\]
Correspondingly, $\Phi_m(\boldsymbol{\beta})$ for $m \in [N]$ can be  approximated by the following two functions 
 \[
 \begin{array}{ll}
\Phi^1_m(\boldsymbol{\beta}; \varepsilon) \triangleq  \displaystyle{\sum_{s=m}^N} \left(2(N-s)+1 - 2 (C_{s})^2 \onebld_{[0, \infty)}(h^1_{D_s}(X^s, \boldsymbol{\beta}; \varepsilon))  \right), \quad 
\mbox{(l.s.c. function)}  \\[0.1in]
\Phi^2_m(\boldsymbol{\beta}; \varepsilon) \triangleq  \displaystyle{\sum_{s=m}^N} \left(2(N-s)+1 - 2 (C_{s})^2 \onebld_{(0, \infty)}(h^2_{D_s}(X^s, \boldsymbol{\beta}; \varepsilon))  \right), \quad 
\mbox{(u.s.c. function)} 
\end{array} \]  
such that $\Phi^2_m(\boldsymbol{\beta}; \varepsilon) \leq \Phi_m(\boldsymbol{\beta}) \leq \Phi^1_m(\boldsymbol{\beta}; \varepsilon)$. 
% for any $\varepsilon \geq 0$. 
Letting $\Phi_{N+1}^1(\boldsymbol{
\beta};\varepsilon) = \Phi_{N+1}^2(\boldsymbol{
\beta};\varepsilon) =0 $, we define 
\begin{equation}
\label{eq:obj_upper_simplified_beta}
{\small
\begin{array}{ll} 
& \wh \Psi_{\varepsilon}(\boldsymbol{\beta})  \triangleq \displaystyle{ \sum_{s=1}^N} \sum_{j=1}^J \, \wh{\mu}(X^s,j) \cdot \onebld_{[ \, 0,\infty )} (h^{1}_j(X^s, \boldsymbol{\beta};\varepsilon))   \\[0.2in]
 & + \displaystyle{ \sum_{s=1}^{m^*-1}}   C_{s} [ Y_{s} -\wh{\mu}(X^{s}, g(X^{s}) ) ]_+    \times    \,\onebld_{[\, 0,\infty )}\left(h^{1}_{D_{s}}(X^{s}, \boldsymbol{\beta};\varepsilon)\right)   \\[0.2in]
  & - \displaystyle{ \sum_{s=1}^{m^*-1}}   C_{s} [ Y_{s} -\wh{\mu}(X^{s}, g(X^{s}) ) ]_-    \times   \,\onebld_{(\, 0,\infty )}\left(h^{2}_{D_{s}}(X^{s}, \boldsymbol{\beta};\varepsilon)\right)   
  \\[0.2in]
     & + \displaystyle{ \sum_{s=m^*}^{N}}  C_{s} [ Y_{s} -\wh{\mu}(X^{s}, g(X^{s}) )  ]_+    \times  \onebld_{[ \, 0,\infty )} (h^{1}_{D_{s}}(X^{s}, \boldsymbol{\beta};\varepsilon)) \times    \displaystyle{\onebld_{[0,\infty)}} \left(  
\min_{s  \geq t \geq m^*}  \Phi^2_{t}(\boldsymbol{\beta}; \varepsilon) - \displaystyle{\min_{N+1 \geq m  \geq s+1}}\Phi^1_{m}(\boldsymbol{\beta}; \varepsilon)  \right)    \\[0.2in]
     & - \displaystyle{ \sum_{s=m^*}^{N}}  C_{s} [ Y_{s} -\wh{\mu}(X^{s}, g(X^{s}) )  ]_- \times  \onebld_{(\, 0,\infty )} (h^{2}_{D_{s}}(X^{s}, \boldsymbol{\beta};\varepsilon)) \times    \onebld_{(0,\infty)} \left( 
\min_{s  \geq t \geq m^*}  \Phi^1_{t}(\boldsymbol{\beta}; \varepsilon) - \displaystyle{\min_{N+1 \geq m  \geq s+1}}\Phi^2_{m}(\boldsymbol{\beta}; \varepsilon)  \right)  
\end{array}}%
\end{equation}

Thus, we propose to employ the following 
u.s.c.\ maximization problem:
\begin{equation}
\label{eq:CDR_policy_approx_optimization}
\max_{\boldsymbol{\beta} \in \mathbb R^p} 
\quad 
\wh{\Psi}_\varepsilon(\boldsymbol{\beta}),
\end{equation}
as the computational surrogate for the 
the policy optimization problem 
\eqref{eq:CDR_policy_optimization}. It is not difficult to see that 
$\wh \Psi_\varepsilon(\boldsymbol{\beta})$
is a u.s.c. function
lower bounding the function 
$\wh{\Psi}_{\rm HSC}(\boldsymbol \beta)$. 

In Proposition \ref{pr:locmax_AHSCOP_approximation}, we further give the conditions under which $\wh \Psi_\varepsilon(\boldsymbol{\beta})$ takes the same value as 
$\wh{\Psi}_{\rm HSC}(\boldsymbol \beta)$ in the neighborhood of a given point $\widebar{ \boldsymbol \beta}$, which yields that the maximization problem \eqref{eq:CDR_policy_approx_optimization} is equivalent to  the policy optimization problem \eqref{eq:CDR_policy_optimization}  in terms of local optimality.

\begin{proposition} \label{pr:locmax_AHSCOP_approximation} \rm 
For any given vector $\overline{\boldsymbol{\beta}} \in \mathbb R^p$, and any index $j \in [J]$,
\gap

\noindent \textbf{(A)} for  any $\xi \in \cal{X}$ with $\min_{ 1 \leq i \leq j-1}h_{j,i}(\xi,\overline{\boldsymbol{\beta}})\neq 0$, then there exist a scalar
$\bar \varepsilon_1 > 0$ and a neighborhood ${\cal N}_1$ of $\overline{\boldsymbol{\beta}}$
such that for any $\boldsymbol{\beta} \in {\cal N}_1$ and $0 < \varepsilon \leq \bar{\varepsilon}_1$, 
\[\begin{array}{ll}
{\onebld_{[ \, 0,\infty )} (h^{1}_j(\xi, \boldsymbol{\beta};\varepsilon))} & = \onebld_{( \, 0,\infty )}\left(\min_{ 1 \leq i \leq j-1}h_{j,i}(\xi,\boldsymbol{\beta})\right) \times  
\onebld_{[ \, 0,\infty )}\left(\min_{ j+1 \leq i \leq J } h_{j,i}(\xi,\boldsymbol{\beta})  \right)\\[0.1in]
& = \onebld_{( \, 0,\infty )}\left(\min_{ 1 \leq i \leq j-1}h_{j,i}(\xi,\overline{\boldsymbol{\beta}})\right) \times  
\onebld_{[ \, 0,\infty )}\left(\min_{ j+1 \leq i \leq J } h_{j,i}(\xi,\overline{\boldsymbol{\beta}})  \right)
\end{array}\]
Similar statement holds for $h^{2}_j(\xi, \boldsymbol{\beta};\varepsilon)$ for any $0 < \varepsilon \leq \bar{\varepsilon}_2$ with $\min_{ j+1 \leq i \leq J } h_{j,i}(\xi,\overline{\boldsymbol{\beta}})  \neq 0$ and some $\bar{\varepsilon}_2 >0$.

\noindent \textbf{(B)}  Let $\bar{\varepsilon}=\min\{\bar{\varepsilon}_1, \bar{\varepsilon}_2\}$. If $\overline{\boldsymbol  \beta}$ is a local maximizer to the optimization problem  \eqref{eq:CDR_policy_optimization}, there exits $\bar{\varepsilon}$ such that  $\overline{\boldsymbol \beta}$ is a local maximizer of \eqref{eq:CDR_policy_approx_optimization}  
for every $\varepsilon \in ( \, 0,\bar{\varepsilon} \, ]$.

\noindent 
\textbf{(C)}  Conversely, suppose that   $\overline{\boldsymbol{\beta}}$ satisfies the following local sign invariance
property: there exist  a neighborhood $\mathcal N$ of $\overline{\boldsymbol  \beta}$ such that for any ${\boldsymbol  \beta} \in \mathcal N$, 
\[ \begin{array}{ll}
 \min_{ 1 \leq i \leq j-1}h_{j,i}(X^s,\boldsymbol{\beta}) \leq 0, \forall j \times s \in \{(\ell,t): \ell \in [J],  \min_{ 1 \leq i \leq \ell-1}h_{\ell,i}(X^t,\overline{\boldsymbol{\beta}}) =0\} \\[0.1in]
  \min_{ D_s+1 \leq i \leq J}h_{D_s,i}(X^s,\boldsymbol{\beta}) \geq 0, \forall s \mbox{ satisfying } \min_{  D_s+1 \leq i \leq J}h_{D_s,i}(X^s,\overline{\boldsymbol{\beta}}) =0. 
 \end{array}\]  If $\overline{\boldsymbol{\beta}}$ is a local maximizer of \eqref{eq:CDR_policy_approx_optimization}, then $\overline{\boldsymbol{\beta}}$ is a local maximizer of  \eqref{eq:CDR_policy_optimization}.
\end{proposition}

\section{Algorithms for Policy Optimization}\label{sec:computation}
We now discuss how to algorithmically compute the learned policy from $\Gcal_{\rm LPC}$ by solving the surrogate policy optimization problem~\eqref{eq:CDR_policy_approx_optimization}. First, we discuss a full MIP formulation (given in Equation \eqref{eq:full_mip}), which can be conveniently solved using off-the-shelf MIP solvers. However, the full MIP becomes computationally inefficient as the problem scale grows. To address this limitation, we then adapt the PIP method, initially proposed in \cite{fang2025treatment} to our setting (presented in Algorithm \ref{alg: PIP}). The PIP algorithm aims to solve a reduced-scale version of the full MIP in Equation~\eqref{eq:full_mip},providing a computational enhancement. Finally, we embed the PIP algorithm, which serves as the computational workhorse, into the policy learning pipeline to compute the learned policy from $\mathcal{G}_{\mathrm{LPC}}$.The proposed policy learning algorithm is summarized in Algorithm \ref{alg: Policy Learning}.

\subsection{Computation via Mixed Integer Programming}
Below, we reformulate problem~\eqref{eq:CDR_policy_approx_optimization}
as a MIP by introducing binary variables $z^{1}_{j,s}, z^{2}_{j,s}$ to represent the Heaviside functions $\onebld_{[ \, 0,\infty )} (h^{1}_j(X^s, \boldsymbol{\beta};\varepsilon))$ and $\onebld_{( \, 0,\infty )} (h^{2}_j(X^s, \boldsymbol{\beta};\varepsilon)) $, respectively. In addition, we use binary variables $w_s^1$ and $w_s^2$ to represent $\displaystyle{\onebld_{[0,\infty)}}   
\left( \min_{s  \geq t \geq m^*}  \Phi^2_{t}(\boldsymbol{\beta}; \varepsilon) - \displaystyle{\min_{N+1 \geq m  \geq s+1}}\Phi^1_{m}(\boldsymbol{\beta}; \varepsilon) \right)$ and $\onebld_{(\, 0,\infty )} (h^{2}_{D_{s}}(X^{s}, \boldsymbol{\beta};\varepsilon)) \times    \onebld_{(0,\infty)} \left( 
\min_{s  \geq t \geq m^*}  \Phi^1_{t}(\boldsymbol{\beta}; \varepsilon) - \displaystyle{\min_{N+1 \geq m  \geq s+1}}\Phi^2_{m}(\boldsymbol{\beta}; \varepsilon)  \right)$, respectively. 
\begin{equation}\label{eq:full_mip}
{\small
\begin{array}{ll}
\displaystyle{ \underset{\boldsymbol{\beta} \in \mathbb R^p}{\mbox{maximize}}}  & \, \displaystyle{ \sum_{s=1}^N} \sum_{j=1}^J \, \wh{\mu}(X^s,j) z^{1}_{j,s} 
+ \displaystyle{ \sum_{s=1}^{m^*-1}}   C_{s} [ Y_{s} -\wh{\mu} (X^{s}, g(X^{s}) ) ]_+    \cdot    z_{D_s,s}^1 - \displaystyle{ \sum_{s=1}^{m^*-1}}   C_{s} [ Y_{s} -
\wh{\mu}(X^{s}, g(X^{s}) ) ]_- z_{D_s,s}^2   
\\[0.2in] 
& \displaystyle{ \sum_{s=m^*}^{N}} C_{s} \left( \begin{array}{ll} [ Y_{s} -
\wh{\mu}(X^{s}, g(X^{s}) )  ]_+  \min\{ z_{D_s,s}^{1}, w^1_{s}\}\\[0.1in]  -   
[ Y_{s} - \wh{\mu}(X^{s}, g(X^{s}) )]_- 
\min\{ z_{D_s,s}^{2}, w^2_{s}\} 
\end{array} \right) \\ [0.3in]
\mbox{subject to} & \, h^{1}_j(X^s, \boldsymbol{\beta};\varepsilon) \geq \underline{B}(1-z^1_{j,s}), \quad  \forall j\in [J], s\in [N]\\ [0.1in] 
& \, h^{2}_{D_s}(X^s, \boldsymbol{\beta};\varepsilon) \leq \overline{B}\, z^2_{D_s,s}, \quad  \forall  s\in [N]\\ [0.1in]
& \,      
 \left( \begin{array}{ll} \displaystyle{\min_{s \geq t \geq m^*}}  \left( \displaystyle{\sum_{k=t}^N} \left(2(N-k)+1 - 2   (C_{k})^2 
 z_{j_{k},k}^2\, \right) \right)\\ [0.2in]
 -  \min \left\{ \displaystyle{\min_{N \geq m  \geq s+1}}   \left( \displaystyle{\sum_{k=m}^N} \left(2(N-k)+1 - 2   (C_{k})^2 z_{j_{k},k}^1\, \right) \right), 0 \right\} \end{array}  \right) \geq \underline{B} (1- w_s^1), \, \forall  N \geq s  \geq m^* \\ [0.5in]
& \,        
 \left(  \begin{array}{ll} \displaystyle{\min_{s \geq t \geq m^*}}  \left( \displaystyle{\sum_{k=t}^N} \left(2(N-k)+1 - 2   (C_{k})^2 \cdot z_{D_k,k}^1\, \right) \right)\\ [0.2in]
 - \min \left\{ \displaystyle{\min_{N \geq m  \geq s+1}} \left( \displaystyle{\sum_{k=m}^N} \left(2(N-k)+1 - 2   (C_{k})^2 \cdot z_{D_k,k}^2\, \right) \right), 0 \right\} \end{array}  \right) \leq \overline{B}  w_s^2, \, \forall  N \geq s  \geq m^* \\ [0.5in]
 &  \,   z_{j,s}^1 \in \{0,1\}, \quad \forall j\in [J], s \in [N]\\[0.1in]
 & \,  z_{D_s,s}^2 \in \{0,1\}, \quad \forall  s \in [N]\\[0.1in]
& \, w_{s}^1,w_{s}^2 \in \{0,1\}, \quad \forall   N \geq s  \geq m^*\\[0.1in] 
& \, w_{s}^1 \geq w_{s'}^1, \quad \forall s' \geq s,  \forall s ', s\in [N].
\end{array}}%
\end{equation}

A key disadvantage of the MIP-based method is its inability to solve large-scale policy optimization problems. As seen in the formulation above, the total number of binary variables is in the order of the sample size $N$ and the number of treatments $J$. This creates an obstacle for achieving satisfactory learning performance, because statistical efficiency typically requires a large amount of data. As such, the MIP-based method may become computationally infeasible in such regimes. Motivated by this issue, we develop an algorithm grounded in the PIP method and specifically tailored to the structural characteristics of our surrogate policy optimization problem~\eqref{eq:CDR_policy_approx_optimization}.

\subsection{The Progressive Integer Programming Method}\label{sec:pip_formulation}
We now introduce our PIP method for solving the surrogate policy optimization problem \eqref{eq:CDR_policy_approx_optimization}. First, we construct the following index sets, used to identify the signs of $h_j^1(X^s, \widebar{\bm{\beta}}; \varepsilon)$ and $h_{D_s}^2(X^s, \widebar{\bm{\beta}}; \varepsilon)$, for a given $\widebar{\bm{\beta}}$ feasible to problem \eqref{eq:CDR_policy_approx_optimization}. Let $\delta_1, \delta_2 >0$ be two given nonnegative (small) scalars, define the complementary index sets:
\[
\begin{array}{ll}
\mathcal J_{1;<}^{-\delta_1}(\widebar{\boldsymbol{\beta}}) = \{(s,j) \in [N]\times[J]: h_j^1(X^s, \widebar{\boldsymbol{\beta}}; \varepsilon)< -\delta_1\}, \\[0.1in] 
\mathcal J_{1;>}^{+\delta_1}(\widebar{\boldsymbol{\beta}}) = \{(s,j)\in [N]\times[J]: h_j^1(X^s, \widebar{\boldsymbol{\beta}}; \varepsilon) > \delta_1\},\\[0.1in] 
  \mathcal J_{1;0}^{\pm\delta_1}(\widebar{\boldsymbol{\beta}}) = \{(s,j)\in [N]\times[J]:  -\delta_1 \leq   h_j^1(X^s, \widebar{\boldsymbol{\beta}};\varepsilon) \leq  \delta_1\},  \\[0.1in] 
  \mathcal J_{2;<}^{-\delta_2}(\widebar{\boldsymbol{\beta}}) = \{s \in [N]: h_{D_s}^2(X^s, \widebar{\boldsymbol{\beta}}; \varepsilon)< -\delta_2\}, \\[0.1in] 
  \mathcal J_{2;>}^{+\delta_2}(\widebar{\boldsymbol{\beta}}) = \{s\in [N]: h_{D_s}^2(X^s, \widebar{\boldsymbol{\beta}}; \varepsilon) > \delta_2\},\\[0.1in] 
  \mathcal J_{2;0}^{\pm\delta_2}(\widebar{\boldsymbol{\beta}}) = \{s\in [N]:  -\delta_2 \leq   h_{D_s}^2(X^s, \widebar{\boldsymbol{\beta}};\varepsilon)\leq \delta_2\}.
\end{array}
\]
For small $\delta_1$ and $\delta_2$, one can expect that the sets $\mathcal{J}_{1;0}^{\pm\delta_1}(\widebar{\boldsymbol{\beta}})$ and $\mathcal{J}_{2;0}^{\pm\delta_2}(\widebar{\boldsymbol{\beta}})$ (referred to as the ``in-between" sets) contains only a small number of elements. Therefore, based on these index sets, we can construct a subproblem of \eqref{eq:CDR_policy_approx_optimization} with a smaller scale, by fixing values of a portion of the Heaviside functions 
$\onebld_{[0,\infty)}(h_j^1(X^s, \widebar{\boldsymbol{\beta}}; \varepsilon))$ and 
$\onebld_{[0,\infty)}(h_{D_s}^2(X^s, \widebar{\boldsymbol{\beta}}; \varepsilon))$. Let $\boldsymbol{\delta} = [\delta_1, \delta_2]$, we state such a partial Heaviside formulation in terms of $\boldsymbol{\delta}$-modified 
objective: given a vector $\widebar{\boldsymbol{\beta}}$ feasible to problem \eqref{eq:CDR_policy_approx_optimization},

{\small
\[ % begin{equation}
% \label{eq:partial_obj_upper_beta}
\begin{array}{ll} 
& \wh{\Psi}_{\varepsilon}^{\, \boldsymbol{\delta}}(\boldsymbol{\beta})  = \displaystyle{ \sum_{(s,j) \in \mathcal J_{1;0}^{\pm\delta_1}(\widebar{\boldsymbol{\beta}})  }}  \, \wh{\mu}(X^s,j) \cdot \onebld_{[ \, 0,\infty )} (h^{1}_j(X^s, \boldsymbol{\beta};\varepsilon))  + \displaystyle{ \sum_{(s,j) \in \mathcal J_{1;>}^{+\delta_1}(\widebar{\boldsymbol{\beta}})  }}  \, \wh{\mu}(X^s,j)    \\[0.2in]
  &  + \displaystyle{ \sum_{\{s \in [m^*-1]: \, (s,D_s) \in \mathcal J_{1;0}^{\pm\delta_1}(\widebar{\boldsymbol{\beta}})\}}}   C_{s} [ Y_{s} -
  \wh{\mu}(X^{s}, g(X^{s}) ) ]_+ \cdot \,\onebld_{[\, 0,\infty )}\left(h^{1}_{D_{s}}(X^{s}, \boldsymbol{\beta};\varepsilon)\right)   +  \displaystyle{ \sum_{\{s \in [m^*-1]: \, (s,D_s) \in \mathcal J_{1;>}^{+\delta_1}(\widebar{\boldsymbol{\beta}})\}}}   C_{s} [ Y_{s} -\wh{\mu}(X^{s}, g(X^{s}) ) ]_+     
  \\[0.2in]
  & - \displaystyle{ \sum_{\{s \in [m^*-1]: \, s  \in \mathcal J_{2;0}^{\pm\delta_2}(\widebar{\boldsymbol{\beta}})\}}} C_{s} [ Y_{s} -
  \wh{\mu}(X^{s}, g(X^{s}) ) ]_- \cdot \,\onebld_{(\, 0,\infty )}\left(h^{2}_{D_{s}}(X^{s}, \boldsymbol{\beta};\varepsilon)\right) - \displaystyle{ \sum_{\{s \in [m^*-1]: \, s  \in \mathcal J_{2;>}^{+\delta_2}(\widebar{\boldsymbol{\beta}})\}}} C_{s} [ Y_{s} -
  \wh{\mu}(X^{s}, g(X^{s}) ) ]_-   \\[0.2in]
     & + \displaystyle{ \sum_{\{s:\, s \geq m^*, \, (s, D_s) \in \mathcal J_{1;0}^{\pm\delta_1}(\widebar{\boldsymbol{\beta}}) \}}} C_{s} [ Y_{s} -
     \wh{\mu}(X^{s}, g(X^{s}) )  ]_+ \cdot  \onebld_{[ \, 0,\infty )} (h^{1}_{D_{s}}(X^{s}, \boldsymbol{\beta};\varepsilon)) \cdot \displaystyle{\onebld_{[0,\infty)}} \left(  
\min_{s  \geq t \geq m^*}  \Phi^2_{t}(\boldsymbol{\beta}; \varepsilon) - \displaystyle{\min_{N+1 \geq m  \geq s+1}}\Phi^1_{m}(\boldsymbol{\beta}; \varepsilon)  \right)    \\[0.2in]
    & + \displaystyle{ \sum_{\{s:\, s \geq m^*, \, (s,D_s) \in \mathcal J_{1;>}^{\delta_1} (\widebar{\boldsymbol{\beta}})\}}} C_{s} [ Y_{s} -
    \wh{\mu}(X^{s}, g(X^{s}) )  ]_+    \cdot     \displaystyle{\onebld_{[0,\infty)}} \left(  
\min_{s  \geq t \geq m^*}  \Phi^2_{t}(\boldsymbol{\beta}; \varepsilon) - \displaystyle{\min_{N+1 \geq m  \geq s+1}}\Phi^1_{m}(\boldsymbol{\beta}; \varepsilon)  \right)    \\[0.2in]
& - \displaystyle{ \sum_{\{ s: \, s\geq m^*, s\in \mathcal J_{2;0}^{\pm\delta_2}(\widebar{\boldsymbol{\beta}}) \}}}   C_{s} [ Y_{s} -\wh{\mu}(X^{s}, g(X^{s}) )  ]_- \cdot  \onebld_{(\, 0,\infty )} (h^{2}_{D_{s}}(X^{s}, \boldsymbol{\beta};\varepsilon)) \cdot  \onebld_{(0,\infty)} \left( 
\min_{s  \geq t \geq m^*}  \Phi^1_{t}(\boldsymbol{\beta}; \varepsilon) - \displaystyle{\min_{N+1 \geq m  \geq s+1}}\Phi^2_{m}(\boldsymbol{\beta}; \varepsilon)  \right)  \\[0.2in]
& - \displaystyle{ \sum_{\{ s: \, s\geq m^*, s\in \mathcal J_{2;>}^{\delta_2}(\widebar{\boldsymbol{\beta}}) \}}}   
C_{s} [ Y_{s} -\wh{\mu}(X^{s}, g(X^{s}) )  ]_- \cdot    \onebld_{(0,\infty)} \left( 
\min_{s  \geq t \geq m^*}  \Phi^1_{t}(\boldsymbol{\beta}; \varepsilon) - \displaystyle{\min_{N+1 \geq m  \geq s+1}}\Phi^2_{m}(\boldsymbol{\beta}; \varepsilon)  \right).  
\end{array}
\]} % end{equation}}

Building on the objective function $\wh{\Psi}_{\varepsilon}^{\boldsymbol{\delta}}(\boldsymbol{\beta})$, we construct a mixed integer subproblem, with a key modification from Equation~\eqref{eq:full_mip}: the binary variables $z^1_{j,s}$ (resp. $z^2_{j,s}$) corresponding to indices not in $\mathcal{J}_{1;0}^{\pm\delta_1}(\widebar{\boldsymbol{\beta}})$ (resp.\ $\mathcal{J}_{2;0}^{\pm\delta_2}(\widebar{\boldsymbol{\beta}})$) are fixed to constants $0$ or $1$. The detailed formulation is provided in Section \ref{sec:restricted_ip} of the online appendix; relative to the full MIP, it involves fewer active binary variables.

Our PIP algorithm begins with a feasible solution $\widebar{\boldsymbol{\beta}}$ of problem 
\eqref{eq:CDR_policy_approx_optimization} and then iteratively solves a sequence of subproblems obtained by maximizing $\wh{\Psi}_{\varepsilon}^{\boldsymbol{\delta}}(\boldsymbol{\beta})$. The parameters $\delta_1$ and $\delta_2$ are updated adaptively during the execution of the algorithm, dynamically controlling the number of binary variables in each subproblem. In particular, at each iteration $\nu+1$, the index sets that determine the binary variables in the next subproblem are defined based on the feasible solution $\boldsymbol{\beta}^{\nu}$ obtained in iteration $\nu$. This design ensures that $\boldsymbol{\beta}^{\nu}$ remains feasible for the subsequent subproblem and that the objective values are nondecreasing across iterations (see Proposition 5 of \cite{fang2025treatment} for a formal characterization of this property). The improvement in objective values can be achieved through either expanding or shrinking the set of binary variables. Incorporating more binary variables tightens the MIP subproblems, while decreasing the number of binary variables reduces computational effort. Accordingly, we specify the following update rule: if the objective value does not improve from the previous iteration, we expand the set of binary variables; if the objective value improves, we shrink this set. The PIP algorithm terminates when the objective value remains unchanged after several consecutive expansions of the binary variable set. 

We summarize the PIP algorithm applied to our policy optimization problem in Algorithm~\ref{alg: PIP}, which is developed based on the conceptual ideas discussed above. The detailed settings of the input parameters will be provided when we practically implement this algorithm in the numerical experiment part (Section \ref{sec:numerical}). A slight difference from the conceptual introduction is that we construct the in-between sets using an asymmetric band, where the left and right cutoffs are mapped to quantiles of the corresponding index-set elements. This prevents the algorithm from accidentally including or excluding all binary variables $z_{j,s}^1$ and $z_{j,s}^2$, since the empirical distributions of $h_j^1(X^s, \widebar{\boldsymbol{\beta}}; \varepsilon)$ and $h_{D_s}^2(X^s, \widebar{\bm{\beta}}; \varepsilon)$ might not be centered at zero and might be skewed.

The key advantage of the PIP algorithm over the MIP method is improved computational efficiency when the policy optimization problem involves a large sample size $N$, which will be shown empirically in Section \ref{sec:numerical}.

\begin{algorithm}[H]
\begin{algorithmic}[1]
\caption{The Progressive Integer Programming (PIP) Method}
\Require Initial feasible solution $\bm{\beta}^0$, initial ratio $r^0$, maximum ratio $r_{\text{max}}$, minimum ratio $r_{\text{min}}$, expansion ratio $\Delta_{\text{expand}}$, shrinkage ratio $\Delta_{\text{shrink}}$, maximum number of iterations $\nu_{\max}$,  maximum number of iterations with unchanged objective value $\Tilde{\nu}_{\max}$.\\
$\nu \leftarrow 0$, $\Tilde{\nu} \leftarrow 0$.\\
\textbf{While} $\nu < \nu_{\max}$ and $\Tilde{\nu} < \Tilde{\nu}_{\max}$ \textbf{do}:\\
~~ $\Omega_1^{+}, \Omega_1^{-},\Omega_2^{+}, \Omega_2^{-} \leftarrow \emptyset$;\\
~~ \textbf{For} $(s,j)$ in $[S] \times [J]$ \textbf{do}:\\
~~~~ \textbf{if} $h_j^1(X^s, \bm{\beta}^\nu; \varepsilon) > 0$, \textbf{then} $\Omega_1^{+} \leftarrow \Omega_1^{+} \cup \{h_j^1(X^s, \bm{\beta}^\nu; \varepsilon)\}$;\\
~~~~ \textbf{else if} $h_j^1(X^s, \bm{\beta}^\nu; \varepsilon) < 0$, \textbf{then} $\Omega_1^{-} \leftarrow \Omega_1^{-} \cup \{h_j^1(X^s, \bm{\beta}^\nu; \varepsilon)\}$;\\
~~~~ \textbf{else} randomly choose $\Omega_1 \in \{\Omega_1^{+}, \Omega_1^{-}\}$ and $\Omega_1 \leftarrow \Omega_1 \cup \{h_j^1(X^s, \bm{\beta}^\nu; \varepsilon)\}$;\\
~~ \textbf{end}\\
~~ \textbf{For} $s$ in $[S]$ \textbf{do}:\\
~~~~ \textbf{if} $h_{D_s}^2(X^s, \bm{\beta}^\nu; \varepsilon) > 0$, \textbf{then} $\Omega_2^{+} \leftarrow \Omega_2^{+} \cup \{h_{D_s}^2(X^s, \bm{\beta}^\nu; \varepsilon)\}$;\\
~~~~ \textbf{else if} $h_{D_s}^2(X^s, \bm{\beta}^\nu; \varepsilon) < 0$, \textbf{then} $\Omega_2^{-} \leftarrow \Omega_2^{-} \cup \{h_{D_s}^2(X^s, \bm{\beta}^\nu; \varepsilon)\}$;\\
~~~~ \textbf{else} randomly choose $\Omega_2 \in \{\Omega_2^{+}, \Omega_2^{-}\}$ and $\Omega_2 \leftarrow \Omega_2 \cup \{h_{D_s}^2(X^s, \bm{\beta}^\nu; \varepsilon)\}$;\\
~~ \textbf{end}\\
~~ $\delta_{1,+}^\nu \leftarrow \Omega_1^{+}.\mathrm{lower\_quantile}(r^{\nu})$, $\delta_{1,-}^\nu \leftarrow -\Omega_1^{-}.\mathrm{upper\_quantile}(r^{\nu})$; \\
~~ $\delta_{2,+}^\nu \leftarrow \Omega_2^{+}.\mathrm{lower\_quantile}(r^{\nu})$, $\delta_{2,-}^\nu \leftarrow -\Omega_2^{-}.\mathrm{upper\_quantile}(r^{\nu})$;\\
~~ \textbf{Form} $\mathcal J_{1;<}^{-\delta_{1,-}^\nu}(\boldsymbol{\beta}^\nu)$, $\mathcal J_{1;>}^{+\delta_{1,+}^\nu}(\boldsymbol{\beta}^\nu)$, $\mathcal J_{1;0}^{-\delta_{1,-}^\nu, +\delta_{1,+}^\nu}(\boldsymbol{\beta}^\nu)$, $\mathcal J_{2;<}^{-\delta_{2,-}^\nu}(\boldsymbol{\beta}^\nu)$, $\mathcal J_{2;>}^{+\delta_{2,+}^\nu}(\boldsymbol{\beta}^\nu)$, $\mathcal J_{2;0}^{-\delta_{2,-}^\nu, +\delta_{2,+}^\nu}(\boldsymbol{\beta}^\nu)$\\
~~~~\textbf{Solve} a restricted MIP 
by maximizing $\wh{\Psi}_{\varepsilon}^{\, \delta}(\boldsymbol{\beta})$ 
(with the 
$\delta_{i,\pm}^{\nu}, i = 1,2$)
% \eqref{eq:partial_obj_upper_beta} 
to obtain a globally optimal solution 
$\bm{\beta}^{\nu+1}$ with the corresponding optimal objective value $\psi^{\nu+1}$;\\
~~~~\textbf{if} $\psi^{\nu+1} = \psi^{\nu}$ \textbf{then} $r^{\nu+1} \leftarrow \min\{r^{\nu} + \Delta_{\text{expand}}, r_{\max}\}$, $\Tilde{v} \leftarrow \Tilde{v} + 1$;\\
~~~~\textbf{else} $r^{\nu+1} \leftarrow \max\{r^{\nu} - \Delta_{\text{shrink}}, r_{\min}\}$, $\Tilde{\nu} \leftarrow 0$;\\
~~~~$\nu = \nu + 1$;\\
\textbf{end}
\label{alg: PIP}
\end{algorithmic}
\end{algorithm}

\subsection{A Practical Offline Policy Learning Algorithm: OCDRL}\label{sec:policy_learning_algo}
With the ability to efficiently solve problem \eqref{eq:CDR_policy_approx_optimization}, we can develop a practical offline policy learning algorithm, called Optimized Clipped Doubly Robust Learning (OCDRL), summarized in Algorithm~\ref{alg: Policy Learning}.
OCDRL consists of three main steps: (i) estimate the reward model $\mu(\bullet,\bullet)$ with machine learning methods and sample-splitting; \endnote{In the spirit of “orthogonalization,” (see, e.g., \citealt{foster2023orthogonal, zhou2023offline}), the purpose of estimating $\mu(\bullet,\bullet)$ in separate subsamples is to ensure that these nuisance components are independent of the policy value estimation step. This independence reduces the impact of their estimation errors on the target estimand, i.e., the value of a policy $g$.}
(ii) construct the OCDR estimator as in Definition \ref{def:OCDR}, where the clipping threshold is optimized to be $\wh{\tau}(g)$ 
as in Equation \eqref{eq:optimal_tau_simplified_pruning}; and
(iii) solve for the optimal policy within the policy class $\mathcal{G}_{\rm LPC}$. In particular, step (iii) leverages the structural properties of $\wh{\tau}(g)$ and the linear treatment rule to formulate the policy optimization problem into a tractable u.s.c. $\varepsilon$-approximation of an HSCOP, which is solved by the PIP method described in Algorithm \ref{alg: PIP}. 

\begin{algorithm}
\begin{algorithmic}[1]
\caption{Optimized Clipped Doubly Robust Learning: Invoking PIP}
\Require Data set $\cal{D} = \{(X^s, D_s, Y_s)\}_{s=1}^n$.\\
Split $\cal{D} $ into subset $\cal{D} _1$ and $\cal{D} _2$.\\
Fit plug-in estimators $\wh{\mu}(\bullet ,j)$ using data in $\cal{D}_1$ for all $j \in [J]$. \\
Use data in $\cal{D}_2$ to construct the OCDR estimator $\wh{\Psi}^{\mathrm{OCDR}}(g)$ as in Definition \ref{def:OCDR}.\\
Reformulate policy optimization 
problem as $\wh{\Psi}_\varepsilon(g)$ (problem \eqref{eq:CDR_policy_approx_optimization}) and call PIP in Algorithm \ref{alg: PIP} to solve for policy $\wh{g} \in \displaystyle{\argmax_{g \in \Gcal_{\rm LPC}}} \wh{\Psi}_{\varepsilon}^{\bm \delta}(g)$.
\Ensure $\wh{g} \in \Gcal_{\rm LPC}$, which can be explicitly represented by parameters $\wh{\bm{\beta}}^j \in \mathbb{R}^d$ for $j \in [J]$.
\label{alg: Policy Learning}
\end{algorithmic}
\end{algorithm}

\section{Statistical Properties}\label{sec:stats_performance}
In this section, we establish statistical guarantees for the proposed OCDRL algorithm (Algorithm~\ref{alg: Policy Learning}). To isolate the statistical properties from computational ones, our analysis treats the learned policy $\wh{g}$ as a global maximizer of the policy optimization problem \eqref{eq:CDR_policy_optimization}. This corresponds to an oracle implementation of the policy optimization step. 
In practice, problem \eqref{eq:CDR_policy_optimization} is solved approximately via the surrogate problem \eqref{eq:CDR_policy_approx_optimization}, with guarantees on local optimality (cf. Proposition~\ref{pr:locmax_AHSCOP_approximation}); we investigate the computational performance empirically in Section~\ref{sec:numerical}, showing that the discrepancy between the achieved objective value and Gurobi’s best objective bound is small.

We note that a standard approach to control the suboptimality gap is by bounding the worst-case policy estimation error (see, e.g., \cite{zhan2024policy}), since
\begin{align*}
    \cal{L}(\wh{g},\Gcal_{\rm LPC}) = \Psi(g^*)-\Psi(\wh{g}) &\leq \Psi(g^*)-
    \wh{\Psi}^{\rm OCDR}(g^*) + 
    \wh{\Psi}^{\rm OCDR}(g^*) - 
    \wh{\Psi}^{\rm OCDR}(\wh{g}) + 
    \wh{\Psi}^{\rm OCDR}(\wh{g}) - \Psi(\wh{g})\\
    &\leq \max_{g \in \mathcal{G}_{\rm LPC}} |\Psi(g) - \wh{\Psi}^{\rm OCDR}(g)|,
\end{align*}
where the second inequality uses the fact that $\wh{g}$ is the maximizer of 
$\wh{\Psi}^{\rm OCDR}(g)$. According to the definition of the OCDR estimator (Definition \ref{def:OCDR}), the worst-case estimation error $\displaystyle{
\max_{g \in \mathcal{G}_{\rm LPC}}}
|\Psi(g) - \wh{\Psi}^{\rm OCDR}(g)|$ is affected by the choice of the clipping threshold. We therefore establish an upper bound on the worst-case estimation error that characterizes its dependence on $\wh{\tau}(g)$. Our result justifies the necessity of optimizing the clipping threshold at 
$\wh{\tau}(g)$. The suboptimality of the proposed OCDRL algorithm is characterized as follows.

\begin{theorem}\label{thm:suboptimality_upper_bound}
For a fixed policy $g \in \Gcal_{\rm LPC}$, define
{\small
\begin{align*}
    \wh{B}(g, \tau(g)) &\triangleq \frac{M^2}{N^2} \left(\left(\sum_{s=1}^N \onebld\left\{\frac{1}{e(X^s, g(X^s))} > \tau(g)\right\}\right)^2 + \sum_{s=1}^N \frac{\onebld\{D_s = g(X^s)\}}{e(X^s, g(X^s))^2} \onebld\left\{\frac{1}{e(X^s, g(X^s))} \leq \tau(g) \right\}\right),\\
    \wt{B}(g, \tau(g)) &\triangleq \frac{M^2}{N^2} \left(\left(\sum_{s=1}^N \onebld\left\{\frac{1}{e(X^s, g(X^s))} > \tau(g)\right\}\right)^2 + \sum_{s=1}^N \frac{1}{e(X^s, g(X^s))} \onebld\left\{\frac{1}{e(X^s, g(X^s))} \leq \tau(g) \right\}\right),\\
    \Delta(g, \tau(g)) &\triangleq M^2 \left(\frac{\tau(g)}{N}\right)^{3/2},
\end{align*}}
where $\wh{\tau}(g)$ and $\wt{\tau}(g)$ denote the minimizers of $\wh{B}(g, \tau(g))$ and $\wt{B}(g, \tau(g))$, respectively.

Suppose $N \geq \frac{2\sigma^2}{\omega} \ln\left(\frac{2N}{\delta}\right)$, where $\omega>0$ and $\delta \in [0,1]$. With probability at least $1-\delta$, the suboptimality is upper bounded, such that
\begin{align*}
    \cal{L}(\wh{g}, \Gcal_{\rm LPC}) \leq \min\left\{2\upsilon \cdot \max\left\{ \wh{B}(g^\dagger, \wh{\tau}(g^\dagger))^{1/2}, \wt{B}(g^\dagger,\wt{\tau}(g^\dagger))^{1/2}+\omega^{1/2}, \Delta(g^\dagger, \wh{\tau}(g^\dagger))^{1/2}\right\}, M\right\},
\end{align*}
where $\upsilon \geq 12M \sqrt{ 2|J|(p+1) \log(N|J|^2) + 2\log(16/\delta)}$ and $g^\dagger = \argmax_{g \in \Gcal_{\rm LPC}} \left|\wh{\Psi}(g) - \Psi(g)\right|$.
\end{theorem}

In Theorem \ref{thm:suboptimality_upper_bound}, the suboptimality bound is dominated by three terms. First, $\wh{B}(g^\dagger, \wh{\tau}(g^\dagger))^{1/2}$ can be explained as the upper bound for the worst-case empirical policy estimation error. Notably, by fixing the policy at $g^\dagger$, it coincides with the objective function of Problem \eqref{program: optimal tau} up to a constant term. Therefore, the minimizer of Problem \eqref{program: optimal tau} $\wh{\tau}(g^\dagger)$ also minimizes $\wh{B}(g^\dagger, \tau(g^\dagger))$. Second, $\wt{B}(g^\dagger, \wt{\tau}(g^\dagger))^{1/2} + \omega^{1/2}$ upper bounds the worst-case expected policy estimation error $\wt{B}(g^\dagger, \wh{\tau}(g^\dagger))$ (conditional on $\{X^s\}_{s=1}^N$), where we highlight that $\wt{\tau}(g^\dagger)$ is the minimizer of $\wt{B}(g^\dagger, \tau(g^\dagger))$. Finally, $\Delta(g^\dagger,\wh{\tau}(g^\dagger))$ is the upper bound of the deviation between $\wh{B}(g^\dagger,\wh{\tau}(g^\dagger))$ and $\wt{B}(g^\dagger,\wh{\tau}(g^\dagger))$, penalizing an excessively large clipping threshold that would otherwise amplify both the discrepancy between empirical and population errors and the estimation variance. 

In general, Theorem \ref{thm:suboptimality_upper_bound} implies that optimizing the clipping threshold $\wh{\tau}(g^\dagger)$ not only minimizes the policy estimation error but also directly tightens the suboptimality gap in policy learning. Since our OCDR estimator is equipped with an optimal clipping threshold, the suboptimality gap of the proposed OCDRL algorithm can be effectively controlled.

To provide a more interpretable and comparable result, we state the following corollary under a setting where the overlap condition holds with sufficient (but mild) strength.

\begin{corollary}\label{cor:convergence_rate}
    For any $\delta \in (0, \exp(-1))$. Suppose $N \geq \frac{2\sigma^2}{\omega} \ln\left(\frac{2N}{\delta}\right)$, $\eta \geq \frac{1}{N}$, and $\nu = c \cdot M \sqrt{ 2|J|(p+1) \log(N|J|^2) + 2\log(16/\delta)}$ for some constant $c \geq 12\sqrt{2}$. Let $g^\dagger = \argmax_{g \in \Gcal_{\rm LPC}} \left|\wh{\Psi}(g) - \Psi(g)\right|$. We have that with probability at least $1-2\delta$,
    {\small
    \begin{align*}
        \cal{L}(\wh{g}, \cal{G}_{\rm LPC}) \leq 4c \cdot M \sqrt{|J|(p+1) \log(N|J|^2)} \cdot \left(\log(2/\delta)^{3/2} \sqrt{\frac{\tau^*(g^\dagger) \vee \wh{\tau}(g^\dagger)}{N}} + \log(2/\delta)^{1/2}\wh{P}(\wh{\tau}(g^\dagger))^{1/2} + \omega^{1/2}\right),
    \end{align*}}
where $\wh{P}(\wh{\tau}(g^\dagger)) = \frac{1}{N}\sum_{s=1}^N \onebld\left\{e(X^s,g^\dagger(X^s))^{-1} > \wh{\tau}(g^\dagger)\right\}$.
\end{corollary}

In the suboptimality upper bound presented in Corollary \ref{cor:convergence_rate}, the first term in the brackets, $\log(2/\delta)^{3/2} \sqrt{\frac{\tau^*(g^\dagger) \vee \wh{\tau}(g^\dagger)}{N}}$, upper bounds the impact of variance. Notably, this variance-related term achieves a sharper statistical rate than the classic result in \cite{kitagawa2018should}, which shows that $\cal{L}(\wh{g}, \Gcal_{\rm LPC}) \lesssim \sqrt{|J|(p+1)/\eta N}$. By contrast, the variance-related term in our bound is smaller, as it depends on $\tau^*(g^\dagger) \vee \wh{\tau}(g^\dagger)$, which is no larger than $1/\eta$. This immediately implies that weight clipping improves policy learning performance by mitigating the impact of high variance in policy estimation, thus preventing spuriously inflated policy value estimates. 

The second term in the brackets upper bounds the impact of bias induced by weight clipping and is proportional to the empirical probability that the IPS exceeds the clipping threshold, i.e.,$ \wh{P}(\wh{\tau}(g^\dagger))$. Note that $\wh{P}(\wh{\tau}(g^\dagger))$ is uniformly bounded by 1, implying that the impact of bias in OCDRL is strictly reduced compared with learning algorithms that use the DM estimator (for which $\wh{P}(\wh{\tau}(g^\dagger)) \equiv 1$, since $\wh{\tau}(g^\dagger)$ is always set to $1/\eta$). Moreover, because the clipping threshold $\wh{\tau}(g^{\dagger})$ is chosen by minimizing the empirical MSE, it naturally balances the trade-off between bias and variance. As a consequence, the contribution of the bias term is controlled and does not deteriorate the overall $O(N^{-1/2})$ convergence rate.

\section{Numerical Experiments}\label{sec:numerical}
In this section, we present numerical results from two synthetic experiments to evaluate both the computational and statistical performance of the proposed policy learning framework.\endnote{ All experiments are performed on a Lenovo laptop with a 1.70 GHz Intel Core Ultra 5 processor and 32GB of memory.} The first experiment demonstrates the computational efficiency of the PIP algorithm, and the second experiment is designed to illustrate the statistical superiority of the proposed OCDRL algorithm relative to benchmark algorithms. We also apply the proposed algorithm to a real-world data set.

\subsection{Assessment of Computational Performance}\label{sec:computational_performance}
\subsubsection{Setup}
The synthetic data set is represented by triples $\{(X^s, D_s, Y_s)\}_{s=1}^N$, generated by mimicking the data-generating process in \cite{fang2025treatment}. First, we sample $|\mathcal{X}|$ distinct vectors from a 20-dimensional unit cube $[0,1]^{20}$. Next, each $X^s$ is randomly drawn from these $|\mathcal{X}|$ vectors, which gives $N$ 20-dimensional vectors $\{(X^s_d)_{d=1}^{20}\}_{s=1}^N$. We consider different combinations of $(|\mathcal{X}|, N)$, specified as $(100, 500)$, $(200, 1000)$, $(300, 1500)$, and $(400,2000)$. The treatments $D_s$ are randomly sampled from $[J] = \{1,2,3,4\}$ following multinomial logit propensity scores, such that $e(X^s;\Theta) = \operatorname{softmax}(\Theta^\top X^s)$, where \(\Theta = [\theta_1,\ldots,\theta_J]\) and each $\theta_j \sim N(0,1)$. Each reward $Y_s = Y(D_s)$ is generated as $Y(D) = Y(1)\onebld\{D=1\} + Y(2)\onebld\{D=2\} + Y(3)\onebld\{D=3\} + Y(4)\onebld\{D=4\} + \varepsilon$, where each random $Y(j)$ is an exponential function (the detailed functional forms are provided in Section \ref{sec:omitted_exp_details} of the online appendix) and $\varepsilon \sim  \mathrm{Lognormal}(0,0.001)$.

For each $(|\mathcal{X}|, N)$ pair, we generate five synthetic data sets using different random seeds. Based on these data sets, we compare the computational performance of the PIP algorithm (Algorithm~\ref{alg: PIP}) with that of the MIP-based method in solving the policy optimization problem given in Equation~\eqref{eq:CDR_policy_approx_optimization}. Note that since the dimensionality of the covariates is high, we employ a regularizer to induce sparsity in parameter selection (See Equation \eqref{eq:regularized_policy_optimization} in the online appendix for the regularized policy optimization formulation.). 

For solving the MIP, we set the GUROBI time limit to 1.5 hours when $N < 1500$ and to 2 hours when $N \geq 1500$. We provide the detailed parameter settings of the PIP algorithm below:

\noindent $\bullet$ \textbf{Initial solutions:} To obtain the initial solution $\bm{\beta}^0$ at which PIP is initiated, we solve the LP relaxation of the full MIP problem formulated in Equation \eqref{eq:full_mip}. 

\noindent  $\bullet$ \textbf{Stopping criterion:} Our PIP algorithm terminates when either of the following conditions is met:
(i) the number of updates to the in-between interval reaches 15, or (ii) the objective value remains unchanged for three consecutive iterations. In addition, each PIP subproblem is solved with a time limit: 100 seconds for small-scale problems ($N < 1500$) and 120 seconds for large-scale problems ($N \geq 1500$).

\subsubsection{Results of Computational Performance}
The results are summarized in Table \ref{tab:computation}. We also report the best objective upper bound reported by the output log of GUROBI. If the MIP solver fails to output a solution within the given time limit, we mark the instance as ``infeasible". The results show that, for small-scale problems ($N = 500$), the rewards achieved by PIP are close to (though slightly lower than) those obtained by the MIP, but with significantly shorter computation time. When problems exceed a moderate scale ($N \geq 1000$), the advantage of PIP becomes dominant: it spends only less than $25\%$ of the computation time to obtain a higher reward compared to the MIP, which in some cases fails to find a feasible solution within two hours. Overall, this experiment highlights the necessity of employing an efficient computational method for policy optimization --- a role effectively fulfilled by the PIP algorithm, especially for large-scale problems.

\begin{table*}[t]
\centering
\caption{Comparisons Between PIP and MIP in terms of Computational Performance}
\label{tab:computation}
\footnotesize
\setlength{\tabcolsep}{4pt}
\renewcommand{\arraystretch}{1.08}

% Tight two-line cells
\newcommand{\hdr}[2]{\makecell[c]{#1\\[-0.35ex]#2}}
\newcommand{\obj}[2]{\makecell[c]{#1\\[-0.35ex](#2)}}
\newcommand{\objb}[2]{\makecell[c]{\textbf{#1}\\[-0.35ex]\textbf{(#2)}}}

\begin{threeparttable}
\begin{tabular}{@{}c c *{4}{c c c}@{}}
\toprule
& & \multicolumn{3}{c}{$(|\mathcal{X}|,N)=(100,500)$}
  & \multicolumn{3}{c}{$(|\mathcal{X}|,N)=(200,1000)$}
  & \multicolumn{3}{c}{$(|\mathcal{X}|,N)=(300,1500)$}
  & \multicolumn{3}{c}{$(|\mathcal{X}|,N)=(400,2000)$} \\
\cmidrule(lr){3-5}\cmidrule(lr){6-8}\cmidrule(lr){9-11}\cmidrule(lr){12-14}
\hdr{Seed}{} & \hdr{Method}{}
& \hdr{Best}{bound} & \hdr{Objective}{(reward)} & \hdr{Time}{}
& \hdr{Best}{bound} & \hdr{Objective}{(reward)} & \hdr{Time}{}
& \hdr{Best}{bound} & \hdr{Objective}{(reward)} & \hdr{Time}{}
& \hdr{Best}{bound} & \hdr{Objective}{(reward)} & \hdr{Time}{} \\
\midrule

\multirow{2}{*}{1}
& MIP & 27.357 & \obj{27.356}{27.389} & 1684s
      & 29.689 & \obj{28.470}{28.516} & 5401s
      & 29.472 & \obj{21.082}{21.166} & 7201s
      & 30.156 & \obj{23.764}{23.833} & 7201s \\
& PIP &      & \obj{26.302}{26.326} & \textbf{247s}
      &      & \objb{28.686}{28.731} & \textbf{1462s}
      &      & \objb{27.730}{27.767} & \textbf{1444s}
      &      & \objb{25.813}{25.863} & \textbf{1804s} \\

\addlinespace[0.6ex]
\multirow{2}{*}{2}
& MIP & 28.186 & \obj{28.186}{28.211} & 685s
      & 28.140 & \obj{21.786}{21.859} & 5401s
      & 27.139 & \obj{24.626}{24.705} & 7202s
      & 32.174 & \makecell[c]{infeas.} & 7200s \\
& PIP &      & \obj{27.240}{27.252} & \textbf{416s}
      &      & \objb{27.626}{27.659} & \textbf{750s}
      &      & \objb{26.351}{26.424} & \textbf{1443s}
      &      & \objb{27.646}{27.684} & \textbf{1779s} \\

\addlinespace[0.6ex]
\multirow{2}{*}{3}
& MIP & 28.693 & \obj{28.639}{28.693} & 5401s
      & 28.114 & \obj{22.755}{22.886} & 5401s
      & 27.091 & \obj{21.004}{21.071} & 7200s
      & 28.215 & \obj{21.210}{21.271} & 7202s \\
& PIP &     & \obj{28.611}{28.653} & \textbf{351s}
      &     & \objb{26.371}{26.443} & \textbf{1165s}
      &     & \objb{24.603}{24.640} & \textbf{1528s}
      &     & \objb{25.645}{25.678} & \textbf{1435s} \\

\addlinespace[0.6ex]
\multirow{2}{*}{4}
& MIP & 30.517 & \obj{30.516}{30.533} & 436s
      & 29.561 & \obj{17.554}{17.693} & 5401s
      & 28.283 & \makecell[c]{infeas.} & 7200s
      & 29.762 & \obj{27.092}{27.123} & 7201s \\
& PIP &     & \obj{30.509}{30.515} & \textbf{67s}
      &     & \objb{28.529}{28.596} & \textbf{741s}
      &     & \objb{24.791}{24.817} & \textbf{1415s}
      &     & \objb{27.608}{27.666} & \textbf{962s} \\

\addlinespace[0.6ex]
\multirow{2}{*}{5}
& MIP & 28.708 & \obj{28.542}{28.608} & 5401s
      & 28.800 & \obj{26.095}{26.176} & 5401s
      & 30.973 & \obj{30.336}{30.355} & 7201s
      & 31.750 & \obj{24.040}{24.075} & 7201s \\
& PIP &    & \obj{28.539}{28.583} & \textbf{675s}
      &    & \objb{26.298}{26.338} & \textbf{1401s}
      &    & \objb{30.710}{30.740} & \textbf{1575s}
      &    & \objb{26.778}{26.807} & \textbf{486s} \\

\bottomrule
\end{tabular}

\begin{tablenotes}[flushleft]
\footnotesize
\item \textit{Note.} ``Best bound'' refers to the best objective upper bound reported by Gurobi, ``Objective'' refers to the value of the full objective function (including the regularization term), and ``Reward'' corresponds exactly to the objective value without regularization. Objective values and rewards obtained by PIP are shown in bold if they are greater than those obtained by MIP, and the computation time for PIP is highlighted in bold when it is shorter than that of MIP. 
\end{tablenotes}
\end{threeparttable}
\end{table*}

\subsection{Assessment of Statistical Performance}
We now move on to the experiment designed to illustrate the statistical performance of our proposed OCDRL framework. In particular, we demonstrate that the OCDR estimator plays an important role in mitigating the high-variance issue and in reducing the suboptimality gap of the learned policy. In contrast, learning frameworks equipped with classic estimators such as DR or IPW are vulnerable to high variance and therefore suffer from inferior policy selection.

\subsubsection{Setup} We consider a policy learning problem with the goal of maximizing reward. We fix the number of treatments at $|J| = 3$ and the data dimension at $p = 2$. In the offline data set $\{(X^s, D_s, Y_s)\}_{s=1}^N$, each $X^s$ is drawn uniformly at random i.i.d. from $[0,1]^2$. The reward, which depends on the covariates and treatment, is generated as $Y_s = \mu(X^s, D_s) + \epsilon_s =  0.2 + \onebld\{D_s = j\} \boldsymbol{\theta}_j^\top X^s + \epsilon_s,~\text{for}~ j \in [J] = \{0,1,2\}$, where $\epsilon_s \sim N(0,0.01)$, $\boldsymbol{\theta}_0 = [1, 0.5]$, $\boldsymbol{\theta}_1 = [-0.5, 1.0]$, and $\boldsymbol{\theta}_2 = [-0.5, -0.5]$. It is easy to verify that treatment $2$ is uniformly dominated and is therefore the worst treatment for all covariate values. The optimal treatment is always either treatment $0$ or treatment $1$, with the remaining one being suboptimal: treatment $0$ is optimal in the region $\{x_2 \le 3x_1\}$, while treatment $1$ is optimal in the complementary region $\{x_2 > 3x_1\}$. We design the historical policy of sampling treatments such that the probabilities of selecting the optimal, suboptimal, and worst treatments are $0.9$, $0.185$, and $0.015$, respectively. Note that in this setting, both the suboptimal and the worst treatments are sampled with low probability. Consequently, their rewards may be estimated with high variance when using DR and IPW. As a result, learning algorithms built upon these estimators are expected to be misled by the optimistically biased value estimates. We compare the policies learned by two representative methods, which are widely used in the literature:
\begin{itemize}
    \item[$\bullet$] \textbf{IPW Learning} (IPWL).This algorithm selects a policy by maximizing the IPW value estimate. In particular, \cite{fang2025treatment} uses a PIP algorithm for solving the IPW-based policy optimization problem.
    \item[$\bullet$] \textbf{DR Learning} (DRL). This algorithm selects a policy by maximizing the DR value estimate. 
\end{itemize}
For all algorithms, we use the PIP method to conduct policy optimization, as its computational efficiency for large-scale problems has been demonstrated in Section~\ref{sec:computational_performance}. 

Each algorithm learns a linear policy from the simulated offline data, producing a parameter vector $\widehat{\beta}^j$ for $j \in \{0,1,2\}$. Once a policy is learned, we evaluate its performance by computing the suboptimality gap on a newly generated test set consisting of $10{,}000$ samples. Specifically, given a new covariate vector $X$ drawn from $[0,1]^2$, each algorithm selects a treatment such that $\hat{j} = \displaystyle{\argmax_{j \in [J]}} \langle \widehat{\beta}^{\,j}, X \rangle$. We assume access to an oracle that knows the true reward function and therefore always selects the truly optimal treatment for any given $X$. The (empirical) suboptimality gap of an algorithm is computed as the average difference, over the $10{,}000$ test samples, between the oracle's reward and the reward obtained by the policy produced by the algorithm. We conduct 30 independent runs for each offline data set with sample size $N \in \{400,600,800,1{,}000\}$ and report the empirical suboptimality gaps averaged over all runs.

\subsubsection{Results of Statistical Performance}
\begin{figure}[t]
    \includegraphics[scale = 0.6]{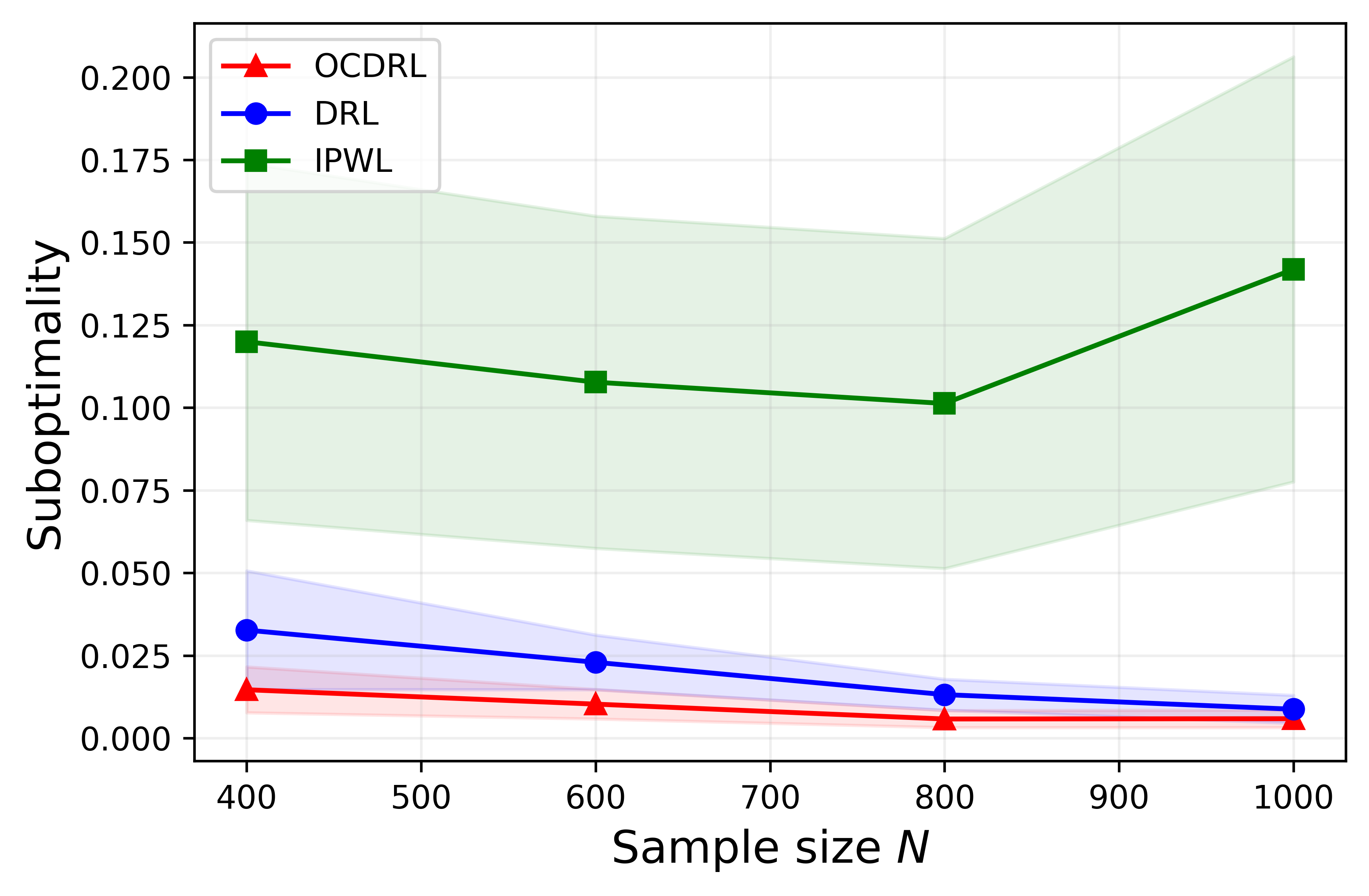}
    \caption{Out-of-Sample Suboptimality Gaps}
    \vspace{-10mm}
    \floatfoot{\textit{Note.} The confidence bands correspond to the mean $\pm$ two times the empirical standard deviation computed over 30 runs.}
    \label{fig:learning_performance}
\end{figure}

We visualize the suboptimality gaps of the three algorithms in Figure~\ref{fig:learning_performance}. The results demonstrate the superior performance of our proposed OCDRL algorithm relative to the two benchmark methods. We highlight two observations. First, IPWL has the worst suboptimality gap and does not show improvement as the sample size grows. This behavior can be attributed to the high variance of the IPW-based policy value estimator, particularly in our experimental setting where certain treatments have very small propensity scores. Second, although DRL is known to improve upon IPWL by partially reducing variance through outcome modeling, its continued reliance on inverse propensity weighting still results in substantially higher variance in policy value estimation compared with our weight-clipping approach. Consequently, its performance is dominated by our proposed OCDRL.

\begin{figure}[tbh]
    \includegraphics[width = \textwidth]{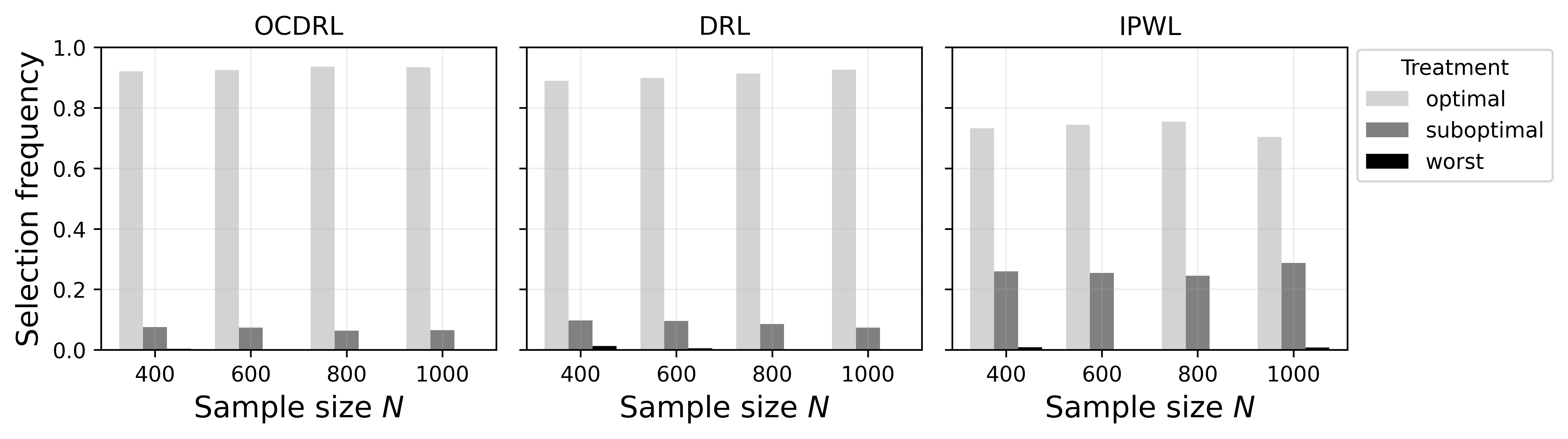}
    \caption{Treatment Selection Behavior of Different Algorithms}
    \vspace{-10mm}
    \floatfoot{\textit{Note. }Each frequency is computed as the average of the corresponding per-run frequencies over 30 independent runs.}
    \label{fig:selection_frequency}
\end{figure}

We further visualize the behavior of the three algorithms in Figure~\ref{fig:selection_frequency}, which provides additional evidence that high variance in policy evaluation can mislead policy learning algorithms. Recall that in the experimental design, the suboptimal and worst treatments are assigned very small propensity scores, which substantially inflate the variance of policy value estimates produced by the DR and IPW estimators. As a consequence, DRL and IPWL tend to select suboptimal and worst treatments more frequently than OCDRL, driven by spuriously inflated policy value estimates. In contrast, OCDRL effectively mitigates the high variance issue via weight clipping. Moreover, the empirical results indicate that the bias introduced by variance reduction is small enough, since it does not cause OCDRL to favor suboptimal or worst treatments more frequently than the benchmark algorithms.

\subsection{Assessment with Real-World Data: Targeting Policy for Insurance Products}
We further evaluate the proposed policy learning framework using a real-world data set from a randomized controlled trial (RCT) conducted by \cite{cai2015social}, which studies the heterogeneous effects of information sessions on farmers’ weather insurance take-up rates.\endnote{The data set is available at: Dean Eckles and Jeremy Yang. Counterfactual policy learning in 15.838. https://kaggle.com/competitions/mitsloancounterfactualpolicy, 2019. Kaggle.} This data set is ideal for performance evaluation because it is collected from an RCT, which implies that, (i) the assumptions on the offline data (Assumption~\ref{asp:dgp}) are naturally satisfied, and (ii) the known treatment probabilities (propensity scores) enable unbiased estimation of counterfactual rewards, providing a reliable measure of performance. In the experiment, households were randomly assigned to an intensive information session that provided detailed explanations of the expected benefits of the insurance product. The treatment is binary, with $J=1$ indicating assignment to the intensive information session. Households’ insurance purchase decisions were recorded, along with a rich set of covariates capturing household-level heterogeneity.

Our policy learning problem is motivated by the observation that individuals respond heterogeneously to insurance information sessions/advertising campaigns (see, e.g., \citealt{aizawa2018advertising, cai2015social, cole2011prices, shapiro2020advertising}), which calls for the design of a cost-effective targeting policy that delivers information sessions to individuals who are most likely to be influenced by such interventions. We assume that the expected profit from selling insurance to a household is 100 RMB, while delivering an intensive information session incurs a cost of 15 RMB per household. This defines the insurer’s reward function: $Y = 100R - 15D$, where $R, D \in \{0,1\}$ denote the purchase decision and the treatment received, respectively. We consider a targeting policy from the linear policy class that depends on three covariates: the household's risk aversion, perceived probability of natural disasters in the following year, and share of rice income in total income. Such information can often be obtained by insurance companies through customer surveys. \endnote{See, e.g., a consumer survey conducted by Swiss Re (one of the world's leading providers of insurance): https://www.swissre.com/reinsurance/insights/asia-protection-gap-consumer-survey.html. We also refer interested readers to \cite{cai2015social} for further details on the measurement of these covariates.}

We split the data set ($N = 4{,}715$) into a training set and an evaluation set in a $1{:}1$ ratio. Policies learned from the training set by our proposed OCDRL and the benchmark methods are evaluated on the evaluation set, where counterfactual rewards are estimated using the IPW estimator, which provides an unbiased estimate of the mean reward given known treatment assignment probabilities. As shown in Table~\ref{tab:real_data_performance}, the policy learned by our proposed OCDRL achieves higher expected rewards than the benchmark methods. We also note that, although our algorithm requires more computation time than DRL, which shares a similar problem structure, this additional effort arises from optimizing the clipping threshold to achieve improved statistical performance.

\begin{table}[h]
\centering
\caption{Rewards and Computation Time of Different Algorithms}
\label{tab:real_data_performance}
\begin{tabular}{cccc}
\hline
Method & OCDRL   & DRL     & IPWL    \\ \hline
Reward & 46.71 & 46.12 & 38.87 \\
Time   & 1.39s   & 0.58s   & 1.91s   \\ \hline
\end{tabular}\\
\vspace{1mm}
\small{
\textit{Note.} All learning algorithms use the PIP as the embedded computation method. Reward is measured by RMB.}
\end{table}

\begin{figure}[h]
    \includegraphics[scale = 0.5]{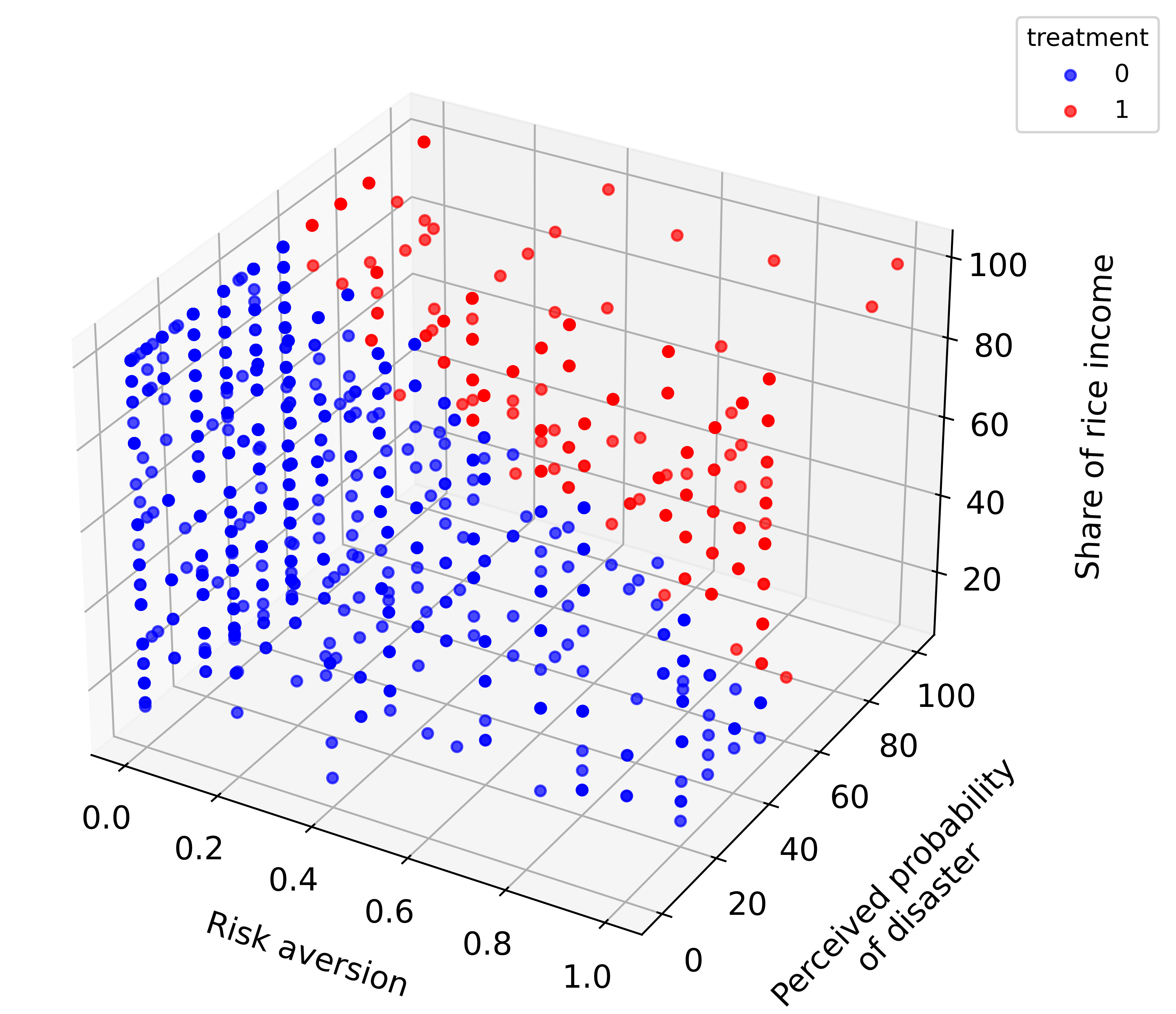}
    \caption{Targeting Policy Learned by OCDRL}
    % \vspace{-20mm}
    \floatfoot{\textit{Note.} Treatment $1$ indicates that a household is assigned to the information sessions.}
    \label{fig:real_data_policy}
\end{figure}

We close this section by discussing the managerial implications of the targeting policy learned by our proposed OCDRL (a visualization of the policy can be found in Figure \ref{fig:real_data_policy} in the online appendix). First, households’ perceived probability of a disaster is a dominant covariate influencing the assignment of information sessions. All households targeted by our learned policy have a perceived disaster probability of at least $75\%$. Second, the covariates interact in the assignment rule. For example, conditional on a high perceived disaster probability, households with lower risk aversion tend to be treated primarily when the share of rice income is high, but risk-averse households should be assigned even when the share of rice income is low.

\section{Concluding Remarks}
In this paper, we systematically develop an offline policy learning framework based on a weight-clipping estimator, with the desire for mitigating the adverse impact of high variance in policy evaluation on downstream policy optimization. Focusing on a linear policy class, we address the optimization challenges arising from the associated bilevel and discontinuous objective by (i) explicitly characterizing the optimal weight-clipping threshold as the minimizer of MSE in policy evaluation, and (ii) reformulating the policy optimization problem as a Heaviside composite optimization problem. We further tackle the computational challenges by developing a progressive integer programming algorithm, adapted from \cite{fang2025treatment}, which offers substantial computational advantages over conventional MIP approaches. The success of our optimization and computational approaches renders weight-clipping-based policy learning practically implementable at scale. We also establish a theoretical suboptimality guarantee for the proposed framework, which reveals how policy learning performance benefits from the optimized choice of the weight-clipping threshold.

%\THEEndNotes
\begingroup \parindent 0pt \parskip 0.0ex \def\enotesize{\small} \theendnotes \endgroup

% Appendix here
% Options are (1) APPENDIX (with or without general title) or
%             (2) APPENDICES (if it has more than one unrelated sections)
% Outcomment the appropriate case if necessary
%
% \begin{APPENDIX}{<Title of the Appendix>}
% \end{APPENDIX}
%
%   or
%
% \begin{APPENDICES}
% \section{<Title of Section A>}

% \end{APPENDICES}

% Acknowledgments here
% \ACKNOWLEDGMENT{}

% References here (outcomment the appropriate case)

% CASE 1: BiBTeX used to constantly update the references
%   (while the paper is being written).
%\bibliographystyle{informs2014} % outcomment this and next line in Case 1
%\bibliography{<your bib file(s)>} % if more than one, comma separated

%\bibliographystyle{informs2014} % outcomment this and next line in Case 1
%\bibliography{sample} % if more than one, comma separated

% CASE 2: BiBTeX used to generate mypaper.bbl (to be further fine tuned)
%\input{mypaper.bbl} % outcomment this line in Case 2

%If you don't use BiBTex, you can manually itemize references as shown below.

%\bibliographystyle{nonumber}
\bibliographystyle{informs2014}
\bibliography{main}

\ECSwitch
\OneAndAHalfSpacedXI
\ECHead{Online Appendices}
\section{Proofs of Statements in Section~\ref{sec:preliminaries}}
\subsection{Proof of Proposition \ref{prop: Dominant Term of Variance}}
\textit{Proof.} Note that when condition on $X^s$, $\wh{\mu}(X^s, g(X^s))$ is a constant and does not contribute to the variance, so it is omitted. For notational simplicity, we drop the subscripts and superscripts $s$.
\begin{align}
    \notag &\Var(\wh{\Gamma}^{\mathrm{DR}}(g) \mid X) = \frac{1}{e(X,g(X))^2} \cdot \Var ( \onebld\{g(X) = D\}(Y -\wh{\mu}(X,g(X))) \mid X)\\
    \label{variance decomposition}&= \frac{1}{e(X,g(X))^2} \cdot \left(\E[\onebld\{g(X) = D\}^2 (Y -\wh{\mu}(X,g(X)))^2 \mid X] - \E[\onebld\{g(X) = D\} (Y -\wh{\mu}(X,g(X))) \mid X]^2\right). 
\end{align}

To rewrite the above equation, we note that
\begin{align}
   \notag &\E[\onebld\{g(X) = D\}^2 (Y -\wh{\mu}(X,g(X)))^2 \mid X] \\
   \notag&= \E[\onebld\{g(X) = D\} (Y -\wh{\mu}(X,g(X)))^2 \mid X]\\
    \notag&= e(X,g(X)) \cdot \E[Y^2 - 2Y \wh{\mu}(X,g(X)) + \wh{\mu}(X,g(X))^2 \mid X]\\
    \label{var decompose 1}&= e(X,g(X)) \cdot (\E[Y^2 \mid X] - 2 \E[Y\mid X] \wh{\mu}(X,g(X)) + \wh{\mu}(X,g(X))^2).
\end{align}

In addition,
\begin{align}
    \label{var decompose 2}e(X,g(X))^2 \cdot (\E[Y \mid X] - \wh{\mu}(X,g(X)))^2 =  e(X,g(X))^2 \cdot (\E[Y \mid X]^2 - 2\E[Y \mid X] \wh{\mu}(X,g(X)) + \wh{\mu}(X,g(X))^2).
\end{align}

Therefore, the inner conditional variance term in Equation \eqref{variance decomposition} can be rewritten as 
{\small
\begin{align}
\notag&\E[\onebld\{g(X) = D\}^2 (Y -\wh{\mu}(X,g(X)))^2 \mid X] - \E[\onebld\{g(X) = j\} (Y -\wh{\mu}(X,g(X))) \mid X]^2 \\
   \notag&=e(X,g(X)) \cdot (\E[Y^2 \mid X] - 2 \E[Y\mid X] \wh{\mu} + \wh{\mu}^2) - e(X,g(X))^2 \cdot (\E[Y \mid X]^2 - 2\E[Y \mid X] \wh{\mu}(X,g(X)) + \wh{\mu}(X,g(X))^2)\\
   \notag&= e(X,g(X)) \cdot (\E[Y^2 \mid X] -\E[Y \mid X]^2 + \E[Y \mid X]^2 - 2 \E[Y\mid X] \wh{\mu}(X,g(X)) + \wh{\mu}(X,g(X))^2)\\
   \notag&~~~~- e(X,g(X))^2 \cdot (\E[Y \mid X]^2 - 2\E[Y \mid X] \wh{\mu}(X,g(X)) + \wh{\mu}(X,g(X))^2)\\
    \label{var decompose 3}&= e(X,g(X)) \Var(Y \mid X) + e(X,g(X))(1-e(X,g(X))) (\E[Y \mid X] - \wh{\mu}(X,g(X)))^2
\end{align}}

Finally, plugging (\ref{var decompose 3}) into (\ref{variance decomposition}) to obtain
\begin{align}
    \label{variance decomposition final}\Var(\wh{\Gamma}(g) \mid X) = \frac{1}{e(X,g(X))} \Var(Y \mid X) + \left(\frac{1}{e(X,g(X))}-1\right) \cdot (\E[Y\mid X]-\wh{\mu}(X,g(X)))^2.
\end{align}

By the assumption that $\Var(Y \mid X) \leq L'$, we can lower bound (\ref{variance decomposition final}) by
\begin{align*}
    \frac{L'}{e(X,g(X))}.
\end{align*}

By the assumption that $\Var(Y \mid X) \leq U'$ and $(\E[Y \mid X] - \wh{\mu}(X,g(X)))^2 \leq M^2$ (which follows from Assumption \ref{asp:dgp}(C)), we can upper bound (\ref{variance decomposition final}) by
\begin{align*}
    \frac{L' + M^2}{e(X,g(X))}.
\end{align*}

Letting $L = L'$ and $U = U' + M^2$ to yield the desired result. \halmos

\section{Proofs of Statements in Section~\ref{sec:formulation}}
\subsection{Proof of Proposition \ref{prop:MSE_CDR}}
\begin{proof}{Proof.}
    First, for any policy $g \in \Gcal$, following the routine bias-variance decomposition, we have
    \begin{align*}
        \E\left[(\wh{\Psi}^{\mathrm{CDR}}(g) - \Psi(g))^2\right] = \underbrace{\E\left[\wh{\Psi}^{\mathrm{CDR}}(g) - \Psi(g)\right]^2}_{\mathrm{Bias}^2} + \underbrace{\Var\left(\wh{\Psi}^{\mathrm{CDR}}(g)\right)}_{\mathrm{Variance}}.
    \end{align*}
     Note that for any given policy $g \in \Gcal$, the bias term is only incurred when $e(X,g(X))^{-1} > \tau$; otherwise, the CDR estimator switches to the DR estimator, which enjoys unbiasedness given that the propensity score is known. To see this,  we show that the bias term can be decomposed as
    \begin{align*}
        &\E\left[\wh{\Psi}^{\mathrm{CDR}}(g) - \Psi(g)\right]\\ 
        &\quad = \underbrace{\E\left[(\wh{\mu}(X, g(X)) - \mu(X,g(X))) \cdot \onebld\left\{\frac{1}{e(X,g(X))} > \tau\right\}\right]}_{\mathrm{(i)}}\\
        &\quad \quad + \underbrace{\E\left[\left(\wh{\mu}(X,g(X)) + \frac{\onebld\{g(X) = D\}}{e(X,g(X))}(Y-\wh{\mu}(X,g(X)))- \mu(X,g(X))\right) \cdot \onebld\left\{\frac{1}{e(X,g(X))} \leq  \tau\right\}\right]}_{\mathrm{(ii)}}.
    \end{align*}

    Term (i) can be upper bounded using the bounded reward assumption (Assumption \ref{asp:dgp}(C)), such that
    \begin{align*}
        (i) \leq M \cdot \E\left[\onebld\left\{\frac{1}{e(X,g(X))} > \tau \right\}\right].
    \end{align*}

As for term (ii), we show that it is unbiased:
\begin{align*}
    &\E\left[\E\left[\left(\wh{\mu}(X,g(X)) + \frac{\onebld\{g(X) = D\}}{e(X,g(X))}(Y-\wh{\mu}(X,g(X)) )- \mu(X,g(X))\right)\bigg| X \right] \cdot \onebld\left\{\frac{1}{e(X,g(X))} \leq \tau\right\}  \right]\\
    &= \E\left[\left(\wh{\mu}(X,g(X)) + \frac{e(X,g(X))}{e(X,g(X))}(\mu(X,g(X))-\wh{\mu}(X,g(X)) )- \mu(X,g(X))\right) \cdot \onebld\left\{\frac{1}{e(X,g(X))} \leq \tau\right\}  \right]\\
    &= 0,
\end{align*}
where the second equality is due to the definition of propensity score, namely, $e(X,g(X)) = \Prob(g(X) \mid X)$.

Therefore, the squared bias term is bounded as follows
\begin{align*}
    \E\left[\wh{\Psi}^{\mathrm{CDR}}(g) - \Psi(g)\right]^2 \leq M^2 \cdot \E\left[\onebld\left\{\frac{1}{e(X,g(X))} > \tau \right\}\right]^2.
\end{align*}

Next, we upper bound the variance. The analysis starts from the fact that $\Var(A+B) \leq 2\Var(A)+2\Var(B)$, where $A$ and $B$ are not independent.
 {\small
 \begin{align*}
     &\Var\left(\wh{\Psi}^{\mathrm{CDR}}(g)\right)\\
     &\quad \leq \underbrace{2 \Var\left(\frac{1}{N} \sum_{s=1}^N \wh{\mu}(X^s,g(X^s))\right)}_{\mathrm{(iii)}} + \underbrace{2\Var\left(\frac{1}{N} \sum_{s=1}^N \frac{\onebld\{g(X^s) = D_s\}}{e(X^s,g(X^s))} \cdot \onebld\left\{\frac{1}{e(X^s,g(X^s))} \leq \tau \right\} (Y_s - \wh{\mu}(X^s, g(X^s)))\right)}_{\mathrm{(iv)}}.
 \end{align*}}%

Term (iii) can be bounded as follows.
\begin{align*}
    2 \Var\left(\frac{1}{N} \sum_{s=1}^N \wh{\mu}(X^s,g(X^s))\right) \leq \frac{2}{N} \E\left[\wh{\mu}(X,g(X))^2\right] \leq \frac{2}{N}M^2,
\end{align*}
where the first inequality uses the fact that $\Var(A) = \E[A^2]-\E[A]^2 \leq \E[A^2]$ and the final inequality uses Assumption \ref{asp:dgp}(C). 

Similarly, term (iv) can be bounded by
\begin{align*}
    &2\Var\left(\frac{1}{N} \sum_{s=1}^N \frac{\onebld\{g(X^s) = D_s\}}{e(X^s,g(X^s))} \cdot \onebld\left\{\frac{1}{e(X^s,g(X^s))} \leq \tau \right\} (Y_s - \wh{\mu}(X^s, g(X^s)))\right)\\
    &\quad  \leq \frac{2}{N}\E\left[\frac{\onebld\{g(X) = D\}}{e(X,g(X))^2} (Y-\wh{\mu}(X,g(X)))^2 \cdot \onebld\left\{\frac{1}{e(X,g(X))} \leq \tau \right\}\right]\\
    &\quad  \leq \frac{2M^2}{N} \E\left[\frac{\onebld\{g(X) = D\}}{e(X,g(X))^2} \cdot\onebld\left\{\frac{1}{e(X,g(X))} \leq \tau \right\}\right].
\end{align*}

Putting pieces together to yield
{\small
\begin{align*}
    &\E\left[(\wh{\Psi}^{\mathrm{CDR}}(g) - \Psi(g))^2\right]\\
    &\quad  \leq M^2 \cdot \E\left[\onebld\left\{\frac{1}{e(X,g(X))} > \tau \right\}\right]^2 + \frac{2M^2}{N} \E\left[\frac{\onebld\{g(X) = D\}}{e(X,g(X))^2} \cdot\onebld\left\{\frac{1}{e(X,g(X))} \leq \tau \right\}\right] + \frac{2}{N}M^2.
    \halmos
\end{align*}}%
\end{proof}

\section{Proofs of Statements in Section \ref{sec:math formulation}}
\subsection{Proof of Proposition \ref{prop:explict_expression_hat_tau}}
Recall that
\[
\wt{\theta}(\tau(g);g) = \begin{cases} \Phi_m(g) & \mbox{ if } C_{(m-1)} \leq \tau(g) < C_{(m)}, \quad \mbox{ for } m=1,2, \ldots, N, \\
 \Phi_{N+1}(g) \equiv 0 & \mbox{ if } \tau(g) \geq C_{(N)},\end{cases}
\]
which implies that the expression of $\wh{\tau}(g)$ can be explicitly computed as follows,
\[
\wh{\tau}(g) = \begin{cases} 
C_{(0)}  & \mbox{if } \Phi_1(g) \leq \min\{ \Phi_{i}(g):  i=2, \ldots, N+1\},\\
C_{(m)}  & \mbox{if } \left\{ \begin{array}{l} \Phi_{m+1}(g) < \min\{\Phi_{i}(g): i=1, \ldots, m\},\\ \mbox{ and } \Phi_{m+1}(g) \leq \min\{\Phi_{i}(g): i=m+2, \ldots, N+1\}\end{array} \right\} \mbox{ for } m=1, \ldots, N-1,  \\
C_{(N)}  & \mbox{if } \Phi_{N+1}(g) < \min\{ \Phi_i(g): i=1,2,\ldots, N\}.    \\ 
\end{cases}
\] 
Based on the above expression of $\wh{\tau}(g)$, in general, we have
\[
\begin{array}{lll}
\wh{\tau}(g) & = & C_{(0)} \cdot \displaystyle{\onebld_{[0,\infty)}} \left(\min_{2 \leq t \leq  N+1} \bar \Phi_{t,1}(g) \right) \\[0.1in]
& &+\displaystyle{\sum_{m=1}^{N-1}} C_{(m)} \cdot \onebld_{[0,\infty)} \left(\min_{m+2 \leq t \leq  N+1} \bar \Phi_{t,m+1}(g) \right)   \cdot \displaystyle{\onebld_{(0,\infty)}} \left(\min_{1 \leq t \leq  m} \bar \Phi_{t,m+1}(g) \right) \\[0.1in]
& & + C_{(N)} \cdot \displaystyle{\onebld_{(0,\infty)}} \left(\min_{1 \leq t \leq  N} \bar \Phi_{t,N+1}(g) \right).
\end{array}
\]
Furthermore, under Assumption \ref{asp:non-integer}, the above expression simplifies to
\begin{equation*} 
\begin{array}{ll}
\wh{\tau}(g)   
% & =\displaystyle{\sum_{m=0}^N C_{(m)}} \cdot \displaystyle{\onebld_{[0,\infty)}} \left(\min_{  t \in[N+1], t \neq m+1} \bar \Phi_{t,m+1}(g) \right) \\
=  \displaystyle{
\sum_{m=0}^N 
} \, C_{(m)} \cdot \displaystyle{\onebld_{(0,\infty)}} \left(\min_{ t \in[N+1], t \neq m+1}  \Phi_{t}(g) - \Phi_{m+1}(g) \right).
\end{array}\halmos
\end{equation*}

\subsection{Proof of Proposition \ref{prop: pruning}}
We prove this proposition by contradiction. Suppose there exists an optimal $m^\prime \neq N+1$ that gives the minimum $\Phi_{m^\prime}(g)$, such that $2\Tilde{\varphi}_{(m^\prime)} < 2(N-m^\prime)+1$. Given that $\Tilde{\varphi}_{(m^\prime)} \geq \varphi_{(m^\prime)}(g)$, we have $\varphi_{(m^\prime)}(g) \leq 2\Tilde{\varphi}_{(m^\prime)} < 2(N-m^\prime)+1$, which implies that $2(N-m^\prime)+1 - 2\varphi_{(m^\prime)}(g) > 0$. By the definition of $\Phi_m(g)$, we know that $\Phi_{m^\prime+1}(g) < \Phi_{m^\prime}(g)$, contradicting $m^\prime$ achieves the minimum.\Halmos 

\subsection{Proof of Proposition \ref{prop:final_expression_hat_tau}}

Under the pruning strategy, the explicit optimal solution $\wh{\tau}(g)$ can be rewritten as
\[
\wh{\tau}(g) = \begin{cases} 
{C}_{(m^*-1)}  & \mbox{if }   {\Phi}_{m^*}(g) \leq \min\{  {\Phi}_{i}(g):  i \geq m^*+1\},\\
{C}_{(m)}  & \mbox{if } \left\{ \begin{array}{l} {\Phi}_{m+1}(g) < \min\{{\Phi}_{i}(g): m \geq i    \geq m^* \}, \mbox{ and } \\ {\Phi}_{m+1}(g) \leq \min\{{\Phi}_{i}(g): N+1 \geq i\geq m+2\}\end{array} \right\} \mbox{ for } N-1 \geq m \geq m^*, \\
{C}_{({N})}  & \mbox{if } m^* = N  \mbox{ or } 
 {\Phi}_{{N}+1}(g) < \min\{ {\Phi}_i(g): N \geq i\geq m^*\}.   \\ 
\end{cases}
\]
Furthermore, Assumption \ref{asp:non-integer} guarantees that the distinct values of $\{\Phi_m^*(g)\}$ for any $g \in \mathcal G$. \halmos

\subsection{Proof of Proposition 6}
\begin{proof}{Proof.} {\bf (A)} We first observe that for any scalar $t_* \neq 0$, there exist a scalar
$\varepsilon_* > 0$ and a neighborhood ${\cal T}_*$ of $t_*$
such that $\onebld_{( \, -\varepsilon,\infty )}(t) \, = \,
\onebld_{[ \, 0, \infty )}(t_*)$ and $\onebld_{[\, \varepsilon,\infty )}(t) \, = \,
\onebld_{( \, 0, \infty )}(t_*)$ for any  $( t,\varepsilon ) 
\in {\cal T}_* \times ( 0,\varepsilon_* ].$
Then,  for  any $\xi \in \Xi$ with $\min_{ 1 \leq i \leq j-1}h_{j,i}(\xi,\overline{\boldsymbol{\beta}})\neq 0$, there exist a scalar
$\bar \varepsilon_1 > 0$ and a neighborhood ${\cal N}_1$ of $\overline{\boldsymbol{\beta}}$
such that for any $\boldsymbol{\beta} \in {\cal N}_1$ and $0 < \varepsilon \leq \bar{\varepsilon}_1$, 
\[\onebld_{( \, 0,\infty )}\left(\min_{ 1 \leq i \leq j-1}h_{j,i}(\xi,\boldsymbol{\beta})\right) = \onebld_{[\,\varepsilon,\infty )}\left(\min_{ 1 \leq i \leq j-1}h_{j,i}(\xi,\boldsymbol{\beta})\right)= \onebld_{[\,0,\infty )}\left(\min_{ 1 \leq i \leq j-1}h_{j,i}(\xi,\boldsymbol{\beta})-\,\varepsilon \right),\] 
which thus yields the claimed statement. 

{\bf (B)} For the local maximizer $\overline{\boldsymbol  \beta}$ to the optimization problem  \eqref{eq:CDR_policy_optimization}, based on statement {\bf (A)}, there exist $\varepsilon>0$ and a neighborhood $\mathcal N$ of $\overline{\boldsymbol  \beta}$ such that for any ${\boldsymbol  \beta} \in \mathcal N$, $
\wh{\Psi}_{\varepsilon}(\boldsymbol{\beta}) \leq \wh{\Psi}_{\rm HSC}(\boldsymbol{\beta}) \leq \wh{\Psi}_{\rm HSC}(\overline{\boldsymbol{\beta}}) = \wh{\Psi}_{\varepsilon}(\overline{\boldsymbol{\beta}})$. Thus $\overline{\boldsymbol{\beta}}$ is a local maximizer of \eqref{eq:CDR_policy_approx_optimization}. 

{\bf (C)} For the local maximizer $\overline{\boldsymbol{\beta}}$ of \eqref{eq:CDR_policy_approx_optimization}, under the sign-invariance condition, there exist $\varepsilon$ and a neighborhood $\mathcal N$ of $\overline{\boldsymbol  \beta}$ such that for any ${\boldsymbol  \beta} \in \mathcal N$, $\forall j \times s \in \{(\ell,t): \ell \in [J],  \min_{ 1 \leq i \leq \ell-1}h_{\ell,i}(X^t,\overline{\boldsymbol{\beta}}) =0\}$,
\[ 
\begin{array}{ll}
 {\onebld_{[ \, 0,\infty )} (h^{1}_j(X^s, \boldsymbol{\beta};\varepsilon))}  = \onebld_{( \, 0,\infty )}\left(\min_{ 1 \leq i \leq j-1}h_{j,i}(X^s,\boldsymbol{\beta})\right) \times  
\onebld_{[ \, 0,\infty )}\left(\min_{ j+1 \leq i \leq J } h_{j,i}(X^s,\boldsymbol{\beta})  \right), 
 \end{array}
\]  
and for all $s \mbox{ satisfying } \min_{  D_s+1 \leq i \leq J}h_{D_s,i}(X^s,\overline{\boldsymbol{\beta}}) =0  $, 
\[ 
\begin{array}{ll}
 \onebld_{( \, 0,\infty )}\left(\min_{ 1 \leq i \leq D_s-1}h_{D_s,i}(X^s,\boldsymbol{\beta})\right) \times  
\onebld_{[ \, 0,\infty )}\left(\min_{ j+1 \leq i \leq J } h_{D_s,i}(X^s,\boldsymbol{\beta})  \right) = {\onebld_{(\, 0,\infty )} (h^{2}_{D_s}(X^s, \boldsymbol{\beta};\varepsilon))}, 
 \end{array}
\]  
Based on \eqref{eq:obj_upper_simplified_beta}, it is then not difficult to see that for any ${\boldsymbol  \beta} \in \mathcal N$, 
\[
\wh{\Psi}_{\rm HSC}(\overline{\boldsymbol{\beta}}) \geq \wh{\Psi}_{\varepsilon}(\overline{\boldsymbol{\beta}}) \geq  \wh{\Psi}_{\varepsilon}({\boldsymbol{\beta}}) = \wh{\Psi}_{\rm HSC}({\boldsymbol{\beta}}),
\]
where the last equality is obtained by the above analysis under the sign-invariance condition. \halmos
\end{proof}

\section{Proofs of Statements in Section \ref{sec:stats_performance}}
\subsection{Proof of Theorem \ref{thm:suboptimality_upper_bound}}
We begin by introducing a standard measure of policy class complexity, which will be used to characterize the performance of policy learning algorithms.

\begin{definition}[Natarajan dimension]
   \textit{ Given a $J$-action policy class $\cal{G} = \{g\in \Gcal: \cal{X} \to [J]\}$, we say a set of $m$ points $\{\xi_1, \cdots, \xi_m\}$ is shattered by $\cal{G}$ if there exist two functions $f_1, f_2: \{\xi_1, \cdots, \xi_m\} \mapsto [J]$, such that
    \begin{itemize}
        \item[(i)] For any $k \in [m]$, $f_{1}(\xi_j) \neq f_2(\xi_j)$;
        \item[(ii)] For any subset $S \subseteq [m]$, there exists some $g \in \cal{G}$, such that $g(\xi_k) = f_1(\xi_k)$ for all $k \in S$ and $g(\xi_k) = f_2(\xi_k)$ for all $k \notin S$.
    \end{itemize}
    The Natarajan dimension of $\Gcal$, defined as $\cal{N}_{\mathrm{dim}}(\Gcal)$, is the largest size of a set of points shattered by $\Gcal$.}
\end{definition}
In particular, $\cal{N}_{\mathrm{dim}}(\Gcal_{\rm LPC}) \leq |J|(p+1)$ \citep[Corollary 29.8]{shalev2014understanding}.

\textbf{Additional notations.} For simplicity, we denote $\wh{\Psi}^{\mathrm{OCDR}}(g)$ and $\wh{\Gamma}_s^{\mathrm{OCDR}}(g)$ as $\wh{\Psi}(g)$ and $\wh{\Gamma}_s(g)$, respectively.

We also define
\begin{align*}
    \wh{B}(g, \wh{\tau}(g)) &\triangleq \frac{M^2}{N^2} \left(\left(\sum_{s=1}^N \onebld\left\{\frac{1}{e(X^s, g(X^s))} > \wh{\tau}(g)\right\}\right)^2 + \sum_{s=1}^N \frac{\onebld\{D_s = g(X^s)\}}{e(X^s, g(X^s))^2} \onebld\left\{\frac{1}{e(X^s, g(X^s))} \leq \wh{\tau}(g) \right\}\right),\\
    \wt{B}(g, \wh{\tau}(g)) &\triangleq \frac{M^2}{N^2} \left(\left(\sum_{s=1}^N \onebld\left\{\frac{1}{e(X^s, g(X^s))} > \wh{\tau}(g)\right\}\right)^2 + \sum_{s=1}^N \frac{1}{e(X^s, g(X^s))} \onebld\left\{\frac{1}{e(X^s, g(X^s))} \leq \wh{\tau}(g) \right\}\right),\\
    \Delta(g, \wh{\tau}(g)) &\triangleq M^2 \left(\frac{\wh{\tau}(g)}{N}\right)^{3/2},\\
    A(g) &\triangleq \max\left\{\wh{B}(g,\wh{\tau}(g))^{1/2}, \wt{B}(g,\wh{\tau}(g))^{1/2}, \Delta(g, \wh{\tau}(g))^{1/2}\right\}.
\end{align*}

In addition, the following notations will be used to develop the symmetrization technique. For each $s \in [N]$, we let $(D_s',Y_s')$ be an independent copy of $(D_s, Y_s)$ conditioned on $X^s$, i.e.,
\begin{align*}
    (D_s', Y_s') \perp (D_s, Y_s) \mid X^s \quad \text{and} \quad (D_s', Y_s') \stackrel{\mathrm{d}}{=} (D_s, Y_s) \mid X^s.
\end{align*}
With such copies, we denote
\begin{align*}
    \wh{\Gamma}_s'(g) &= \wh{\mu}(X^s, g(X^s)) + \frac{\onebld\{D_s' = g(X^s)\}}{e(X^s, g(X^s))} \onebld\left\{\frac{1}{e(X^s, g(X^s))} \leq \wh{\tau}(g) \right\}(Y_s'-\wh{\mu}(X^s, g(X^s))),\\
    \wh{B}'(g, \wh{\tau}(g)) &= \frac{M^2}{N^2} \left(\left(\sum_{s=1}^N \onebld\left\{\frac{1}{e(X^s, g(X^s))} > \wh{\tau}(g)\right\}\right)^2 + \sum_{s=1}^N \frac{\onebld\{D_s' = g(X^s)\}}{e(X^s, g(X^s))^2} \onebld\left\{\frac{1}{e(X^s, g(X^s))} \leq \wh{\tau}(g) \right\}\right)
\end{align*}

We define the realizations of $\wh{\Gamma}_s(g)$ and $\wh{\Gamma}_s'(g)$ as $\wh{\gamma}_s(g)$ and $\wh{\gamma}_s(g)$, respectively, such that
\begin{align*}
    \wh{\gamma}_s(g) &= \wh{\mu}(\xi^s, g(\xi^s)) + \frac{\onebld\{j_s = g(\xi^s)\}}{e(\xi^s, g(\xi^s))} \onebld\left\{\frac{1}{e(\xi^s, g(\xi^s))} \leq \tau \right\}(y_s-\wh{\mu}(\xi^s, g(\xi^s))),\\
    \wh{\gamma}_s'(g) &= \wh{\mu}(\xi^s, g(\xi^s)) + \frac{\onebld\{j_s' = g(\xi^s)\}}{e(\xi^s, g(\xi^s))} \onebld\left\{\frac{1}{e(\xi^s, g(\xi^s))} \leq \tau \right\}(y_s'-\wh{\mu}(\xi^s, g(\xi^s))).
\end{align*}
Here, $\{\xi^s\}_{s=1}^N$, $\{j_s\}_{s=1}^N$, $\{j_s'\}_{s=1}^N$, $\{y_s\}_{s=1}^N$, and $\{y_s'\}_{s=1}^N$ are realizations of $\{X_s\}_{s=1}^N$, $\{D_s\}_{s=1}^N$, $\{D_s'\}_{s=1}^N$, $\{Y_s\}_{s=1}^N$, and $\{Y_s'\}_{s=1}^N$, respectively.

\textbf{We now start proving Theorem \ref{thm:suboptimality_upper_bound}}. First, we decompose the suboptimality in the following way.
\begin{align*}
    \cal{L}(\wh{g}, \Gcal_{\rm LPC}) &= \Psi(g^*) - \Psi(\wh{g})\\
    &= \Psi(g^*) - \wh{\Psi}(g^*) + \wh{\Psi}(g^*) - \wh{\Psi}(\wh{g}) + \wh{\Psi}(\wh{g})- \Psi(\wh{g})\\
    &\leq \Psi(g^*) - \wh{\Psi}(g^*) + \wh{\Psi}(\wh{g})- \Psi(\wh{g})\\
    &\leq 2\max_{g \in \Gcal_{\rm LPC}} \left|\Psi(g) - \wh{\Psi}(g)\right|,
\end{align*}
where the first inequality uses the fact that $\wh{\Psi}(\wh{g}) \geq \wh{\Psi}(g^*)$, since $\wh{g} \in \argmax_{g \in \Gcal_{\rm LPC}} \wh{\Psi}(g)$.

Our next step is to show that the (normalized) estimation error is uniformly bounded over $g \in \Gcal_{\rm LPC}$. To do this, we can decompose it as follows.
\begin{align*}
    \frac{|\wh{\Psi}(g) - \Psi(g)|}{A(g)} &= \left|\frac{1}{N} \sum_{s=1}^N \frac{\wh{\Gamma}_s(g) - \Psi(g)}{A(g)} \right|\\
    &\leq \underbrace{\left|\frac{1}{N} \sum_{s=1}^N  \frac{\wh{\Gamma}_s(g) - \mu(X^s, g(X^s))}{A(g)}\right|}_{\text{(i)}} + \underbrace{\left|\frac{1}{N} \sum_{s=1}^N  \frac{\mu(X^s, g(X^s)) - \Psi(g)}{A(g)}\right|}_{\text{(ii)}}.
\end{align*}

Lemmas \ref{lemma:symmetrization} and \ref{lemma:tail_prob_rad} below give an upper bound for Term (i). The proofs of Lemmas \ref{lemma:symmetrization} and \ref{lemma:tail_prob_rad} are given in Sections \ref{sec:lemma:symmetrization} and \ref{sec:lemma:tail_prob_rad}, respectively.
\begin{lemma}\label{lemma:symmetrization}
    For any constant $\zeta>4$, we have
    {\small
    \begin{align*}
        &\Prob\left(\sup_{g \in \Gcal} \left|\frac{1}{N} \sum_{s=1}^N  \frac{\wh{\Gamma}_s(g) - \mu(X^s, g(X^s))}{A(g)} \right| \geq \zeta \right)\\
        &\quad \leq 2\sup_{\xi,j,j',y,y'} \Prob \Bigg(\sup_{g \in \Gcal}\left|\sum_{s=1}^N \rho_s\left\{ \wh{\gamma}_s(g) - \wh{\gamma}'_s(g)\right\}\right|\\
    &\quad \quad \quad \quad \quad \quad \quad \geq \frac{\zeta M}{8} \cdot \Big(\sum_{s=1}^N \onebld\left\{\frac{1}{e(\xi^s, g(\xi^s))} \leq \wh{\tau}(g) \right\} \cdot\left( \frac{\onebld\{j_s = g(\xi^s)\}}{e(\xi^s, g(\xi^s))^2} + \frac{\onebld\{j_s' = g(\xi^s)\}}{e(\xi^s, g(\xi^s))^2}\right)\Big)^{1/2}~\text{for all $g \in \Gcal_{\rm LPC}$}\Bigg),
    \end{align*}}%
where the probability is over i.i.d. Rademacher random variables $\rho_1, \cdots, \rho_N$ with other quantities fixed at their realizations.
\end{lemma}

\begin{lemma}\label{lemma:tail_prob_rad}
    For any fixed $\delta \in (0,1)$ and any fixed realizations $\{\xi_s\}_{s=1}^N$, $\{j_s\}_{s=1}^N$, $\{j_s'\}_{s=1}^N$, $\{y_s\}_{s=1}^N$, the following event holds with probability at least $1-\delta$ over the randomness of $\rho_1, \cdots, \rho_N$:
    \begin{align*}
        &\left|\sum_{s=1}^N \rho_s \cdot (\wh{\gamma}_s(g) - \wh{\gamma}'_s(g))\right|\\
        &\leq 2 M\sqrt{2|J|(p+1)\log(N|J|^2) + 2\log(2/\delta)}\\
        &\quad \quad \quad \quad\times \Bigg(\sum_{s=1}^N \onebld\left\{\frac{1}{e(\xi^s, g(\xi^s))} \leq \wh{\tau}(g) \right\} \cdot \left(\frac{\onebld\{j_s = g(\xi^s)\}+\onebld\{j_s' = g(\xi^s)\}}{e(\xi^s, g(\xi^s))^2}\right)\Bigg)^{\frac{1}{2}}, \quad \text{for all $g \in \Gcal_{\rm LPC}$}.
    \end{align*}
\end{lemma}
Letting $\zeta = 16M\sqrt{|J|(p+1)\log(N|J|^2) + \log(8/\delta)}$, Lemmas \ref{lemma:symmetrization} and \ref{lemma:tail_prob_rad} together imply that
\begin{align}\label{ineq:subopt_part1}
    \sup_{g \in \Gcal_{\rm LPC}} \left|\frac{1}{N} \sum_{s=1}^N  \frac{\wh{\Gamma}_s(g) - \mu(X^s, g(X^s))}{A(g)} \right| \leq 8M\sqrt{2|J|(p+1)\log(N|J|^2) + 2\log(8/\delta)} ,
\end{align}
with probability at least $1-\delta/2$.

Next, we bound Term (ii) in the following lemma, whose proof can be found in Section \ref{sec:lemma:term2}.

\begin{lemma}\label{lemma:term2}
    For any $\delta \in (0,1)$, with probability at least $1-\delta/2$, we have
\begin{align}\label{ineq:subopt_part2}
\sup_{g \in \Gcal_{\rm LPC}} \left|\frac{1}{N}\sum_{s=1}^N \frac{\mu(X^s, g(X^s))-\Psi(g)}{A(g)}\right| \leq  2M\sqrt{8 |J|(p+1) \log(N|J|^2) + \log(8/\delta)}.
    \end{align}
\end{lemma}

Applying a union bound to \eqref{ineq:subopt_part1} and \eqref{ineq:subopt_part2}, we conclude that
\begin{align*}
    \left|\frac{\wh{\Psi}(g) - \Psi(g)}{A(g)}\right|
    \leq 12M \sqrt{ 2|J|(p+1)\log(N|J|^2) + 2\log(16/\delta)},
\end{align*}
with probability at least $1-\delta$.

Combining Lemmas \ref{lemma:symmetrization}- \ref{lemma:term2}, we conclude that
\begin{align*}
    \cal{L}(\wh{g}, \Gcal_{\rm LPC}) \leq \min\left\{2\upsilon \cdot \left\{ A(g^\dagger)\right\}, M\right\},
\end{align*}
with probability at least $1-\delta$, where $\upsilon = 12M \sqrt{ 2 |J|(p+1)\log(N|J|^2) + 2\log(16/\delta)}$ and $g^\dagger = \argmax_{g \in \Gcal_{\rm LPC}}|\wh{\Psi}(g)-\Psi(g)|$. Here, we compare with the constant $M$ because of the bounded outcome assumption, i.e., $\Psi(g) \in [0,M]$ for all $g \in \Gcal_{\rm LPC}$. 

In $A(g^\dagger)$, we note that $\E[\wh{B}(g^\dagger, \wh{\tau}(g^\dagger)) \mid X] = \wt{B}(g^\dagger, \wh{\tau}(g^\dagger))$, which follows from the definition of the propensity score $e(X,g(X))$. Our next step is to prove that $\wt{B}(g^\dagger, \wh{\tau}(g^\dagger)) \leq \wt{B}(g^\dagger, \wt{\tau}(g^\dagger)) + \omega$, for some $\omega \geq 0$, with high probability. To this end, we introduce some new notations in the following.

For notational simplicity, we use $\wh{b}(\tau)$ and $\wt{b}(\tau)$ to represent $\wh{B}(g^\dagger, \tau(g^\dagger))$ and $\wt{B}(g^\dagger, {\tau}(g^\dagger))$, respectively. Similar to the definitions of $\wh{\tau}(g^\dagger)$ and $\wt{\tau}(g^\dagger)$, the minimizers of $\wh{b}(\tau)$ and $\wt{b}(\tau)$ are denoted by $\wh{\tau}$ and $\wt{\tau}$, respectively. Define the set of $\omega$-optimal solutions for $\wt{b}(\tau)$ as 
\begin{align*}
    S^\omega &\triangleq \{\tau \in \cal{T}: \wt{b}(\tau) \leq b(\wt{\tau}) + \omega\}, \quad \text{where $ \cal{T} \triangleq [0,1/\eta]$}.
\end{align*}
Our goal is equivalent to showing that $\wh{\tau} \in S^\omega$ with high probability. Consider a mapping $u: \cal{T}\setminus S^\omega \to S^\omega$. Note that if the set $\cal{T}\setminus S^\omega$ is empty, then any feasible value $\tau \in \cal{T}$ is an $\omega$-optimal solution of the problem $\displaystyle{\min_{\tau \in \cal{T}}} ~\wt{b}(\tau)$. Therefore, we assume that this set is nonempty. Suppose the mapping $u$ is constructed in a way such that for every $\tau \in \cal{T}\setminus S^\omega$,
\begin{align}\label{eq:u_mapping}
    \wt{b}(u(\tau)) \leq \wt{b}(\tau) - \omega. 
\end{align}
For $\tau \notin S^\omega$, define $Z(\tau) \triangleq \wh{b}(u(\tau))-\wh{b}(\tau)$ as a random variable (whose randomness is from treatments $D_s$). Note that $\E[Z(\tau) \mid X] = \wt{b}(u(\tau))-\wt{b}(\tau) \leq -\omega$, where the inequality is due to \eqref{eq:u_mapping}. We are now ready to state the following lemma.

\begin{lemma}\label{lemma:SAA}
Suppose there exists a constant $\sigma>0$, such that for every $\tau \in \cal{T} \setminus S^\omega$, the random variable $Z(\tau) - \E[Z(\tau) \mid X]$ is $\sigma$-sub-Gaussian. Then for $\omega>0$ and $\delta \in (0,1)$, when the sample size $N$ satisfies $N \geq \frac{2\sigma^2}{\omega^2} \ln \left(\frac{N}{\delta}\right)$, we have
\begin{align*}
    \wt{b}(\wh{\tau}) \leq \wh{b}(\wt{\tau})+\omega,
\end{align*}
with probability at least $1-\delta$.
\end{lemma}
The proof of Lemma \ref{lemma:SAA} is in Section \ref{sec:proof_lemma:SAA}. Lemma \ref{lemma:SAA} implies that when the sample size $N$ satisfies that condition $N \geq \frac{2\sigma^2}{\omega^2} \ln \left(\frac{N}{\delta}\right)$, we have $\wt{B}(g^\dagger, \wh{\tau}(g^\dagger)) \leq \wt{B}(g^\dagger, \wt{\tau}(g^\dagger)) + \omega$. Letting $\delta = \delta/2$ to conclude that
\begin{align*}
    \cal{L}(\wh{g}, \Gcal_{\rm LPC}) \leq \min\left\{2\upsilon \cdot \max\left\{ \wh{B}(g^\dagger, \wh{\tau}(g^\dagger))^{1/2}, \wt{B}(g^\dagger,\wt{\tau}(g^\dagger))^{1/2}+\omega^{1/2}, \Delta(g^\dagger, \wh{\tau}(g^\dagger))^{1/2}\right\}, M\right\},
\end{align*}
with probability at least $1-\delta$. \halmos

\subsection{Proof of Corollary \ref{cor:convergence_rate}}
By definition,
{\small
\begin{align*}
    \wh{B}(g^\dagger, \wh{\tau}(g^\dagger)) &= \frac{M^2}{N^2} \left(\left(\sum_{s=1}^N \onebld\left\{\frac{1}{e(X^s, g^*(X^s))} > \wh{\tau}(g^\dagger)\right\}\right)^2 + \sum_{s=1}^N \frac{\onebld\{D_s = g^\dagger(X^s)\}}{e(X^s, g^\dagger(X^s))^2} \onebld\left\{\frac{1}{e(X^s, g^\dagger(X^s))} \leq \wh{\tau}(g^\dagger) \right\}\right)\\
    &\leq \frac{M^2}{N^2} \left(N \sum_{s=1}^N \onebld\left\{\frac{1}{e(X^s, g^\dagger(X^s))} > \wh{\tau}(g^\dagger)\right\} + \sum_{s=1}^N \frac{\onebld\{D_s = g^\dagger(X^s)\}}{e(X^s, g^\dagger(X^s))^2} \onebld\left\{\frac{1}{e(X^s, g^\dagger(X^s))} \leq \wh{\tau}(g^\dagger) \right\}\right).
\end{align*}}%
Therefore, 
{\small
\begin{align*}
    &\wh{B}(g^\dagger, \wh{\tau}(g^\dagger))^{1/2}\\ &\leq \frac{M}{N} \left(N \sum_{s=1}^N \onebld\left\{\frac{1}{e(X^s, g^\dagger(X^s))} > \wh{\tau}(g^\dagger)\right\} + \sum_{s=1}^N \frac{\onebld\{D_s = g^\dagger(X^s)\}}{e(X^s, g^\dagger(X^s))^2} \onebld\left\{\frac{1}{e(X^s, g^\dagger(X^s))} \leq \wh{\tau}(g^\dagger) \right\}\right)^{1/2}\\
    &\leq \frac{M}{N} \left(N^{1/2} \left(\sum_{s=1}^N \onebld\left\{\frac{1}{e(X^s, g^\dagger(X^s))} > \wh{\tau}(g^\dagger)\right\}\right)^{1/2} + \left(\sum_{s=1}^N \frac{\onebld\{D_s = g^\dagger(X^s)\}}{e(X^s, g^\dagger(X^s))^2} \onebld\left\{\frac{1}{e(X^s, g^\dagger(X^s))} \leq \wh{\tau}(g^\dagger) \right\}\right)^{1/2} \right)\\
    &\leq \frac{M}{N^{1/2}} \left(\sum_{s=1}^N \onebld\left\{\frac{1}{e(X^s, g^\dagger(X^s))} > \wh{\tau}(g^\dagger)\right\}\right)^{1/2} + \frac{M}{N} \cdot \left(\sum_{s=1}^N \frac{\onebld\{D_s = g^\dagger(X^s)\}}{e(X^s, g^\dagger(X^s))^2} \onebld\left\{\frac{1}{e(X^s, g^\dagger(X^s))} \leq \wh{\tau}(g^\dagger) \right\}\right)^{1/2}.
\end{align*}}%

Let
\begin{align*}
    G_s = \frac{\onebld\{D_s = g^\dagger(X^s)\}}{e(X^s, g^\dagger(X^s))^2} \onebld\left\{\frac{1}{e(X^s, g^\dagger(X^s))} \leq \wh{\tau}(g^\dagger) \right\}.
\end{align*}
We have
\begin{align*}
    \E[G_s] &= \E\left[\frac{1}{e(X,g^\dagger(X))} \onebld\{e(X,g^\dagger(X))^{-1} \leq \wh{\tau}(g^\dagger)\}\right],\\
    \Var(G_s) &\leq \E[G_s^2] \leq \wh{\tau}(g^\dagger)^3,\\
    |G_s| &\leq \wh{\tau}(g^\dagger)^2.
\end{align*}

We are now ready to apply Bernstein’s inequality: for any $\varsigma>0$,
\begin{align}\label{prob_bound:Bernstein}
    \Prob\left(\sum_{s=1}^N G_s - N \E[G_s] \geq \varsigma \right) \leq \exp\left(-\frac{\varsigma^2/2}{N\wh{\tau}(g^\dagger)^3 + \varsigma \wh{\tau}(g^\dagger)^2 /3}\right).
\end{align}
By setting $\varsigma = 2\log(1/\delta)\sqrt{N \wh{\tau}(g^\dagger)^3}$, we have 
\begin{align*}
    \frac{\varsigma^2/2}{N\wh{\tau}(g^\dagger)^3 + \varsigma \wh{\tau}(g^\dagger)^2 /3} = \frac{2\log^2(1/\delta)}{1+\frac{2}{3}\log(1/\delta) N^{-1/2}\wh{\tau}(g^\dagger)^{1/2}}.
\end{align*}
Given $\delta \in (0,\exp(-1))$, $\log(1/\delta) > 1$. In addition, given the condition $\eta_* > 1/N$, we have $\wh{\tau}(g^\dagger) \leq \frac{1}{\eta_*}<N$, where the first inequality is due to Proposition \ref{prop: pruning}, which implies that the possible values of $\wh{\tau}(g^\dagger)$ can only be one of the inverse propensity scores in $[1,1/\eta_*]$. Therefore, we have $\wh{\tau}(g^\dagger)/N < 1$, which is a sufficient condition that yields
\begin{align*}
    \frac{2\log^2(1/\delta)}{1+\frac{2}{3}\log(1/\delta) N^{-1/2}\wh{\tau}(g^\dagger)^{1/2}} \geq \log(1/\delta).
\end{align*}

Therefore, $\eqref{prob_bound:Bernstein} \leq \delta$, which implies that with probability at least $1-\delta$,
{\small
\begin{align*}
    \frac{1}{N} \left(\sum_{s=1}^N G_s\right)^{1/2} &\leq \frac{1}{N} \sqrt{N \cdot \E\left[\frac{1}{e(X,g^\dagger(X))} \onebld\{e(X,g^\dagger(X))^{-1} \leq \wh{\tau}(g^\dagger)\}\right] + 2\log(1/\delta)\sqrt{N \wh{\tau}(g^\dagger)^3}}\\
    &\leq \frac{1}{\sqrt{N}} \sqrt{\wh{\tau}(g^\dagger) + 2\log(1/\delta) \wh{\tau}(g^\dagger)^{3/2} N^{-1/2}}\\
    &\leq \frac{1}{\sqrt{N}} \left(\sqrt{\wh{\tau}(g^\dagger)} + \frac{2\log(1/\delta) \wh{\tau}(g^\dagger)^{3/2} }{2\wh{\tau}(g^\dagger)^{1/2}N^{1/2}}\right)\\
    &= \sqrt{\frac{\wh{\tau}(g^\dagger)}{N}} + \log(1/\delta) \frac{\wh{\tau}(g^\dagger)}{N}\\
    &\leq \log(1/\delta) \sqrt{\frac{\wh{\tau}(g^\dagger)}{N}} + \log(1/\delta) \sqrt{\frac{\wh{\tau}(g^\dagger)}{N}} = 2\log(1/\delta) \sqrt{\frac{\wh{\tau}(g^\dagger)}{N}}.
\end{align*}}%
where the third inequality uses the fact that $\sqrt{a+b} \leq \sqrt{a} \cdot \sqrt{1+b/a} \leq \sqrt{a}+\frac{b}{2\sqrt{2}}$ and the final inequality is due to $\log(1/\delta)>1$ (given $\delta \in (0, \exp(-1))$) and $\wh{\tau}(g^\dagger)/N < 1$.

Based on the above results, we conclude that
\begin{align*}
    \wh{B}(g^\dagger, \wh{\tau}(g^\dagger))^{1/2} &\leq M \left(\sqrt{\frac{\sum_{s=1}^N \onebld\left\{\frac{1}{e(X^s, g^\dagger(X^s))} > \wh{\tau}(g^\dagger)\right\}}{N}} + 2\log(1/\delta) \sqrt{\frac{\wh{\tau}(g^\dagger)}{N}}\right).
\end{align*}

To bound $\wt{B}(g, \wt{\tau})$, we use the fact that 
\begin{align*}
    \frac{1}{e(X^s, g(X^s))} \onebld\left\{\frac{1}{e(X^s,g(X^s))} \leq \wt{\tau}\right\} \leq \wt{\tau},
\end{align*}
which gives that
\begin{align*}
    \wt{B}(g, \wt{\tau}) \leq M \left(\sqrt{\frac{\sum_{s=1}^N \onebld\left\{\frac{1}{e(X^s, g^\dagger(X^s))} > \tau^*(g^\dagger)\right\}}{N}} + \sqrt{\frac{\tau^*(g^\dagger)}{N}}\right). 
\end{align*}

In addition, we need to bound the term $ \Delta(g^\dagger, \wh{\tau}(g^\dagger))^{1/2}$. Note that
\begin{align*}
    \Delta(g^\dagger, \wh{\tau}(g^\dagger))^{1/2} \leq M \left(\frac{\wh{\tau}(g^\dagger)}{N} \right)^{3/4} \leq M \left(\frac{\wh{\tau}(g^\dagger)}{N} \right)^{1/2},
\end{align*}
where the second inequality is due to $\wh{\tau}(g^\dagger)/N <1$.

Putting pieces together, we show that
\begin{align*}
    &\max\left\{ \wh{B}(g^\dagger, \wh{\tau}(g^\dagger))^{1/2}, \wt{B}(g^\dagger, \tau^*(g^\dagger))^{1/2}, \Delta(g^\dagger,\wh{\tau}(g^\dagger))^{1/2} \right\}\\
    &\quad \leq  M \left(\sqrt{\frac{\sum_{s=1}^N \onebld\left\{\frac{1}{e(X^s, g^\dagger(X^s))} > \wh{\tau}(g^\dagger)\right\}}{N}} + 2\log(1/\delta) \sqrt{\frac{\tau^*(g^\dagger) \vee \wh{\tau}(g)}{N}} + \omega\right),
\end{align*}
as long as $\delta \in [0,\exp(-1)]$. 

Finally, we take a union bound with the event $\cal{L}(\wh{g}, \Gcal_{\rm LPC}) \leq 2\nu \max\left\{ \wh{B}(g^\dagger, \wh{\tau}(g^\dagger))^{1/2}, \wt{B}(g^\dagger,\tau^*(g^\dagger))^{1/2}+\omega^{1/2}, \Delta(g^\dagger, \wh{\tau}(g^\dagger))^{1/2}\right\}$ in Theorem \ref{thm:suboptimality_upper_bound}. We know that with probability at least $1-\delta$,
{\small
\begin{align*}
    &\cal{L}(\wh{g}, \Gcal_{\rm LPC}) \\
    &\quad \leq c \cdot M \sqrt{\cal{N}_{\mathrm{dim}}(\Gcal_{\rm LPC}) \log(N|J|^2) + \log(2/\delta)} \cdot \left(\sqrt{\frac{\sum_{s=1}^N\onebld\left\{e(X^s, g^\dagger(X^s)) > \wh{\tau}(g^\dagger)\right\}}{N}} + 2\log(1/\delta) \sqrt{\frac{\tau^*+(g^\dagger)\vee \wh{\tau}(g)}{N}}\right)\\
    &\quad \leq 4c \cdot M \sqrt{\cal{N}_{\mathrm{dim}}(\Gcal_{\rm LPC}) \log(N|J|^2)} \left(\log(2/\delta)^{1/2}\sqrt{\frac{\sum_{s=1}^N\onebld\left\{e(X^s, g^\dagger(X^s)) > \wh{\tau}(g^\dagger)\right\}}{N}} + \log(2/\delta)^{3/2} \sqrt{\frac{\tau^*(g^\dagger)\vee \wh{\tau}(g)}{N}} \right)\\
    &\quad \leq 4c \cdot M \sqrt{|J|(p+1) \log(N|J|^2)} \left(\log(2/\delta)^{1/2}\sqrt{\frac{\sum_{s=1}^N\onebld\left\{e(X^s, g^\dagger(X^s)) > \wh{\tau}(g^\dagger)\right\}}{N}} + \log(2/\delta)^{3/2} \sqrt{\frac{\tau^*(g^\dagger)\vee \wh{\tau}(g)}{N}} \right).\halmos
\end{align*}}%

\section{Proofs of Auxiliary Lemmas}
\subsection{Proof of Lemma \ref{lemma:symmetrization}}\label{sec:lemma:symmetrization}
For any constant $\zeta>4$, we would like to control the following ``bad" event:
\begin{align*}
    \cal{E} &\triangleq \left\{\sup_{g \in \cal{G}}\left|\frac{1}{N} \frac{\sum_{s=1}^N  \wh{\Gamma}_s(g) - \mu(X^s, g(X^s))}{A(g)}\right| \geq \zeta \right\}\\
    &= \left\{\left|\frac{1}{N} \sum_{s=1}^N  \frac{\wh{\Gamma}_s(g^\dagger) - \mu(X^s, g^\dagger(X^s))}{A(g^\dagger)}\right| \geq \zeta  \right\},
\end{align*}
where $g^\dagger \triangleq \argmax_{g \in \Gcal_{\rm LPC}} \left|\frac{1}{N}\sum_{s=1}^N \frac{\wh{\Gamma}_s(g) - \mu(X^s, g(X^s))}{A(g)}\right|$. To this end, we define two auxiliary events:
\begin{align*}
    \cal{E}_1 &\triangleq \left\{ \left|\frac{1}{N} \sum_{s=1}^N  \wh{\Gamma}'_s(g^\dagger) - \mu(X^s, g^\dagger(X^s))\right|  \geq 2 \wt{B}(g^\dagger, \wh{\tau}(g^\dagger))^{1/2} \right\},\\
    \cal{E}_2 &\triangleq \left\{\wh{B}'(g^\dagger,\wh{\tau}(g^\dagger))^{1/2} \geq 2 \max\left\{\wh{B}(g^\dagger, \wh{\tau}(g^\dagger))^{1/2},\Delta(g^\dagger,\eta)^{1/2}\right\} \right\}.
\end{align*}
We will control the probability of $\cal{E}$ by showing that
\begin{align*}
    \Prob(\cal{E} \mid \{X^s\}_{s=1}^N) \leq 2 \Prob(\cal{E} \cap \cal{E}_1^c \cap \cal{E}_2^c \mid \{X^s\}_{s=1}^N).
\end{align*}
In order to achieve this goal, we will apply a symmetrization technique using Rademacher random variables to upper bound $\Prob(\cal{E} \cap \cal{E}_1^c \cap \cal{E}_2^c \mid \{X^s\}_{s=1}^N)$ by the tail probability of a Rademacher process. Lemma \ref{lemma:upper_bound_event2} in the following controls the probability of events $\cal{E}_1^c$ and $\cal{E}_2^c$ for a fixed policy $g \in \Gcal_{\rm LPC}$, conditioned on $\{X^s, D_s, Y_s\}_{s=1}^N$. The proof of Lemma \ref{lemma:upper_bound_event2} is deferred to Section \ref{sec:proof_lemma:upper_bound_event2}.

\begin{lemma}\label{lemma:upper_bound_event2}
For a fixed policy $g \in \Gcal_{\rm LPC}$, we have
\begin{align*}
    &\Prob\left(\left|\frac{1}{N} \sum_{s=1}^N  \wh{\Gamma}'_s(g) - \mu(X^s, g(X^s))\right|  \geq 2\wt{B}(g, \wh{\tau}(g))^{1/2} \mid \{X^s, D_s, Y_s\}_{s=1}^N\right) \leq \frac{1}{4},\\
    &\Prob\left(\wh{B}'(g,\wh{\tau}(g))^{1/2} \geq 2 \max\left\{\wt{B}(g, \wh{\tau}(g))^{1/2},\Delta(g, \wh{\tau}(g))^{1/2}\right\} \mid \{X^s, D_s, Y_s\}_{s=1}^N\right) \leq \frac{1}{4}
\end{align*}
\end{lemma}
Based on Lemma \ref{lemma:upper_bound_event2}, the tower property gives that $\Prob(\cal{E}_1^c \cap \cal{E}_2^c \mid \{X^s\}_{s=1}^N, \cal{E}) \geq \frac{1}{2}$. Multiplying both sides by $\Prob(\cal{E} \mid \{X^s\}_{s=1}^N)$ and applying Bayes' rule to obtain
\begin{align}\label{prob_bound:Bayes}
    \frac{1}{2} \Prob(\cal{E} \mid \{X^s\}_{s=1}^N) \leq \Prob(\cal{E}_1^c \cap \cal{E}_2^c \mid \{X^s\}_{s=1}^N, \cal{E}) \times \Prob(\cal{E} \mid \{X^s\}_{s=1}^N) = \Prob(\cal{E} \cap \cal{E}_1^c \cap \cal{E}_2^c \mid \{X^s\}_{s=1}^N).
\end{align}
It remains to upper bound the probability of event $\cal{E} \cap \cal{E}_1^c \cap \cal{E}_2^c$ conditioned on $\{X^s\}_{s=1}^N$. Note that under the event $\cal{E} \cap \cal{E}_1^c \cap \cal{E}_2^c$, we have
{\small
\begin{align*}
    \frac{1}{N}\left|\sum_{s=1}^N \wh{\Gamma}_s(g^\dagger) - \wh{\Gamma}'_s(g^\dagger)\right| &\geq \frac{1}{N}\left|\sum_{s=1}^N \wh{\Gamma}_s(g^\dagger) - \mu(X^s, g^\dagger(X^s))\right| - \frac{1}{N}\left|\sum_{s=1}^N \mu(X^s, g^\dagger(X^s)) - \wh{\Gamma}'_s(g^\dagger)\right|\\
    &\geq \frac{1}{N}\left|\sum_{s=1}^N \wh{\Gamma}_s(g^\dagger) - \mu(X^s, g^\dagger(X^s))\right|-2\wt{B}(g^\dagger,\wh{\tau}(g^\dagger))\\
    &\geq \frac{1}{N}\left|\sum_{s=1}^N \wh{\Gamma}_s(g^\dagger) - \mu(X^s, g^\dagger(X^s))\right|-\frac{\zeta}{2} \max\left\{\wh{B}(g^\dagger, \wh{\tau}(g^\dagger))^{1/2}, \wt{B}(g^\dagger,\wh{\tau}(g^\dagger))^{1/2}, \Delta(g^\dagger,\eta)^{1/2}\right\}\\
    &\geq \frac{\zeta}{2} \max\left\{\wh{B}(g^\dagger, \wh{\tau}(g^\dagger))^{1/2}, \wt{B}(g^\dagger,\wh{\tau}(g^\dagger))^{1/2}, \Delta(g^\dagger,\eta)^{1/2}\right\},
\end{align*}}%
where the third inequality is due to $\zeta>4$ and the final inequality is due to event $\cal{E}$. Based on the above inequality and event $\cal{E}_2$, we have
{\small
\begin{align*}
    &\frac{1}{N}\left|\sum_{s=1}^N \wh{\Gamma}_s(g^\dagger) - \wh{\Gamma}'_s(g^\dagger)\right|\\
    &\geq \frac{\zeta}{2} \max\left\{\wh{B}(g^\dagger, \wh{\tau}(g^\dagger))^{1/2}, \wt{B}(g^\dagger,\wh{\tau}(g^\dagger))^{1/2}, \Delta(g^\dagger,\eta)^{1/2}\right\}\\
    &\geq \frac{\zeta}{4} \max\left\{2\wh{B}(g^\dagger, \wh{\tau}(g^\dagger))^{1/2}, \wh{B}'(g^\dagger, \wh{\tau}(g^\dagger))^{1/2}\right\}\\
    &\geq \frac{\zeta}{8} \left(\wh{B}(g^\dagger, \wh{\tau}(g^\dagger))^{1/2}+ \wh{B}'(g^\dagger, \wh{\tau}(g^\dagger))^{1/2}\right)\\
    &\geq \frac{\zeta}{8} \cdot \frac{M}{N} \left( \left(\sum_{s=1}^N \frac{\onebld\{D_s = g^\dagger(X^s)\}}{e(X^s, g(X^s))^2}\onebld\left\{\frac{1}{e(X^s, g^\dagger(X^s))} \leq \wh{\tau}(g^\dagger) \right\}\right)^{1/2} + \left(\sum_{s=1}^N\frac{\onebld\{D_s' = g^\dagger(X^s)\}}{e(X^s, g^\dagger(X^s))^2}\onebld\left\{\frac{1}{e(X^s, g^\dagger(X^s))} \leq \wh{\tau}(g^\dagger) \right\} \right)^{1/2}\right)\\
    &\geq \frac{\zeta}{8}\left(\frac{M}{N} \cdot \left(\sum_{s=1}^N \frac{\onebld\{D_s = g^\dagger(X^s)\}}{e(X^s, g^\dagger(X^s))^2}\onebld\left\{\frac{1}{e(X^s, g^\dagger(X^s))} \leq \wh{\tau}(g^\dagger) \right\} + \frac{\onebld\{D_s' = g^\dagger(X^s)\}}{e(X^s, g^\dagger(X^s))^2}\onebld\left\{\frac{1}{e(X^s, g^\dagger(X^s))} \leq \wh{\tau}(g^\dagger) \right\}\right)^{1/2}\right),
\end{align*}}%
where the second last inequality is due to the definitions of $\wh{B}(g, \wh{\tau}(g))$ and $\wh{B}'(g, \wh{\tau}(g))$ and the last inequality uses the fact that $(a^{1/2}+b^{1/2}) \geq (a+b)^{1/2}$. Therefore, 
{\small
\begin{align*}
    \eqref{prob_bound:Bayes} &\leq \Prob\Bigg(\left|\sum_{s=1}^N \wh{\Gamma}_s(g^\dagger) - \wh{\Gamma}'_s(g^\dagger)\right|\\
    &\quad \quad \geq \frac{\zeta M}{8}\Big(\sum_{s=1}^N \onebld\left\{\frac{1}{e(X^s, g(X^s))} \leq \wh{\tau}(g) \right\} \cdot\left( \frac{\onebld\{D_s = g(X^s)\}}{e(X^s, g(X^s))^2} + \frac{\onebld\{D_s' = g(X^s)\}}{e(X^s, g(X^s))^2}\right)\Big)^{1/2} \bigg| \{X^s\}_{s=1}^N\Bigg)\\
    &\leq \Prob\Bigg(\sup_{g \in \Gcal_{\rm LPC}}\left|\sum_{s=1}^N \wh{\Gamma}_s(g) - \wh{\Gamma}'_s(g)\right| \cdot \Big(\sum_{s=1}^N \onebld\left\{\frac{1}{e(X^s, g(X^s))} \leq \wh{\tau}(g) \right\} \cdot\left( \frac{\onebld\{D_s = g(X^s)\}}{e(X^s, g(X^s))^2} + \frac{\onebld\{D_s' = g(X^s)\}}{e(X^s, g(X^s))^2}\right)\Big)^{-1/2}\\
    &\quad \quad \geq \frac{\zeta M}{8} \bigg| \{X^s\}_{s=1}^N\Bigg)\\
    &\leq \Prob\Bigg(\sup_{g \in \Gcal_{\rm LPC}}\left|\rho_s \left\{\sum_{s=1}^N \wh{\Gamma}_s(g) - \wh{\Gamma}'_s(g)\right\}\right| \cdot \Big(\sum_{s=1}^N \onebld\left\{\frac{1}{e(X^s, g(X^s))} \leq \wh{\tau}(g) \right\} \cdot\left( \frac{\onebld\{D_s = g(X^s)\}}{e(X^s, g(X^s))^2} + \frac{\onebld\{D_s' = g(X^s)\}}{e(X^s, g(X^s))^2}\right)\Big)^{-1/2}\\
    &\quad \quad \geq \frac{\zeta M}{8} \bigg| \{X^s\}_{s=1}^N\Bigg),
\end{align*}}%
where $\rho_s \in \{-1, +1\}$ for $s = 1, \cdots, N$ are i.i.d. Rademacher variables, taking the value $1$ with probability $1/2$ and $-1$ with probability $1/2$, and the above probability is taken over $\{D_s, Y_s, D_s', Y_s', \rho_s\}_{s=1}^N$. Given that $\{D_s,Y_s\}$ and $\{D_s', Y_s\}_{s=1}^N$ are exchangeable when conditioning on $\{X^s\}_{s=1}^N$, we apply the tower property to obtain
{\small
\begin{align*}
    &\Prob\Bigg(\sup_{g \in \Gcal_{\rm LPC}}\left|\rho_s \left\{\sum_{s=1}^N \wh{\Gamma}_s(g) - \wh{\Gamma}'_s(g)\right\}\right| \cdot \Big(\sum_{s=1}^N \onebld\left\{\frac{1}{e(X^s, g(X^s))} \leq \wh{\tau}(g) \right\} \cdot\left( \frac{\onebld\{D_s = g(X^s)\}}{e(X^s, g(X^s))^2} + \frac{\onebld\{D_s' = g(X^s)\}}{e(X^s, g(X^s))^2}\right)\Big)^{-1/2}\\
    &\quad \quad \geq \frac{\zeta M}{8} \bigg| \{X^s\}_{s=1}^N\Bigg)\\
    &= \E\Bigg[\Prob\Bigg(\sup_{g \in \Gcal_{\rm LPC}}\left|\rho_s \left\{\sum_{s=1}^N \wh{\Gamma}_s(g) - \wh{\Gamma}'_s(g)\right\}\right| \cdot \Big(\sum_{s=1}^N \onebld\left\{\frac{1}{e(X^s, g(X^s))} \leq \wh{\tau}(g) \right\} \cdot\left( \frac{\onebld\{D_s = g(X^s)\}}{e(X^s, g(X^s))^2} + \frac{\onebld\{D_s' = g(X^s)\}}{e(X^s, g(X^s))^2}\right)\Big)^{-1/2}\\
    &\quad \quad \geq \frac{\zeta M}{8} \bigg| \{X^s,D_s,D_s',Y_s, Y_s'\}_{s=1}^N\Bigg) \bigg| \{X^s\}_{s=1}^N\Bigg],
\end{align*}}%
where the inner conditional probability is over the randomness of the Rademacher random variables. Taking a supremum over all possible realizations of $\{X^s, D_s, D_s', Y_s, Y_s'\}_{s=1}^N$, we conclude that
{\small
\begin{align*}
\eqref{prob_bound:Bayes} &\leq \\
    \sup_{\xi,j,j',y,y'} &\Prob\Bigg(\sup_{g \in \Gcal_{\rm LPC}}\left|\rho_s \left\{\sum_{s=1}^N \wh{\gamma}_s(g) - \wh{\gamma}'_s(g)\right\}\right|\\
    &\quad \quad \geq \frac{\zeta M}{8} \cdot \Big(\sum_{s=1}^N \onebld\left\{\frac{1}{e(\xi^s, g(\xi^s))} \leq \wh{\tau}(g) \right\} \cdot\left( \frac{\onebld\{j_s = g(\xi^s)\}}{e(\xi^s, g(\xi^s))^2} + \frac{\onebld\{j_s' = g(\xi^s)\}}{e(\xi^s, g(\xi^s))^2}\right)\Big)^{1/2}~\text{for all $g \in \Gcal_{\rm LPC}$}\Bigg). \halmos
\end{align*}}%

\subsection{Proof of Lemma \ref{lemma:tail_prob_rad}}\label{sec:lemma:tail_prob_rad}
For given realizations $\xi,j,j',y,y'$ and policy $g \in \Gcal_{\rm LPC}$,
{\small
\begin{align*}
    &(\wh{\gamma}_s(g) - \wh{\gamma}_s'(g))^2\\
    &=\left(\frac{\onebld{D_s = g(\xi_s)}}{e(\xi_s, g(\xi_s))} \onebld \left\{\frac{1}{e(\xi_s, g(\xi_s))} \leq \wh{\tau}(g) \right\}(y_s-\wh{\mu}(\xi_s, g(\xi_s))) - \frac{\onebld{j'_s = g(\xi_s)}}{e(\xi_s, g(\xi_s))} \onebld \left\{\frac{1}{e(\xi_s, g(\xi_s))} \leq \wh{\tau}(g) \right\}(y'_s-\wh{\mu}(\xi_s, g(\xi_s)))\right)^2\\
    &\leq 2 \cdot \onebld \left\{\frac{1}{e(\xi_s, g(\xi_s))} \leq \wh{\tau}(g) \right\} \cdot \Bigg(\wh{\mu}(\xi_s, g(\xi_s))^2 \left(\frac{\onebld\{j'_s = g(\xi_s)\}}{e(\xi_s, g(\xi_s))} -\frac{\onebld\{D_s = g(\xi_s)\}}{e(\xi_s, g(\xi_s))}\right)^2\\
    &\hspace{5cm} + \left(\frac{\onebld\{D_s = g(\xi_s)\}}{e(\xi_s, g(\xi_s))}y_s -\frac{\onebld\{j'_s = g(\xi_s)\}}{e(\xi_s, g(\xi_s))}y'_s\right)^2\Bigg)\\
    &\leq 2 \cdot \onebld \left\{\frac{1}{e(\xi_s, g(\xi_s))} \leq \wh{\tau}(g) \right\} \cdot \Bigg(\wh{\mu}(\xi_s, g(\xi_s))^2 \left(\frac{\onebld\{j'_s = g(\xi_s)\}}{e(\xi_s, g(\xi_s))^2} +\frac{\onebld\{D_s = g(\xi_s)\}}{e(\xi_s, g(\xi_s))^2}\right)\\
    &\hspace{5cm} +\left(\frac{\onebld\{D_s = g(\xi_s)\}}{e(\xi_s, g(\xi_s))^2}y_s^2 +\frac{\onebld\{j'_s = g(\xi_s)\}}{e(\xi_s, g(\xi_s))^2}y_{s}^{'2}\right)\Bigg)\\
    &\leq 4M^2 \onebld \left\{\frac{1}{e(\xi_s, g(\xi_s))} \leq \wh{\tau}(g) \right\} \left(\frac{\onebld\{j'_s = g(\xi_s)+ \onebld\{D_s = g(\xi_s)}{e(\xi_s, g(\xi_s))^2}\right)
\end{align*}}%
where the first inequality uses the fact that $(a+b)^2 \leq 2a^2 + 2b^2$, the second inequality uses $(a-b)^2 \leq a^2 + b^2$, and the final inequality uses the assumption that $y,\wh{\mu} \in [0,M]$ (Assumption \ref{asp:dgp}(C)). By Hoeffding's inequality with respect to the randomness of $\rho$, we have
\begin{align*}
    &\mathbb{P}_\rho \left(\big|\sum_{s=1}^N \rho_s \cdot (\wh{\gamma}_s(g) - \wh{\gamma}_s'(g))\big| \geq \frac{\zeta }{8} \Big(\sum_{s=1}^N \onebld\left\{\frac{1}{e(\xi^s, g(\xi^s))} \leq \wh{\tau}(g) \right\} \cdot\left( \frac{\onebld\{D_s = g(\xi^s)\}}{e(\xi^s, g(\xi^s))^2} + \frac{\onebld\{D_s' = g(\xi^s)\}}{e(\xi^s, g(\xi^s))^2}\right)\Big)^{1/2} \right)\\
    &\leq 2\exp \left(-\frac{\zeta^2}{512M^2}\right),
\end{align*}
for any $\zeta > 0$. 

Let $\Gcal(\xi) \triangleq \{g(\xi^1), \cdots, g(\xi^N): g \in \Gcal_{\rm LPC}\}$. From Lemma \ref{lemma:Natarajan}, we have $|\Gcal(\xi)| \leq N^{\cal{N}_{\mathrm{dim}}(\Gcal_{\rm LPC})} |J|^{2\cal{N}_{\mathrm{dim}}(\Gcal_{\rm LPC})}$, where $|J|$ is the number of treatments. Taking a union bound over $\Gcal(\xi)$, we have
\begin{align*}
    &\mathbb{P}_{\rho}\left(\big|\sum_{s=1}^N \rho_s \cdot (\wh{\gamma}_s(g) - \wh{\gamma}_s'(g))\big| \geq \frac{\zeta}{8} \cdot \onebld \left\{\frac{1}{e(\xi_s, g(\xi_s))} \leq \tau \right\} \left(\frac{\onebld\{j'_s = g(\xi_s)+ \onebld\{j_s = g(\xi_s)}{e(\xi_s, g(\xi_s))^2} \right)^{\frac{1}{2}} ~\text{for all $g \in \Gcal_{\rm LPC}$}\right)\\
    &\leq 2 |\Gcal(\xi)|\exp \left(-\frac{\zeta^2}{512M^2}\right)\\
    &\leq 2N^{\cal{N}_{\mathrm{dim}}(\Gcal_{\rm LPC})} J^{2\cal{N}_{\mathrm{dim}}(\Gcal_{\rm LPC})}\exp \left(-\frac{\zeta^2}{512M^2}\right)\\
    &\leq 2N^{|J|(p+1)} |J|^{2|J|(p+1)}\exp \left(-\frac{\zeta^2}{512M^2}\right),
\end{align*}
where the final inequality is due to Lemma \ref{lemma:n_dim_lpc}.

Letting $\zeta = 16 M\sqrt{2|J|(p+1)\log(N|J|^2) + 2\log(2/\delta)}$ to yield the desired result. \halmos

\subsection{Proof of Lemma \ref{lemma:term2}}\label{sec:lemma:term2}
Let $X^{s'}$ be an i.i.d. copy of $X^s$, for $s = 1, \cdots, N$. Note that $\Psi(g) = \E[\mu(X, g(X))] = \E[Y(g(X))]$ is the policy value. For any fixed $\kappa > \sqrt{2}M$, we use Markov's inequality to show that
\begin{align*}
    \Prob\left(\left|\sum_{s=1}^N \mu(X^{s'},g(X^{s'})) - \Psi(g)\right| \geq \kappa \sqrt{N}\right) &\leq \frac{1}{\kappa^2N} \E\left[\left(\sum_{s=1}^N \mu(X^{s'},g(X^{s'})) - \Psi(g)\right)^2\right]\\
    &= \frac{1}{\kappa^2N} \sum_{s=1}^N \E\left[\left(\mu(X^{s'},g(X^{s'})) - \Psi(g)\right)^2\right]\\
    &\leq \frac{M^2}{\kappa^2} \leq \frac{1}{2}.
\end{align*}

Let $g^\star = \argmax_{g \in \Gcal_{\rm LPC}}|\sum_{s=1}^N \mu(X^s,g(X^s))-\Psi(g)|$, then we have
\begin{align}\label{prob_bound:greater_than_1/2}
    \Prob\left(\left|\sum_{s=1}^N \mu(X^{s'},g^\star(X^{s'})) - \Psi(g^\star)\right| < \kappa \sqrt{N} \bigg| \{X^s\}_{s=1}^N\right) \geq \frac{1}{2}.
\end{align}

For any $\kappa \geq \sqrt{2}M$, we have
\begin{align}\label{prob_bound:term2}
    \notag &\frac{1}{2} \Prob\left(\sup_{g \in \Gcal_{\rm LPC}} \left|\sum_{s=1}^N\mu(X^s, g(X^s)) - \Psi(g) \right| \geq 2\kappa \sqrt{N} \right)\\
    \notag &\quad = \frac{1}{2} \Prob\left(\left|\sum_{s=1}^N\mu(X^s, g^\star(X^s)) - \Psi(g^\star) \right| \geq 2\kappa \sqrt{N} \right)\\
    \notag &\quad \leq \Prob\left(\left|\sum_{s=1}^N\mu(X^{s'}, g^\star(X^{s'})) - \Psi(g^\star) \right| \geq \kappa \sqrt{N} \bigg| \{X^s\}_{s=1}^N\right) \times \Prob\left(\sum_{s=1}^N\left|\mu(X^{s}, g^\star(X^{s})) - \Psi(g^\star)\right| \geq 2\kappa \sqrt{N}\right)\\
    \notag &\quad = \Prob\left(\left|\sum_{s=1}^N\mu(X^{s'}, g^\star(X^{s'})) - \Psi(g^\star) \right| \geq \kappa \sqrt{N} ,\left|\sum_{s=1}^N\mu(X^{s}, g^\star(X^{s})) - \Psi(g^\star)\right| \geq 2\kappa \sqrt{N}\right)\\
    \notag &\quad \leq \Prob\left(\left|\sum_{s=1}^N\mu(X^{s'}, g^\star(X^{s'})) - \mu(X^{s}, g^\star(X^{s})) \right| \geq \kappa \sqrt{N} \right)\\
    &\quad \leq \Prob\left(\sup_{g \in \Gcal_{\rm LPC}}\left|\sum_{s=1}^N\mu(X^{s'}, g(X^{s'})) - \mu(X^{s}, g(X^{s})) \right| \geq \kappa \sqrt{N} \right)
\end{align}
where the first inequality is due to \eqref{prob_bound:greater_than_1/2}, the second equality uses the tower property, and the second last inequality is due to the triangle inequality (as we have seen in the Proof of Lemma \ref{lemma:symmetrization}).

Again, by letting $\rho_s$ be i.i.d. Rademacher random variables, and by the exchangeability of $X^{s'}$ and $X^s$, we have
\begin{align*}
    \eqref{prob_bound:term2} &= \Prob\left(\sup_{g \in \Gcal_{\rm LPC}}\left|\sum_{s=1}^N \rho_s \left(\mu(X^{s'}, g(X^{s'})) - \mu(X^{s}, g(X^{s}))\right) \right| \geq \kappa \sqrt{N} \right)\\
    &\quad \leq 2\Prob\left(\sup_{g \in \Gcal_{\rm LPC}}\left|\sum_{s=1}^N \rho_s \left(\mu(X^{s'}, g(X^{s'}))\right) \right| \geq \kappa \sqrt{N} / 2\right)\\
    &\quad = 2\E\left[\Prob_{\rho}\left(\sup_{g \in \Gcal_{\rm LPC}}\left|\sum_{s=1}^N \rho_s \left(\mu(X^{s'}, g(X^{s'}))\right) \right| \geq \kappa \sqrt{N} / 2 \bigg| \{X^s\}_{s=1}^N\right) \right]\\
    &\quad \leq 2\sup_{\xi \in \cal{X}} \Prob_{\rho} \left(\sup_{g \in \Gcal_{\rm LPC}} \left|\sum_{s=1}^N \rho_s \cdot \mu(\xi^s, g(\xi^s)) \right| \geq \kappa\sqrt{N}/2 \right).
\end{align*}

Similar to the proof of Lemma \ref{lemma:tail_prob_rad}, we apply a union bound to obtain
\begin{align*}
    \Prob_{\rho} \left(\sup_{g \in \Gcal_{\rm LPC}} \left|\sum_{s=1}^N \rho_s \cdot \mu(\xi^s, g(\xi^s)) \right| \geq \kappa\sqrt{N}/2 \right) &\leq |\Gcal(\xi)| \cdot \Prob_{\rho} \left(\sup_{g \in \Gcal_{\rm LPC}} \left|\sum_{s=1}^N \rho_s \cdot \mu(\xi^s, g(\xi^s)) \right| \geq \kappa\sqrt{N}/2 \right)\\
    &\leq N^{\cal{N}_{\mathrm{dim}}(\Gcal_{\rm LPC})} J^{2\cal{N}_{\mathrm{dim}}(\Gcal_{\rm LPC})} \Prob_{\rho} \left(\sup_{g \in \Gcal_{\rm LPC}} \left|\sum_{s=1}^N \rho_s \cdot \mu(\xi^s, g(\xi^s)) \right| \geq \kappa\sqrt{N}/2 \right)\\
    &\leq 2N^{\cal{N}_{\mathrm{dim}}(\Gcal_{\rm LPC})} J^{2\cal{N}_{\mathrm{dim}}(\Gcal_{\rm LPC})} \exp(-\kappa^2/8M^2)\\
    &\leq 2N^{|J|(p+1)} J^{2|J|(p+1)} \exp(-\kappa^2/8M^2),
\end{align*}
where the second inequality uses Lemma \ref{lemma:Natarajan}, the third inequality uses Hoeffding's inequality, and the final inequality uses Lemma \ref{lemma:n_dim_lpc}.

Taking $\kappa = M\sqrt{8 |J|(p+1) \log(N|J|^2) + \log(8/\delta)}$, we have proved that
\begin{align*}
    \sup_{g \in \Gcal_{\rm LPC}}\left|\sum_{s=1}^N \mu(X^s, g(X^s)) - \Psi(g) \right| \leq 2\kappa \sqrt{N},
\end{align*}
with probability at least $1-\delta/2$. Divide both sides by $N$, we have
\begin{align}\label{prob_bound:term2_intermediate}
    \sup_{g \in \Gcal_{\rm LPC}}\frac{1}{N} \left|\sum_{s=1}^N \mu(X^s, g(X^s)) - \Psi(g) \right| \leq 2\kappa \frac{\sqrt{N}}{N} = 2\kappa \frac{1}{\sqrt{N}},
\end{align}
with probability at least $1-\delta/2$. Finally, we note that
\begin{align*}
    A(g) = \max\left\{\wh{B}(g, \wh{\tau}(g))^{1/2}, \wt{B}(g,\wh{\tau}(g))^{1/2}, \Delta(g, \wh{\tau}(g))^{1/2}\right\} \geq \wt{B}(g,\wh{\tau}(g))^{1/2} \geq \frac{1}{\sqrt{N}},
\end{align*}
where the last inequality uses $M \geq 1$ (Assumption \ref{asp:dgp}(C)). Therefore, by \eqref{prob_bound:term2_intermediate}, we conclude that
\begin{align*}
     \sup_{g \in \Gcal_{\rm LPC}}\frac{1}{N} \left|\sum_{s=1}^N \frac{\mu(X^s, g(X^s)) - \Psi(g)}{A(g)} \right| \leq 2\kappa \frac{1}{\sqrt{N}} \frac{1}{A(g)} \leq 2\kappa \frac{1}{\sqrt{N}} \sqrt{N} = 2\kappa, 
\end{align*}
with probability at least $1-\delta/2$. \halmos

\subsection{Proof of Lemma \ref{lemma:SAA}}\label{sec:proof_lemma:SAA}
Since $\wh{\tau} \in \argmax_{\tau \in \cal{T}}\wh{b}(\tau)$, we have that $\wh{b}(\wh{\tau}) \leq \wh{b}(u(\wh{\tau}))$. Therefore, if $\hat{\tau} \notin S^\omega$, then $Z(\wh{\tau}) = \wh{B}(u(\wh{\tau})) - \wh{B}(\wh{\tau}) \geq 0$, which implies that $\{\hat{\tau} \notin S^\omega\} \subseteq \{\exists~\tau \notin S^\omega, Z(\tau) \geq 0\}$ or equivalently, $\Prob(\hat{\tau} \notin S^\omega \mid X) \leq \Prob(\exists~\tau \notin S^\omega, Z(\tau) \geq 0 \mid X)$. For a fixed $\tau \notin S^\omega$, we recall that $\E[Z(\tau) \mid X] \leq -\omega$. This implies that, on the event $\{Z(\tau) \geq 0\}$, $Z(\tau) - \E[Z(\tau) \mid X] \geq 0 - \E[Z(\tau) \mid X] \geq \omega$. Hence, we have that $\{Z(\tau) \geq 0\} \subseteq \{Z(\tau) - \E[Z(\tau)] \geq \omega \mid X\}$, or equivalently, for all $\tau \notin S^\omega$,
\begin{align}\label{ineq:z}
    \Prob(Z(\tau) \geq 0 \mid X) \leq \Prob(Z(\tau) - \E[Z(\tau)] \geq \omega \mid X) \leq \exp(-\frac{N\omega^2}{2\sigma^2}),
\end{align}
where the second inequality is due to the assumption that $Z(\tau) - \E[Z(\tau)]$ is a $\sigma$-sub-Gaussian random variable. We will show that this assumption holds at the end of the proof.

Using the union bound and \eqref{ineq:z}, we have
\begin{align*}
    \Prob(\exists~\tau \notin S^\omega, Z(\tau) \geq 0 \mid X) \leq \sum_{\tau \in \cal{T} \setminus S^\omega}  \Prob(Z(\tau) \geq 0 \mid X) \leq N \exp(-\frac{N\omega^2}{2\sigma^2}).
\end{align*}
Therefore, we conclude that
\begin{align*}
    \Prob(\wh{\tau} \notin S^\omega \mid X) \leq  N \exp(-\frac{N\omega^2}{2\sigma^2}).
\end{align*}

Since we hope that $ N \exp(-\frac{N\omega^2}{2\sigma^2}) \leq \delta$, with some algebra, we solve for the condition of $N$, such that $N \geq 2\sigma^2/\omega^2 \ln \left(\frac{N}{\delta}\right)$. Taking expectation over $X$ gives that
\begin{align*}
    \Prob(\wt{b}(\wh{\tau}) \leq \wt{b}(\wt{\tau}) + \omega) \geq 1-\delta.
\end{align*}

Finally, we conclude the proof by showing that $Z(\tau)-\E[Z(\tau) \mid X]$ is a $\sigma$-sub-Gaussian random variable. Recall that when condition on $X$, the randomness of $\wh{b}(\tau)$ is from the term $\sum_{s=1}^N \frac{\onebld \{D_s = g(X^s)\}}{e(X^S, g(X^s))} \onebld \left\{\frac{1}{e(X^s, g(X^s)) \leq \tau(g^\dagger)}\right\}$. Let $R(\tau) \triangleq \sum_{s=1}^N r_s(\tau)$, where $r_s(\tau) \triangleq \frac{\onebld \{D_s = g(X^s)\}}{e(X^S, g(X^s))} \onebld \left\{\frac{1}{e(X^s, g(X^s)) \leq \tau(g^\dagger)}\right\}$. Let $z_s(\tau) = r_s(u(\tau))- r_s(\tau)$, we have $R(u(\tau)) - R(\tau) = \sum_{s=1}^N z_s(\tau)$. Note that $z_1(\tau), \cdots, z_N(\tau)$ are independent. In addition, each $z_s(\tau)$ is bounded such that $|z_s(\tau)| \leq 1/\eta^2$.

Define $\bar{z}_s(\tau) \triangleq z_s(\tau) - \E[z_s(\tau) \mid X]$. Observe that $\E[\bar{z}_s(\tau) \mid X] = 0$, $|\bar{z}_s(\tau)| \leq |z_s(\tau)| + |\E[z_s(\tau) \mid X]| \leq 2/\eta^2$, and $\bar{z}_1, \cdots, \bar{z}_N$ are independent. 

In our case, we can invoke Lemma \ref{lemma:Hoeffding} by replacing $S_s = \frac{M^2}{N^2}\bar{z}_s(\tau)$ to show that
\begin{align*}
    \E[\exp(\lambda(Z(\tau)-\E[Z(\tau) \mid X]))] \leq \exp\left(\lambda^2 N \left(\frac{2M^2}{N^2\eta^2}\right)^2\right) = \exp\left(\lambda^2 \frac{4M^4}{N^3 \eta^4} \right).
\end{align*}
Therefore, conditional on $X$, $Z(\tau)-\E[Z(\tau) \mid X]$ is sub-Gaussian with variance proxy $\sigma^2 = \frac{8M^4}{N^3 \eta^4}$. \halmos

\subsection{Proof of Lemma \ref{lemma:upper_bound_event2}}\label{sec:proof_lemma:upper_bound_event2}
For any constant $\alpha_1 >0$, we apply Chebyshev's inequality to show that
    \begin{align*}
        &\Prob \left(\big|\sum_{s=1}^N \wh{\Gamma}_s'(g)-\mu(X^s, g(X^s))\big| \geq \alpha_1 N ~\bigg|~ \{X^s, D_s, Y_s\}_{s=1}^N\right)\\
        &\leq \frac{1}{\alpha_1^2 N^2} \E\left[\left(\sum_{s=1}^N \wh{\Gamma}_s'(g)-\mu(X^s, g(X^s)) \right)^2 ~\bigg|~ \{X^s, D_s, Y_s\}_{s=1}^N\right]\\
        &= \frac{1}{\alpha_1^2 N^2} \sum_{s=1}^N \E\left[\left( \wh{\Gamma}_s'(g)-\mu(X^s, g(X^s)) \right)^2 ~\bigg|~ \{X^s, D_s, Y_s\}_{s=1}^N\right].
    \end{align*} 
Next, we note that each term in the last summation can be decomposed using routine (conditional) MSE decomposition, such that
    \begin{align*}
        &\E\left[\left( \wh{\Gamma}_s'(g)-\mu(X^s, g(X^s)) \right)^2 ~\bigg|~ \{X^s, D_s, Y_s\}_{s=1}^N\right]\\
        &= \Var(\wh{\Gamma}_s'(g) \mid X^s) + \E\left[\wh{\Gamma}_s'(g)-\mu(X^s, g(X^s)) \mid \{X^s, D_s,Y_s\}_{s=1}^N\right]^2.
    \end{align*}

We bound the conditional variance and conditional bias terms as follows.
\begin{align*}
    &\Var(\wh{\Gamma}_s'(g) \mid X^s) \\
    &\leq 2\Var(\wh{\mu}(X^s, g(X^s)) \mid X^s) + 2\Var\left( \frac{\onebld\{D_s' = g(X^s)\}}{e(X^s, g(X^s))} \onebld\left\{\frac{1}{e(X^s, g(X^s))} \leq \wh{\tau}(g) \right\}(Y_s'-\wh{\mu}(X^s, g(X^s))) ~\bigg|~ X^s\right)\\
    &\leq \E\left[\frac{\onebld\{D_s' = g(X^s)\}}{e(X^s, g(X^s))^2} \onebld\left\{\frac{1}{e(X^s, g(X^s))} \leq \wh{\tau}(g) \right\} (Y_s'-\wh{\mu}(X^s, g(X^s)))^2 ~\bigg|~ X^s\right]\\
    &= M^2 \frac{1}{e(X^s,g(X^s))} \onebld\left\{\frac{1}{e(X^s, g(X^s))} \leq \wh{\tau}(g) \right\}.
\end{align*}

As for the conditional bias term, we apply results from the proof of Proposition \ref{prop:MSE_CDR} to obtain
\begin{align*}
    \E\left[\wh{\Gamma}_s'(g)-\mu(X^s, g(X^s)) \mid \{X^s, D_s,Y_s\}_{s=1}^N\right]^2 \leq M^2 \cdot \onebld\left\{\frac{1}{e(X^s, g(X^s))} > \wh{\tau}(g)\right\}.
\end{align*}

Combining the results, we have
\begin{align*}
    &\E\left[\left( \wh{\Gamma}_s'(g)-\mu(X^s, g(X^s)) \right)^2 ~\bigg|~ \{X^s, D_s, Y_s\}_{s=1}^N\right]\\
    &\leq M^2 \cdot \onebld\left\{\frac{1}{e(X^s, g(X^s))} > \wh{\tau}(g)\right\} + M^2 \frac{1}{e(X^s,g(X^s))} \onebld\left\{\frac{1}{e(X^s, g(X^s))} \leq \wh{\tau}(g) \right\}\\
    &= M^2 \left(\onebld\left\{\frac{1}{e(X^s, g(X^s))} > \wh{\tau}(g)\right\} + \frac{1}{e(X^s,g(X^s))} \onebld\left\{\frac{1}{e(X^s,g(X^s))} \leq \wh{\tau}(g) \right\}\right).
\end{align*}

Putting pieces together, we have
\begin{align*}
    &\Prob \left(\big|\sum_{s=1}^N \wh{\Gamma}_s'(g)-\mu(X^s, g(X^s))\big| \geq \alpha N ~\bigg|~ \{X^s, D_s, Y_s\}_{s=1}^N\right) \\
    &\leq \frac{M^2}{\alpha_1^2 N^2} \sum_{s=1}^N \left(\onebld\left\{\frac{1}{e(X^s, g(X^s))} > \wh{\tau}(g) \right\} + \frac{1}{e(X^s, g(X^s))} \onebld\left\{\frac{1}{e(X^s, g(X^s))} \leq \wh{\tau}(g) \right\}\right)\\
    &\leq \frac{M^2}{\alpha_1^2 N^2} \left(\left(\sum_{s=1}^N \onebld\left\{\frac{1}{e(X^s, g(X^s))} > \wh{\tau}(g) \right\}\right)^2 + \sum_{s=1}^N \frac{1}{e(X^s, g(X^s))} \onebld\left\{\frac{1}{e(X^s, g(X^s))} \leq \wh{\tau}(g) \right\}\right)\\
    &= \frac{\wt{B}(g, \wh{\tau}(g))}{\alpha_1^2}.
\end{align*}
Letting $\alpha_1 = 2\wt{B}(g, \wh{\tau}(g))^{1/2}$ to yield the desired result.

We now proceed to prove the second probability bound in Lemma \ref{lemma:upper_bound_event2}. For a fixed policy $g \in \Gcal_{\rm LPC}$ and any constant $\alpha_2>0$, we have

{\small
\begin{align}\label{prob_bound:Chebyshev1}
    \notag&\Prob\left(\sum_{s=1}^N \frac{\onebld\{D_s' = g(X^s)\}}{e(X^s,g(X^s))^2} \onebld\left\{\frac{1}{e(X^s,g(X^s))} \leq \wh{\tau}(g)\right\} \geq \alpha_2 + \sum_{s=1}^N \frac{1}{e(X^s,g(X^s))} \onebld\left\{\frac{1}{e(X^s,g(X^s))} \leq \wh{\tau}(g)\right\} \bigg| \{X^s,D_s,Y_s\}_{s=1}^N\right)\\
    \notag&\quad = \Prob\Bigg(\sum_{s=1}^N \frac{\onebld\{D_s' = g(X^s)\}}{e(X^s,g(X^s))^2} \onebld\left\{\frac{1}{e(X^s,g(X^s))} \leq \wh{\tau}(g)\right\} - \E\left[\frac{\onebld\{D_s' = g(X^s)\}}{e(X^s,g(X^s))^2} \onebld\left\{\frac{1}{e(X^s,g(X^s))} \leq \wh{\tau}(g)\right\} \bigg| \{X^s,D_s,Y_s\}_{s=1}^N\right] \\
    \notag&~\quad\quad\quad\quad \geq \alpha_2 \bigg| \{X^s,D_s,Y_s\}_{s=1}^N\Bigg)\\
    \notag&\quad \leq \frac{1}{\alpha_2^2} \sum_{s=1}^N \Var\left(\frac{\onebld\{D_s' = g(X^s)\}}{e(X^s,g(X^s))^2} \onebld\left\{\frac{1}{e(X^s,g(X^s))} \leq \wh{\tau}(g)\right\} \bigg| \{X^s, D_s,Y_s\}_{s=1}^N \right)\\
    \notag&\quad \leq \frac{1}{\alpha_2^2} \sum_{s=1}^N \E\left[\left(\frac{\onebld\{D_s' = g(X^s)\}}{e(X^s,g(X^s))^2} \onebld\left\{\frac{1}{e(X^s,g(X^s))} \leq \wh{\tau}(g)\right\}\right)^2 \bigg| \{X^s, D_s,Y_s\}_{s=1}^N \right]\\
    &\quad \leq \frac{1}{\alpha_2^2} \sum_{s=1}^N \wh{\tau}(g)^3 = \frac{N}{\alpha_2^2} \wh{\tau}(g)^3,
\end{align}}%
where the first inequality uses Chebyshev's inequality. By letting $\alpha_2 = 2\sqrt{N} \wh{\tau}(g)^{3/2}$, \eqref{prob_bound:Chebyshev1} is upper bounded by $1/4$.

Note that when condition on $\{X^s,D_s,Y_s\}_{s=1}^N$, the term $\sum_{s=1}^N \onebld\left\{\frac{1}{e(X^s, g(X^s))} > \wh{\tau}(g)\right\}$ is a constant, therefore, the probability bound in \eqref{prob_bound:Chebyshev1} implies that
\begin{align*}
    &\Prob\left(\wh{B}'(g, \wh{\tau}(g)) \geq \frac{M^2}{N^2} \cdot 2\sqrt{N} \wh{\tau}(g)^{3/2} + \wt{B}(g, \wh{\tau}(g)) \bigg| \{X^s, D_s,Y_s\}_{s=1}^N\right)\\
    &\quad \leq \Prob\left(\wh{B}'(g, \wh{\tau}(g)) \geq 2\Delta(g, \wh{\tau}(g)) + \wt{B}(g, \wh{\tau}(g)) \bigg| \{X^s, D_s,Y_s\}_{s=1}^N\right),
\end{align*}
which implies that
\begin{align*}
    \Prob\left(\wh{B}'(g, \wh{\tau}(g)) \geq 4 \max\left\{\Delta(g, \wh{\tau}(g)), \wt{B}(g, \wh{\tau}(g))\right\} \bigg| \{X^s, D_s,Y_s\}_{s=1}^N\right) \leq\frac{1}{4}.
\end{align*}
The proof of Lemma \ref{lemma:upper_bound_event2} is completed. \halmos

\section{Supporting Results}
\begin{lemma}[\cite{natarajan1989learning}]\label{lemma:Natarajan}
    Let $\cal{G} \triangleq \{g: S \mapsto \cal{A}\}$ for a finite set $S \in \cal{X}^{|S|}$. Then we have
    \begin{align*}
        |\Gcal| \leq |\cal{X}|^{\cal{N}_{dim} (\Gcal)} |\cal{A}|^{2\cal{N}_{dim}(\Gcal)}.
    \end{align*}
\end{lemma}
The next lemma gives an upper bound for the Natarajan dimension of the linear policy class $\Gcal_{\rm LPC}$ considered in the paper.

\begin{lemma}[Corollary 29.8 of \cite{shalev2014understanding}]\label{lemma:n_dim_lpc}
For a class of linear policies, defined as
\begin{equation*}
\mathcal{G}_{\mathrm{LPC}} \triangleq \left \{g: g(\xi) = \argmax_{j \in [J]}~(\xi^\top \beta^j + b_j), \xi \in \R^{p}, \beta^j \in \R^{p}, b_j \in \R \right\},
\end{equation*}
the Natarajan dimension of $\Gcal_{\rm LPC}$ satisfies
\begin{align*}
    \cal{N}_{\rm dim}(\Gcal_{\rm LPC}) \leq |J|(p+1).
\end{align*}

\end{lemma}

Hoeffding’s lemma \citep{hoeffding1963probability} implies that all random variables that are bounded uniformly are actually sub-Gaussian. The lemma below shows that the sum of independent sub-Gaussian random variables is also sub-Gaussian.
\begin{lemma}[Theorem 1.6 of \cite{rigollet2023high}]\label{lemma:Hoeffding}
Suppose $S_1, \cdots, S_N$ are independent, mean-zero, and $|S_s| \leq c$, then $\sum_{s=1}^NS_s$ is sub-Gaussian with variance proxy at most $2Nc$, i.e., for all $\lambda \in \R$,
\begin{align*}
    \E\left[\exp\left(\lambda \sum_{s=1}^N S_s\right) \big| X\right] \leq \exp\left(\frac{\lambda^2}{2} \cdot 2Nc^2\right) \leq \exp(\lambda^2 Nc^2).
\end{align*}
\end{lemma}

\newpage
\section{The Restricted Integer Program}\label{sec:restricted_ip}
Given the vector $\widebar{\boldsymbol{\beta}}$ and the index sets $\mathcal J_{1;<}^{-\delta_1}(\widebar{\boldsymbol{\beta}}), \mathcal J_{1;>}^{+\delta_1}(\widebar{\boldsymbol{\beta}}),  \mathcal J_{1;0}^{\pm\delta_1}(\widebar{\boldsymbol{\beta}}), \mathcal J_{2;<}^{-\delta_2}(\widebar{\boldsymbol{\beta}}), \mathcal J_{2;>}^{+\delta_2}(\widebar{\boldsymbol{\beta}}), \mathcal J_{2;0}^{\pm\delta_2}(\widebar{\boldsymbol{\beta}})$, the restricted integer program has an objective function $\wh{\Psi}_{\varepsilon}^{\, \boldsymbol{\delta}}(\boldsymbol{\beta})$, subject to the following constraints:
{\small
$\begin{array}{ll}
& \, h^{1}_j(X^s, \boldsymbol{\beta};\varepsilon) \geq \underline{B}(1-z^1_{j,s}), \quad \forall (j,s) \in \mathcal J_{1;0}^{\pm\delta_1}(\widebar{\boldsymbol{\beta}})\\ [0.1in]
& \, h^{1}_j(X^s, \boldsymbol{\beta};\varepsilon) \geq 0, \quad  \forall (j,s) \in \mathcal J_{1;>}^{+\delta_1}(\widebar{\boldsymbol{\beta}})\\ [0.1in] 
& \, h^{1}_j(X^s, \boldsymbol{\beta};\varepsilon) \geq \underline{B}, \quad  \forall (j,s) \in \mathcal J_{1;<}^{-\delta_1}(\widebar{\boldsymbol{\beta}})\\ [0.1in] 
& \, h^{2}_{D_s}(X^s, \boldsymbol{\beta};\varepsilon) \leq \overline{B}\, z^2_{D_s,s}, \quad  \forall  s\in \mathcal J_{2;0}^{\pm\delta_2}(\widebar{\boldsymbol{\beta}})\\ [0.1in]
& \, h^{2}_{D_s}(X^s, \boldsymbol{\beta};\varepsilon) \leq 1, \quad  \forall  s\in \mathcal J_{2;>}^{+\delta_2}(\widebar{\boldsymbol{\beta}})\\ [0.1in]
& \, h^{2}_{D_s}(X^s, \boldsymbol{\beta};\varepsilon) \leq 0, \quad  \forall  s\in \mathcal J_{2;<}^{-\delta_2}(\widebar{\boldsymbol{\beta}})\\ [0.1in]
& \,      
 \left( \begin{array}{ll} \displaystyle{\min_{s \geq t \geq m^*}}  \left( \displaystyle{\sum_{k=t}^N} \left(2(N-k)+1\right) - \sum_{k \in [t,N] \cap \mathcal J_{2;0}^{\pm\delta_2}(\widebar{\boldsymbol{\beta}})}2   (C_{k})^2 
 z_{D_{k},k}^2 - \sum_{k \in [t,N] \cap \mathcal J_{2;>}^{+\delta_2}(\widebar{\boldsymbol{\beta}})} 2   (C_{k})^2 \, \right) \\ [0.2in]
 -  \min \left\{ \displaystyle{\min_{N \geq m  \geq s+1}}   \left( \sum_{k=m}^N \left(2(N-k)+1\right) -  \sum_{k \in [m,N] \cap \mathcal J_{1;>}^{+\delta_1}(\widebar{\boldsymbol{\beta}})} 2   (C_{k})^2 z_{D_{k},k}^1 -  \sum_{k \in [m,N] \cap \mathcal J_{1;0}^{\pm\delta_1}(\widebar{\boldsymbol{\beta}})} 2   (C_{k})^2 \, \right) , 0 \right\} \end{array}  \right) \geq \underline{B} (1- w_s^1), \, \forall  N \geq s  \geq m^* \\ [0.5in]
& \,        
 \left(  \begin{array}{ll} \displaystyle{\min_{s \geq t \geq m^*}}  \left( \displaystyle{\sum_{k=t}^N} \left(2(N-k)+1\right) - \sum_{k \in [t,N] \cap \mathcal J_{1;0}^{\pm\delta_1}(\widebar{\boldsymbol{\beta}})}2   (C_{k})^2 \cdot z_{D_k,k}^1 - \sum_{k \in [t,N] \cap \mathcal J_{1;>}^{+\delta_1}(\widebar{\boldsymbol{\beta}})}2   (C_{k})^2\, \right) \\ [0.2in]
 - \min \left\{ \displaystyle{\min_{N \geq m  \geq s+1}} \left( \displaystyle{\sum_{k=m}^N} \left(2(N-k)+1\right) - \sum_{k \in [m,N] \cap \mathcal J_{2;0}^{\pm\delta_2}(\widebar{\boldsymbol{\beta}})}2 (C_{k})^2 \cdot z_{D_k,k}^2 - \sum_{k \in [m,N] \cap \mathcal J_{2;>}^{+\delta_2}(\widebar{\boldsymbol{\beta}})}2 (C_{k})^2\, \right) , 0 \right\} \end{array}  \right) \leq \overline{B}  w_s^2, \, \forall  N \geq s  \geq m^* \\ [0.5in]  
 &  \,   z_{j,s}^1 \in \{0,1\}, \quad \forall (j,s) \in \mathcal J_{1;0}^{\pm\delta_1}(\widebar{\boldsymbol{\beta}})\\[0.1in]
 & \,  z_{D_s,s}^2 \in \{0,1\}, \quad \forall  s \in \mathcal J_{2;0}^{\pm\delta_2}(\widebar{\boldsymbol{\beta}})\\[0.1in]
& \, w_{s}^1,w_{s}^2 \in \{0,1\}, \quad \forall   N \geq s  \geq m^*\\[0.1in] 
& \, w_{s}^1 \geq w_{s'}^1, \quad \forall s' \geq s,  \forall s ', s\in [N].
\end{array}$}%

\section{Additional Details of Numerical Experiments}\label{sec:omitted_exp_details}
\subsection{Implementation Details of the PIP Algorithm}\label{appx:implementation_PIP}
\noindent $\bullet$ \textbf{Initial solutions:} To obtain the initial solution $\bm{\beta}^0$, we solve the LP relaxation of the full MIP problem formulated in Section \ref{sec:computation}. This yields the initial $\bm{\beta}^0$ at which PIP is initiated.

\noindent  $\bullet$ \textbf{Stopping criterion:} Our PIP algorithm terminates when either of the following conditions is met:
(i) the number of updates to the in-between interval reaches 15, or (ii) the objective value remains unchanged for 3 consecutive iterations. In addition, each PIP subproblem is solved with a time limit: 100 seconds for small-scale problems ($N < 1,500$) and 120 seconds for large-scale problems ($N \geq 1,500$).

\subsection{Omitted Details in Section \ref{sec:computational_performance}}
\textbf{Reward function.} The reward function is specified as
\begin{equation*}
    Y(D) = Y(1)\onebld\{D=1\}+Y(2)\onebld\{D=2\}+Y(3)\onebld\{D=3\}+Y(4)\onebld\{D=4\}+ \varepsilon,
\end{equation*}
where
\begin{align*}
Y(1) &= \exp\{1.2 + 0.2X_0 + 1.7X_1 - 0.2X_2 + 2X_0X_1\} + X_r + \varepsilon,\\
Y(2) &= \exp\{1-X_0 + 2X_1 + 2X_0X_1\} + X_r + \varepsilon,\\
Y(3) &= \exp\{1.2 + 0.2X_0 + 1.7X_1 - 0.1 X_2 + 2X_0X_1 + 1.3X_0X_1\} + X_r + \varepsilon,\\
Y(4) &= \exp\{1.6 + 2X_0 -0.1X_1 + 2X_0X_1 -1.2X_1X_2\} + X_r + \varepsilon,
\end{align*}
with $r$ in $X_r$ being a random index in $\{1, \cdots, 20\}$ that may vary across the data sets, and the noise term $\varepsilon$ follows a $\mathrm{Lognormal}(0,0.001)$ distribution.

\textbf{Regularized policy optimization.} Given the high-dimensional data setting, we employ a regularizer to induce sparsity in parameter selection. Specifically, the optimization problem in Equation \eqref{eq:CDR_policy_approx_optimization} is adapted to
\begin{equation}\label{eq:regularized_policy_optimization}
\max_{\boldsymbol{\beta} \in \mathbb R^p} 
\quad 
\wh{\Psi}_\varepsilon(\boldsymbol{\beta}) - \lambda \sum_{j=1}^J R_j(\beta^{j}),
\end{equation}
where $R_j(\bullet)$ is the regularizer, e.g., $R_j(\beta^{j}) = \|\beta^{j}\|_1$. Note that problem \eqref{eq:regularized_policy_optimization} fits the form of HSCOP in Definition \ref{def:HSCOP} by setting the Bouligand differentiable function as $f_1(\boldsymbol{\beta}) = -\lambda\sum_{j=1}^J R_j(\beta^j)$. We note that many commonly used regularizers are not differentiable; see e.g., \cite{ahn2017difference} and \cite{cui2018composite}.

\subsection{Visualization of the Learned Policy}

\end{document}